\documentclass[12pt]{article}%
\usepackage{amssymb}
\usepackage{graphicx}
\usepackage{amsmath}%
\usepackage{amsfonts}
\newtheorem{theorem}{Theorem}

\newtheorem{definition}[theorem]{Definition}

\newtheorem{lemma}[theorem]{Lemma}

\newtheorem{proposition}[theorem]{Proposition}

\newenvironment{proof}[1][Proof]{\textbf{#1.} }{\ \rule{0.5em}{0.5em}}
\begin{document}

\title{Poisson Hypothesis for Information Networks \\(A study in non-linear Markov processes)\\I. Domain of Validity }
\author{Alexander Rybko\\Institute for the Information Transmission Problems,\\Russian Academy of Sciences, Moscow, Russia,\\rybko@iitp.ru
\and Senya Shlosman\\Centre de Physique Theorique, CNRS, \\Luminy, 13288 Marseille, France\\shlosman@cpt.univ-mrs.fr}
\maketitle

\begin{abstract}
In this paper we study the Poisson Hypothesis, which is a device to analyze
approximately the behavior of large queueing networks. We prove it in some
simple limiting cases. We show in particular that the corresponding dynamical
system, defined by the non-linear Markov process, has a line of fixed points
which are global attractors. To do this we derive the corresponding non-linear
equation and we explore its self-averaging properties. We also argue that in
cases of havy-tail service times the PH can be violated.

MSC-class: 82C20 (Primary), 60J25 (Secondary)

\end{abstract}

\section{Introduction}

The Poisson Hypothesis deals with large queueing systems. For general systems
one can not compute exactly the quantities of interest, so various
approximations are used in practice. The Poisson Hypothesis was formulated
first by L. Kleinrock in \cite{K}. It is the statement that certain
approximation becomes exact in the appropriate limit. It concerns the
following situation. Suppose we have a large network of servers, through which
many customers are travelling, being served at different nodes of the network.
If the node is busy, the customers wait in the queue. Customers are entering
into the systems via some nods, and the external flows of customers from the
outside are Poissonian. The service time at each node is random, with some
fixed distribution, depending on the node. We are interested in the stationary
distribution $\pi_{\mathcal{N}}$ at a given node $\mathcal{N} $: what is the
distribution of the queue at it, what is the average waiting time, etc. Except
for a very few special cases, when the service times are exponential, the
distributions $\pi_{\mathcal{N}}$ in general can not be computed. The recipe
of the Poisson Hypothesis for approximate computation of $\pi_{\mathcal{N}}$
is the following:

\begin{itemize}
\item consider the total flow $\mathcal{F}$ of customers to the node
$\mathcal{N}.$ (In general, $\mathcal{F}$ is not Poissonian, of course.)
Replace $\mathcal{F}$ with a constant rate Poisson flow $\mathcal{P},$ the
rate being equal to the average rate of $\mathcal{F}.$ Compute the stationary
distribution $\hat{\pi}_{\mathcal{N}}$ at $\mathcal{N},$ corresponding to the
inflow $\mathcal{P}.$ (These computations are the subject of classical
queueing theory and usually provide explicit formulas.) The claim is that
$\hat{\pi}_{\mathcal{N}}\approx\pi_{\mathcal{N}}.$
\end{itemize}

The Poisson Hypothesis is supposed to give a good estimate if the internal
flow to every node $\mathcal{N}$ is a sum of flows from many other nodes, and
each of these flows constitute only a small fraction of the total flow to
$\mathcal{N}.$

Clearly, the Poisson Hypothesis can not be literally true. It can hopefully
hold only after some kind of ``thermodynamic''\ limit is taken. Its meaning is
that in the long run the different nodes become virtually independent, i.e.
propagation of chaos takes place. The reason for that should be that any
synchronization of the nodes, if initially present, dissolves with time, due
to the randomness of the service times.

In the present paper we prove the Poisson Hypothesis for the information
networks in some simple cases. Namely, we will consider the following closed
queueing network. Let there be $M$ servers and $N$ customers to be served. The
distribution of the service time is given by some fixed random variable
$\eta.$ Upon being served, the customer chooses one of $M$ servers with
probability $\frac{1}{M},$ and goes for the service there. If there is a
queue, he waits for his turn. Then in the limit $M,N\rightarrow\infty,$ with
$\frac{N}{M}\rightarrow\rho,$ the Poisson Hypothesis holds, under certain
general restrictions on $\eta$. More precisely, \textbf{Poisson Hypothesis}
\textit{for our model} states that:

\begin{enumerate}
\item every $k$-tuple of servers becomes asymptotically mutually independent
and identically distributed, as $M,N\rightarrow\infty,$ for any $k,$

\item the total flow $\mathcal{F}_{M,N}$ of customers to a given node goes, as
$M,N\rightarrow\infty,$ to a Poisson flow $\mathcal{P},$

\item the rate function $\lambda\left(  t\right)  $ of the Poisson flow
$\mathcal{P}$ goes to a constant limit $c\left(  \rho\right)  $, as
$t\rightarrow\infty,$\ which depends only on the load $\rho$\ and thus can be
(easily) computed apriori.

{\small To be precise, one needs also to put conditions on convergence of the
sequence of initial states }$\nu_{M,N}$ {\small of our servers, in order the
above to hold, see Section \ref{KR1} for more details.}
\end{enumerate}

In this paper we prove the Poisson Hypothesis (PH) in the following cases:

\begin{itemize}
\item for a special class of the service times $\eta$ (with some exponential
moments finite; actually a bit more is needed, see $\left(  \ref{11}\right)
$) -- with no extra conditions;

\item for general $\eta$-s with heavy tails -- provided the initial state of
our system possesses certain desynchronization property (valid once the
expected service time is finite, see $\left(  \ref{017}\right)  $).
\end{itemize}

In a subsequent paper \cite{RS} we will show that for $\eta$-s with only
polynomial moments and for certain initial states the initial synchronization
of the nodes may not vanish with time. So PH\ can be violated, we can not
predict the long time behavior of a single server, and we have the phenomenon
of phase transition, which manifests itself in the strong dependence on the
initial state of the system.

An important step in establishing the validity of PH was made in the paper
\cite{KR1}. Namely, the properties 1 and 2 were obtained there. However, the
technique of \cite{KR1} was not enough to prove the relaxation property
$\lambda\left(  t\right)  \rightarrow c,$ and moreover it does not hold in
general. It was proven there that the situation at a given single server is
described by the so-called non-linear Markov process $\mu_{t}$ with Poissonian
input with rate $\lambda\left(  t\right)  ,$ and the (non-Poissonian) output
with the same rate $\lambda\left(  t\right)  $. The remaining problem can be
formulated as follows: this non-linear Markov process defines some complicated
dynamical system, and the question is about its invariant measures. Namely,
this system has one parameter family of fixed points, and one needs to show
that every trajectory converges to one of them.

In the present paper we complete the program, showing that the above
relaxation $\lambda\left(  t\right)  \rightarrow c$ indeed takes place for
certain class of the service times $\eta$ and certain class of initial states,
and so $\mu_{t}\rightarrow\mu_{c},$ where $\mu_{c}$ is the stationary
distribution of the stationary Markov process with the Poisson input,
corresponding to constant rate $\lambda\left(  t\right)  =c$. In the language
of dynamical systems, we find conditions under which there are no other
invariant measures except these defined by the fixed points.

The central discovery of the present paper, which seems to be the key to the
solution of the problem, is that, roughly speaking, the function
$\lambda\left(  t\right)  $ has to satisfy the following non-linear equation:
\begin{equation}
\lambda\left(  t\right)  =\left[  \lambda\left(  \cdot\right)  \ast
q_{\lambda,t}\left(  \cdot\right)  \right]  \left(  t\right)  .\label{200}%
\end{equation}
Here $\ast$ stays for convolution: for two functions $a\left(  \cdot\right)
,b\left(  \cdot\right)  $ it is defined as
\[
\left[  a\left(  \cdot\right)  \ast b\left(  \cdot\right)  \right]  \left(
t\right)  =\int a\left(  t-x\right)  b\left(  x\right)  \,dx,
\]
while $q_{\lambda,t}\left(  \cdot\right)  $ is a family of probability
densities with $t$ real, which depends also in an implicit way on the unknown
function $\lambda\left(  \cdot\right)  $. We call $\left(  \ref{200}\right)  $
the self-averaging property. The present paper consists therefore of two
parts: we prove that indeed the self-averaging relation holds, and we prove
then that in certain cases it implies relaxation. The self-averaging property
is a special case of a more general averaging relation $\left(  \ref{007}%
\right)  ,$\ which relates the output flow rate to the input flow rate in a
queuing process with one server (called $M(t)/GI/1$\ in queueing theory
jargon). It appears to be new.

It is amazing that the relation $\left(  \ref{200}\right)  $ depends crucially
on the validity of some purely combinatorial statement concerning certain
problem of the placement of the rods on the line $\mathbb{R}^{1},$ see Section
6. This validity seems not at all obvious or to be expected, so to have it is
our luck.\medskip

To fix the terminology, we remind the reader here what we mean by the
\textbf{non-linear Markov process} (see \cite{M1}, \cite{M2}). We do this for
the simplest case of discrete time Markov chains, taking values in a finite
set $S,$ $\left\vert S\right\vert =k.$ In such a case the set of states of
this Markov chain is a simplex $\Delta_{k}$ of all probability measures on
$S,$ $\Delta_{k}=\left\{  \mu=\left(  p_{1},...,p_{k}\right)  :p_{i}%
\geq0,p_{1}+...+p_{k}=1\right\}  ,$ while the Markov evolution defines a map
$P:\Delta_{k}\rightarrow\Delta_{k}.$ In the case of usual Markov chain $P$ is
affine, and this is why we will call it \textit{linear}\textbf{\ }chain. In
this case the matrix of transition probabilities coincides with $P.$ The
non-linear Markov chain is defined by a family of transition probability
matrices $P_{\mu},$ $\mu\in\Delta_{k},$ so that matrix element $P_{\mu}\left(
i,j\right)  $ is a probability of going from $i$ to $j $ in one step, starting
\textit{in the state}\textbf{\ }$\mu.$ The \textit{(non-linear)} map $P$ is
then defined by $P\left(  \mu\right)  =\mu P_{\mu}.$

The ergodic properties of the linear Markov chains are settled by the
Perron-Frobenius theorem. In particular, if the linear map $P$ is such that
the image $P\left(  \Delta_{k}\right)  $ belongs to the interior
$\mathrm{Int\,}\left(  \Delta_{k}\right)  $ of $\Delta_{k},$ then there is
precisely one point $\mu\in\mathrm{Int\,}\left(  \Delta_{k}\right)  ,$ such
that $P\left(  \mu\right)  =\mu,$ and for every $\nu\in\mathrm{\,}\Delta_{k}$
we have the convergence $P^{n}\left(  \nu\right)  \rightarrow\mu$ as
$n\rightarrow\infty.$

In case $P$ is non-linear, we are dealing with more or less arbitrary
dynamical system on $\Delta_{k}$, and the question about stationary states of
the chain or about measures on $\Delta_{k}$ invariant under $P$ can not be
settled in general.\medskip

Therefore it is natural to ask about the specific features of our dynamical
system, which permit us to find all its invariant measures. We explain this in
the following subsection. The reader who is not interested in this aspect of
the problem can safely skip it.\medskip

\textbf{Dynamical systems aspect}{\small . Here we will use the notation of
the paper, though in fact the situation of the paper is more complicated; in
particular the underlying space is not a manifold, but a space of all
probability measures over some non-compact set.}

{\small Let }$M${\small \ be a manifold, supplied with the following
structures:}

\begin{itemize}
\item {\small for every point }$\mu\in M${\small \ and every }$\lambda
>0${\small \ a tangent vector }$X\left(  \mu,\lambda\right)  ${\small \ at
}$\mu${\small \ is defined,}

\item {\small a function }$b:M\rightarrow R^{+}${\small \ is fixed.}
\end{itemize}

{\small We want to study the dynamical system }
\begin{equation}
\frac{d}{dt}\mu\left(  t\right)  =X\left(  \mu\left(  t\right)  ,b\left(
\mu\left(  t\right)  \right)  \right)  .\label{001}%
\end{equation}
{\small Its flow conserves another given function, }$N:M\rightarrow R^{+}%
,${\small \ and we want to prove that our dynamical system has one-parameter
family of fixed points - each corresponding to one value of }$N${\small \ -
and no other invariant measures.}

{\small We have the following extra properties of our dynamical system:}

{\small Let }$\lambda\left(  t\right)  >0;${\small \ consider the differential
equation }
\begin{equation}
\frac{d}{dt}\mu\left(  t\right)  =X\left(  \mu\left(  t\right)  ,\lambda
\left(  t\right)  \right)  ,\;t\geq0,\label{002}%
\end{equation}
{\small with }$\mu\left(  0\right)  =\nu.${\small \ We denote the solution to
it by }$\mu_{\nu,\lambda\left(  \cdot\right)  }\left(  t\right)
.${\small \ We know that}

{\small for every }$c>0${\small \ and every initial data }$\nu,${\small \ the
solution }$\mu_{\nu,\lambda\left(  \cdot\right)  }\left(  t\right)
${\small \ to }$\left(  \ref{002}\right)  ${\small \ converges to some
stationary point }$\nu_{c}\in M,${\small \ }
\begin{equation}
\mu_{\nu,\lambda\left(  \cdot\right)  }\left(  t\right)  \rightarrow\nu
_{c},\text{{\small provided} }\lambda\left(  t\right)  \rightarrow
c\text{{\small as} }t\rightarrow\infty,\label{004}%
\end{equation}

\begin{itemize}
\item {\small for the function }$N${\small \ we have }
\[
\frac{d}{dt}N\left(  \mu_{\nu,\lambda\left(  \cdot\right)  }\left(  t\right)
\right)  =\lambda\left(  t\right)  -b\left(  \mu_{\nu,\lambda\left(
\cdot\right)  }\left(  t\right)  \right)  .
\]
{\small In particular, for every trajectory }$\hat{\mu}_{\nu}\left(  t\right)
${\small \ of }$\left(  \ref{001}\right)  ${\small \ (where }$\hat{\mu}_{\nu
}\left(  0\right)  =\nu${\small ) we have }$N\left(  \hat{\mu}_{\nu}\left(
t\right)  \right)  =N\left(  \nu\right)  .${\small \ Also, }$N\left(  \nu
_{c}\right)  ${\small \ is continuous and increasing in }$c${\small ;}

\item {\small for every }$\nu,\lambda\left(  \cdot\right)  ${\small \ and
every }$t>0${\small \ there exists a probability density }$q_{\nu,\lambda
,t}\left(  x\right)  ,\;x\geq0,${\small \ such that }
\[
b\left(  \mu_{\nu,\lambda\left(  \cdot\right)  }\left(  t\right)  \right)
=\left(  \lambda\ast q_{\nu,\lambda,t}\right)  \left(  t\right)  ,
\]
{\small where }
\[
\left(  \lambda\ast q_{\nu,\lambda,t}\right)  \left(  y\right)  =\int_{x\geq
0}q_{\nu,\lambda,t}\left(  x\right)  \lambda\left(  y-x\right)  \,dx.
\]
{\small The family }$q_{\nu,\lambda,t}\left(  x\right)  ${\small \ satisfies:
}
\[
\int_{0}^{1}q_{\nu,\lambda,t}\left(  x\right)  \,dx=1\text{ for all }%
\nu,\lambda,t,
\]
{\small and }
\[
\inf_{\substack{\nu,\lambda,t \\x\in\left[  0,1\right]  }}q_{\nu,\lambda
,t}\left(  x\right)  >0
\]
{\small (absolute continuity with respect to Lebesgue).}
\end{itemize}

{\small Then for every initial state }$\nu${\small \ }
\begin{equation}
\hat{\mu}_{\nu}\left(  t\right)  \rightarrow\nu_{c},\label{003}%
\end{equation}
{\small where }$c${\small \ satisfies }$N\left(  \nu_{c}\right)  =N\left(
\nu\right)  .$

{\small Our statement follows from the fact that the self-averaging property,
}
\[
f\left(  t\right)  =\left(  f\ast q_{t}\right)  \left(  t\right)  ,
\]
{\small with }$q_{t}\left(  \cdot\right)  ${\small \ being a family of
probability densities on }$\left[  0,1\right]  ${\small , implies that
}$f\left(  t\right)  \rightarrow const${\small \ as }$t\rightarrow\infty
,${\small \ so }$\left(  \ref{003}\right)  ${\small \ follows from }$\left(
\ref{004}\right)  .${\small \ This implication is the subject of the Theorems
\ref{finite range}, \ref{infinite range}, \ref{noisy case}.}

\medskip

We feel that the relation $\left(  \ref{200}\right)  $ is an important feature
of the subject we are interested in. Therefore in the present paper we study
it and the related questions in some generality.

$i)$ We start with the equation
\begin{equation}
f\left(  t\right)  =\left[  f\left(  \cdot\right)  \ast q_{t}\left(
\cdot\right)  \right]  \left(  t\right)  .\label{201}%
\end{equation}
Here we suppose that $q_{t}\left(  \cdot\right)  $ is just some one-parameter
family of probability densities (independent of $f$ ), so $\left(
\ref{200}\right)  $ is a special case of $\left(  \ref{201}\right)  .$ On the
other hand, we suppose additionally that all the distributions $q_{t}\left(
\cdot\right)  $ are supported by some finite interval. We establish relaxation
in this case.

$ii)$ We then do the same for the case of distributions $q_{t}\left(
\cdot\right)  $ with unbounded support.

$iii)$ Last, we treat the true problem, where in addition to the infinite
support, an extra parameter $\mu$ appears and an extra perturbation is added
to convolution term in $\left(  \ref{201}\right)  :$%
\begin{equation}
\lambda\left(  t\right)  =\left(  1-\varepsilon_{\lambda,\mu}\left(  t\right)
\right)  \left[  \lambda\left(  \cdot\right)  \ast q_{\lambda,\mu,t}\left(
\cdot\right)  \right]  \left(  t\right)  +\varepsilon_{\lambda,\mu}\left(
t\right)  Q_{\lambda,\mu}\left(  t\right)  .\label{202}%
\end{equation}
Here the parameter $\varepsilon_{\lambda,\mu}\left(  t\right)  $ is small:
$\varepsilon_{\lambda,\mu}\left(  t\right)  \rightarrow0$ as $t\rightarrow
\infty,$ the term $Q_{\lambda,\mu}\left(  t\right)  $ is uniformly bounded,
while the meaning of $\mu$ will be explained later.

As we proceed from $i)$ to $iii),$ we will have to assume more about the class
of distributions $\left\{  q_{\cdot}\right\}  ,$ for which the self-averaging
implies relaxation.

\medskip

We finish this introduction by a brief discussion of the previous work on the
subject, and their methods.

As we said before, part of the proof of the Poissonian Hypothesis -- the so
called Weak Poissonian Hypothesis -- was obtained in \cite{KR1}. By proving
that the Markov semigroups describing the Markov processes for finite $M,N,$
after factorization by the symmetry group of the model converge, as
$M,N\rightarrow\infty,$ $\frac{N}{M}\rightarrow\rho$, to the semigroup,
describing the non-linear Markov process, the authors have proven that the
limit flows to each node are independent Poisson flows with the same rate
function $\lambda\left(  t\right)  .$ This statement is fairly general, and
can be generalized to other models with the same kind of the symmetry -- the
so-called mean-field models. The general theory -- see, for example, \cite{L}
-- implies then, that all the limit points of the stationary measures of the
Markov processes with finite $M,N$ are invariant measures of the limiting
non-linear Markov process. The remaining step -- the proof that the limiting
dynamical system has no other attractors except the one-parameter family of
the fixed points -- is done in the present paper \textit{for some class of the
service times} $\eta$. Apriori this fact is not at all clear, and one can
construct natural examples of the systems with many complicated attractors,
which are reflected in the complex behavior of the Markov processes with
finite $M,N.$ However, the self-averaging property, explained above, rules out
such a possibility. It seems that the self-averaging property can also be
generalized to other mean-field models. In a forthcoming paper \cite{RS} we
show that for more general service times the corresponding dynamical system
has more complicated attractors, so that the relaxation $\lambda\left(
t\right)  \rightarrow c$\ might or might not hold, depending on the initial state.

The Poisson Hypothesis was fully established in a pioneer paper \cite{St} for
a special case when the service time is non-random. This is a much simpler
case, and the methods of the paper can not be extended to our situation. They
are sufficient for a simpler case of the Poissonian service times, which case
was studied in \cite{KR2}.

The paper \cite{DKV} deals with another mean-field model, describing some open
queueing network. Though the Poisson Hypothesis does not hold for it, the
spirit of the main statement there is the same as in the present paper: the
limiting dynamical system has precisely one global attractor, which
corresponds to the fixed point.

One of specific feature of the method of the paper \cite{DKV}, as well as
related paper \cite{DF}, is that the Markov processes have countable sets of
values. They also correspond to open networks, when the customers are coming
into the system from the outside and leave it after being served. So one can
in principle use monotonicity arguments and stochastic domination. In our
situation the phase space is (one-dimensional) real manifold, while the number
of customers is fixed, and this technique does not seem to be applicable.

\medskip

The importance of the Poisson Hypothesis as the central problem of the theory
of large queueing systems was emphasized, among others, by Roland Dobrushin
\cite{D2} and Alexander Borovkov \cite{B}.

\medskip

The organization of the paper is the following:

In the next Section \ref{notation} we introduce notation used in the rest of
the paper. We formulate the properties needed of the distribution of the
service time $\eta$ -- the main parameter of our model.

In the Section \ref{KR1} we recall more results of the paper \cite{KR1}.

In the next Section \ref{main} we formulate out main result.

In Section \ref{self-av} we derive the self-averaging relation $\left(
\ref{200}\right)  $, our key tool in establishing the Poisson Hypothesis. This
is a statement about single server queue, $M\left(  t\right)  /GI/1,$ in
queueing theory jargon. This and the next sections are self-sufficient, and
there we do not need any condition at all on the service time distribution.

The Section \ref{comb} contains the proof of a combinatorial statement dealing
with the rod placements on $\mathbb{R}^{1}.$ It is a key statement used in the
previous section. It seems also to be of independent interest.

In the next Section \ref{ten} we have collected various technical statements
and estimates used in the sequel. Again, some of them, -- for example, the
calculus Lemma \ref{calcul} -- seem to be of independent interest.

The Section \ref{gencase} contains the derivation of the noisy version of the
self-averaging relation $\left(  \ref{200}\right)  $ -- the relation $\left(
\ref{134}\right)  .$ Again, it uses the combinatorial statement of the Section
\ref{comb}.

In Section \ref{warm-up} we consider a simplified case of our theory, when
there are infinitely many servers at each node, so there are no queues. This
case is presented for pedagogical reasons only. It can be analyzed by
application of renewal theory and proving the Local Limit Theorem for i.i.d.
random variables.

In the next Section \ref{no-go} we explain that in more realistic situation
one needs the Local Limit Theorem for Markov chain, instead of i.i.d.-s, but
we show that it does not hold in general, so the analog of the renewal theory
needed there does not exist as well. This is why the rest of the paper uses
the analytic methods quite a bit.

The next three Sections \ref{frc}, \ref{infrange} and \ref{nc} contain the
main step of our proof: the derivation of the relaxation property from the
self-averaging relation. The Section \ref{frc} deals with the case of finite
range averaging kernels, the Section \ref{infrange} -- with the infinite range
kernels, while the Section \ref{nc} -- with the noisy case. As we proceed from
easier cases to more difficult, the generality of our theorems becomes less
and less.

The last Section \ref{con} contains some conclusions and lists possible
directions of further research.

The reader who wants to go directly to the proof of the general result, can
jump from Section \ref{gencase} straight to Section \ref{appp} and then to
Section \ref{nc}.

\section{Notation \label{notation}}

In this section we will fix the notation for the non-linear Markov process,
which describes a given server in the above described limit.

\textbf{Server.} It is defined by specifying the distribution of the random
serving time $\eta$, i.e. by the function
\[
F\left(  t\right)  =\mathbf{\Pr}\left\{  \text{serving time }\eta\leq
t\right\}  .
\]

We suppose that $\eta$ is such that:

\begin{enumerate}
\item the density function $p\left(  t\right)  $ of $\eta$ is positive on
$t\geq0$ and uniformly bounded from above;

\item $p\left(  t\right)  $ satisfies the following strong Lipschitz
condition: for some $C<\infty$ and for all $t\geq0$
\begin{equation}
\left\vert p\left(  t+\Delta t\right)  -p\left(  t\right)  \right\vert \leq
Cp\left(  t\right)  \left\vert \Delta t\right\vert ,\label{02}%
\end{equation}
provided $t+\Delta t>0$ and $\left\vert \Delta t\right\vert <1;$

\item introducing the random variables
\[
\eta\Bigm|_{\tau}=\left(  \eta-\tau\Bigm|\eta>\tau\right)  ,\tau\geq0,
\]
we suppose that for some $\delta>0,$ $M_{\delta,\tau}<\infty$
\begin{equation}
\mathbb{E}\left(  \eta\Bigm|_{\tau}\right)  ^{2+\delta}<M_{\delta,\tau
}.\label{154}%
\end{equation}
Of course, this condition holds once
\begin{equation}
M_{\delta}\equiv\mathbb{E}\left(  \eta\right)  ^{2+\delta}<\infty.\label{009}%
\end{equation}
In what follows, the function
\begin{equation}
R_{\eta}\left(  \tau\right)  \equiv\mathbb{E}\left(  \eta\Bigm|_{\tau}\right)
<\infty\label{181}%
\end{equation}
will play a crucial role. In particular, if
\begin{equation}
R_{\eta}\left(  \tau\right)  <\bar{C}\label{11}%
\end{equation}
for all $\tau\geq0,$ then PH holds for every initial state, as we will explain
later. However, the condition $\left(  \ref{11}\right)  $ is too restrictive;
it implies that the random variable $\eta$ has some exponential moments finite
(though the opposite is not true). \textit{Generally\ we are not assuming
}$\left(  \ref{11}\right)  $\textit{.}

\item Without loss of generality we can suppose that
\begin{equation}
\mathbb{E}\left(  \eta\right)  =1.\label{153}%
\end{equation}
In what follows, the function $p\left(  \cdot\right)  $ will be fixed.

The remaining two conditions will not be used explicitly in the present paper.
However we have to impose them since they are used in the paper \cite{KR1},
while we are using its results:

\item the probability density $p\left(  t\right)  $ is differentiable in $t,$
with $p^{\prime}\left(  t\right)  $ continuous. Moreover, introducing the
functions $p_{\tau}\left(  t\right)  $ as the densities of the random
variables $\eta\Bigm|_{\tau},$ we require that the function $p_{\tau}\left(
0\right)  $ is bounded uniformly in $\tau\geq0,$
\begin{equation}
p_{\tau}\left(  0\right)  \leq U\left(  \eta\right)  ,\label{158}%
\end{equation}
while the function $\frac{d}{d\tau}p_{\tau}\left(  0\right)  $ is continuous
and bounded uniformly in $\tau\geq0;$

\item the limits $\lim_{\tau\rightarrow\infty}p_{\tau}\left(  0\right)  ,$
$\lim_{\tau\rightarrow\infty}\frac{d}{d\tau}p_{\tau}\left(  0\right)  $ exist
and are finite.\medskip
\end{enumerate}

\textbf{Configurations.}\textit{\ }By a configuration of a server at a given
time moment $t$ we mean the following data:

\begin{itemize}
\item The number $n\geq0$ of customers waiting to be served. The customer who
is served at $t$, is included in the total amount $n.$ This quantity $n$ will
be called \textit{the length of the queue}.

\item The duration $\tau$ of the elapsed service time of the customer under
the service at the moment $t.$
\end{itemize}

Therefore the set of all configurations $\Omega$ is the set of all pairs
$\left(  n,\tau\right)  ,$ with an integer $n>0$ and a real $\tau>0,$ plus the
point $\mathbf{0}$, describing the situation of the server being idle. For a
configuration $\omega=\left(  n,\tau\right)  \in\Omega$ we define $N\left(
\omega\right)  =n.$ We put $N\left(  \mathbf{0}\right)  =0.$

\textbf{States.}\textit{\ }By a state of the system we mean a probability
measure $\mu$ on $\Omega.$ We denote by $\mathcal{M}\left(  \Omega\right)  $
the set of all states on $\Omega.$

\textbf{Observables.}\textit{\ }There are some natural random variables
associated with our system. One is the queue length in the state $\mu$,
$N_{\mu}=N_{\mu}\left(  \omega\right)  .$ We denote by $N\left(  \mu\right)  $
the mean queue length in the state $\mu:$%
\[
N\left(  \mu\right)  =\mathbb{E}\left(  N_{\mu}\right)  \equiv\left\langle
N_{\mu}\left(  \omega\right)  \right\rangle _{\mu},
\]
and we introduce the subsets $\mathcal{M}_{q}\left(  \Omega\right)
\subset\mathcal{M}\left(  \Omega\right)  ,\,q\geq0$ by
\[
\mathcal{M}_{q}\left(  \Omega\right)  =\left\{  \mu\in\mathcal{M}\left(
\Omega\right)  :N\left(  \mu\right)  =q\right\}  .
\]
Another one is the expected service time $S_{\mu},$ corresponding to the
function
\[
S\left(  \omega\right)  =\left\{
\begin{array}
[c]{ll}%
0 & \text{ for }\omega=\mathbf{0,}\\
\left(  n-1\right)  +R_{\eta}\left(  \tau\right)  & \text{for }\omega=\left(
n,\tau\right)  ,\text{with }n>0.
\end{array}
\right.
\]
Again, we define
\begin{equation}
S\left(  \mu\right)  =\mathbb{E}\left(  S_{\mu}\right)  \equiv\left\langle
S\left(  \omega\right)  \right\rangle _{\mu}.\label{017}%
\end{equation}
Clearly, if the condition $\left(  \ref{11}\right)  $ holds, then
\begin{equation}
S\left(  \mu\right)  \leq\bar{C}N\left(  \mu\right)  ,\label{12}%
\end{equation}
and therefore $S\left(  \mu\right)  $ is finite once $\mu\in\mathcal{M}%
_{q}\left(  \Omega\right)  $ for some $q.$ In general, the expected service
time $S\left(  \mu\right)  $ can be infinite for some states $\mu.$ These are
the states for which PH can be violated, as we explain in \cite{RS}.

\textbf{Input flow.}\textit{\ }Let a function $\lambda\left(  t\right)  \geq0
$ is given. We suppose that the input flow to our server is a Poisson process
with rate function $\lambda\left(  t\right)  ,$ which means in particular that
the probabilities $P_{k}\left(  t,s\right)  $ of the events that $k$ new
customers arrive during the time interval $\left[  t,s\right]  $ satisfy
\[
P_{k}\left(  t,t+\Delta t\right)  =\left\{
\begin{array}
[c]{ll}%
\lambda\left(  t\right)  \Delta t+o\left(  \Delta t\right)  & \text{ for
}k=1,\\
1-\lambda\left(  t\right)  \Delta t+o\left(  \Delta t\right)  & \text{ for
}k=0,\\
o\left(  \Delta t\right)  & \text{ for }k>1,
\end{array}
\right.
\]
as $\Delta t\rightarrow0,$ while for non-intersecting time segments $\left[
t_{1},s_{1}\right]  ,\,\,\left[  t_{2},s_{2}\right]  $ the flows are independent.

\textbf{Output flow.}\textit{\ }Suppose the initial state $\nu=\mu\left(
0\right)  ,$ as well as the rate function $\lambda\left(  t\right)  ,$ with
$\lambda\left(  t\right)  =0$ for $t<0,$ of the input flow are given. Then the
system evolves in time, and its state at the moment $t$ is given by the
measure
\[
\mu\left(  t\right)  =\mu_{\nu,\lambda\left(  \cdot\right)  }\left(  t\right)
.
\]
In particular, the probabilities $Q_{k}\left(  t,s\right)  =Q_{k}\left(
t,s;\nu,\lambda\left(  \cdot\right)  ,p\left(  \cdot\right)  \right)  $ of the
events that $k$ customers have finished their service during the time interval
$\left[  t,s\right]  $ are defined. We suppose that the customer, once served,
leaves the system.

The resulting random point process $Q_{\cdot}\left(  \cdot,\cdot\right)  $
need not, of course, be Poissonian. However we still can define its rate
function $b\left(  t\right)  $ as the one satisfying
\[
Q_{k}\left(  t,t+\Delta t\right)  =\left\{
\begin{array}
[c]{ll}%
b\left(  t\right)  \Delta t+o\left(  \Delta t\right)  & \text{ for }k=1,\\
1-b\left(  t\right)  \Delta t+o\left(  \Delta t\right)  & \text{ for }k=0,\\
o\left(  \Delta t\right)  & \text{ for }k>1,
\end{array}
\right.
\]
as $\Delta t\rightarrow0.$ The rate function $b\left(  \cdot\right)  $ of the
output flow is determined once the initial state $\nu=\mu\left(  0\right)  $
and the rate function $\lambda\left(  \cdot\right)  $ of the input flow are
given. Therefore the following (non-linear) operator $A$ is well defined:
\[
b\left(  \cdot\right)  =A\left(  \nu,\lambda\left(  \cdot\right)  \right)  .
\]
We will call the general situation, described by the pair $\nu,\lambda\left(
\cdot\right)  ,$ and $b\left(  \cdot\right)  =A\left(  \nu,\lambda\left(
\cdot\right)  \right)  $, as a General Flow Process (GFP). (This time
inhomogeneous (linear) Markov process is usually labelled in queueing theory
by $M(t)/GI/1$.)

The following is known about the operator $A,$ see \cite{KR1}:

\begin{itemize}
\item For every initial state $\nu$ the equation
\[
A\left(  \nu,\lambda\left(  \cdot\right)  \right)  =\lambda\left(
\cdot\right)
\]
-- in words:\textit{\ the rate of the input equals the rate of the output} --
has exactly one solution $\lambda\left(  \cdot\right)  =\lambda_{\nu}\left(
\cdot\right)  .$ Then the evolving state $\mu_{\nu,\lambda_{\nu}\left(
\cdot\right)  }\left(  t\right)  $ is \textit{the non-linear Markov process,
}which we will abbreviate as NMP.

\item This non-linear Markov process has the following conservation property:
for all $t$%
\[
N\left(  \mu_{\nu,\lambda_{\nu}\left(  \cdot\right)  }\left(  t\right)
\right)  =N\left(  \nu\right)
\]
(because \textquotedblleft the rates of the input flow and the output flow
coincide\textquotedblright). So the spaces $\mathcal{M}_{q}\left(
\Omega\right)  $ are invariant under non-linear Markov evolutions.

\item All the functions $\lambda_{\nu}\left(  \cdot\right)  $ are bounded:
\begin{equation}
\lambda_{\nu}\left(  t\right)  \leq C=C\left(  \eta\right) \label{13}%
\end{equation}
uniformly in $\nu$ and $t.$ (This is clear, since the output flow has its rate
uniformly bounded, by the constant $U\left(  \eta\right)  ,$ see $\left(
\ref{158}\right)  .$ So $\left(  \ref{13}\right)  $ holds with $C\left(
\eta\right)  =U\left(  \eta\right)  .$)

\item For every constant $c\in\lbrack0,1)$ there exists the initial state
$\nu_{c},$ such that
\begin{equation}
A\left(  \nu_{c},c\right)  =c.\label{03}%
\end{equation}
(Here we identify the constant $c$ with the function taking the value $c$
everywhere.) Moreover, this measure $\nu_{c}$ is a stationary state: $\mu
_{\nu_{c},c}\left(  t\right)  =\nu_{c}$ for all $t>0.$ The function
$c\rightsquigarrow N\left(  \nu_{c}\right)  $ is continuous increasing, with
$N\left(  \nu_{0}\right)  =0,$ $N\left(  \nu_{c}\right)  \uparrow\infty$ as
$c\rightarrow1.$
\end{itemize}

The non-linear Markov process $\mu_{\nu,\lambda_{\nu}\left(  \cdot\right)
}\left(  t\right)  $ is the main object of the present paper. Therefore we
will give now another definition of this process, via jump rates of
transitions during the infinitesimal time, $\Delta t.$ So suppose that our
process is in the state $\mu\in\mathcal{M}\left(  \Omega\right)  ,$ and
assumes the value $\omega=\left(  n,\tau\right)  \in\Omega.$ During the time
increment $\Delta t$ the following two transitions can happen with
probabilities of order of $\Delta t:$

\begin{itemize}
\item the customer under the service will finish it and will leave the server,
so the value $\omega=\left(  n,\tau\right)  $ will become $\left(
n-1,\varsigma\right)  ,$ with $\varsigma\leq\Delta t.$ The probability of this
event is
\[
c_{1}\Delta t+o\left(  \Delta t\right)  ,
\]
where
\[
c_{1}=c_{1}\left(  \omega\right)  =\lim_{\Delta t\rightarrow0}\frac{1}{\Delta
t}\frac{\int_{\tau}^{\tau+\Delta t}p\left(  x\right)  \,dx}{\int_{\tau
}^{\infty}p\left(  x\right)  \,dx},
\]
while for $\omega=\mathbf{0}$ we put $c_{1}\left(  \mathbf{0}\right)  =0;$

\item a new customer will arrive to the server, so the value $\omega=\left(
n,\tau\right)  $ will become $\left(  n+1,\tau+\Delta t\right)  .$ The
probability of this event is given by
\[
c_{2}\Delta t+o\left(  \Delta t\right)  ,
\]
where the rate $c_{2}$ depends on the whole state $\mu,$ and does not depend
on $\omega:$
\[
c_{2}=c_{2}\left(  \mu\right)  =\mathbb{E}_{\mu}\left(  c_{1}\left(
\omega\right)  \right)  .
\]
In words, the input rate is the average rate of the output in the state $\mu.$
\end{itemize}

\bigskip

It is curious to note that while for the general nonlinear continuous time
Markov processes its rates depend both on the configuration and on the state
of the process, in our case the rate $c_{1}$ depends only on the
configuration, while the rate $c_{2}$ -- only on the state of the process.

\section{More facts from \cite{KR1} \label{KR1}}

Consider the following continuous time Markov process $\mathfrak{M}$. Let
there be $M$ servers and $N$ customers. The serving times are i.i.d., with
distribution $\eta.$ The configuration of the system consists of specifying
the numbers of customers $n_{i},\;i=1,...,M,$ waiting at each server, plus the
duration $\tau_{i}$ of the service time for every customer under service.
Therefore it is a point in
\[
\Theta_{M,N}=\left\{  \left(  \omega_{1},...,\omega_{M}\right)  \in\Pi
_{i=1}^{M}\Omega_{i}:n_{1}+...+n_{M}=N\right\}  .
\]
Upon being served, the customer goes to one of $M$ servers with equal
probability $1/M,$ and is there the last in the queue.

The permutation group $\mathcal{S}_{M}$ acts on $\Theta_{M,N},$ leaving the
transition probabilities invariant. Therefore we can consider the
factor-process. Its values are (unordered) finite subsets of $\Omega.$ It can
be equivalently described as a measure
\[
\nu=\frac{1}{M}\sum_{i=1}^{M}\delta_{\left(  n_{i},\tau_{i}\right)  }.
\]
We will identify such measures with the configurations of the symmetrized
factor-process. Note that
\[
\left\langle n\right\rangle _{\nu}=\frac{N}{M}.
\]
So if we use the notation $\mathcal{M}_{\rho}\left(  \Omega\right)
\subset\mathcal{M}\left(  \Omega\right)  $ for the measures $\mu$ on $\Omega$
for which $\left\langle n\right\rangle _{\mu}=\rho,$ then we have that $\nu
\in\mathcal{M}_{\frac{N}{M}}\left(  \Omega\right)  .\;$We also introduce the
notation $\mathcal{M}_{\frac{N}{M},M}\left(  \Omega\right)  \subset
\mathcal{M}_{\frac{N}{M}}\left(  \Omega\right)  $ for the family of atomic
measures, such that each atom has a weight $\frac{k}{M}$ for some integer $k.$

A state of our Markov process is a probability measure on the set of
configurations, i.e. an element of $\mathcal{M}\left(  \mathcal{M}\left(
\Omega\right)  \right)  .$ If the initial state of the process is supported by
$\mathcal{M}_{\rho}\left(  \Omega\right)  ,$ then at any positive time it is
still the element of $\mathcal{M}\left(  \mathcal{M}_{\rho}\left(
\Omega\right)  \right)  .$ A natural embedding $\mathcal{M}\left(
\Omega\right)  \subset\mathcal{M}\left(  \mathcal{M}\left(  \Omega\right)
\right)  ,$ which to each configuration $\nu\in\mathcal{M}\left(
\Omega\right)  $ corresponds the atomic measure $\delta_{\nu},$ will be
denoted by $\delta. $

For $\mu_{0}=\delta_{\nu}\in\mathcal{M}\left(  \mathcal{M}_{\frac{N}{M}%
,M}\left(  \Omega\right)  \right)  $ to be the initial state of our Markov
process, we denote by $\mu_{t}$ the evolution of this state. (Clearly, in
general $\mu_{t}\notin\delta\left(  \mathcal{M}\left(  \Omega\right)  \right)
$ for positive $t.)$ This process is ergodic. We denote by $\pi_{M,N}$ the
stationary measure of this process.

Let now $\kappa\in\mathcal{M}_{\rho}\left(  \Omega\right)  $ be some measure,
let the sequences of integers $N_{j},M_{j}\rightarrow\infty$ be such that
$\frac{N_{j}}{M_{j}}\rightarrow\rho,$ and let the measures $\nu^{j}%
\in\mathcal{M}_{\frac{N_{j}}{M_{j}},M_{j}}\left(  \Omega\right)  $ be such
that $\nu_{j}\rightarrow\kappa$ weakly. Consider the Markov processes $\mu
_{t}^{j}\in\mathcal{M}\left(  \mathcal{M}_{\frac{N_{j}}{M_{j}},M_{j}}\left(
\Omega\right)  \right)  ,$ corresponding to the initial conditions
$\delta_{\nu^{j}}.$ As we just said, in general $\mu_{t}^{j}\notin
\delta\left(  \mathcal{M}_{\frac{N_{j}}{M_{j}},M_{j}}\left(  \Omega\right)
\right)  $ for each $j,$ once $t>0.$ However, for the limit $\mu_{t}%
=\lim_{j\rightarrow\infty}\mu_{t}^{j}$ we have that $\mu_{t}\in\mathcal{M}%
\left(  \mathcal{M}_{\rho}\left(  \Omega\right)  \right)  ,$ and moreover
$\mu_{t}\in\delta\left(  \mathcal{M}_{\rho}\left(  \Omega\right)  \right)  ,$
so we can say that the random evolutions $\mu_{t}^{j}$ tend to the non-random
evolution $\mu_{t}$ with $\mu_{0}\equiv\kappa,$ as $N_{j},M_{j}\rightarrow
\infty.$

Therefore we have a dynamical system
\begin{equation}
\mathcal{T}_{t}:\mathcal{M}_{\rho}\left(  \Omega\right)  \rightarrow
\mathcal{M}_{\rho}\left(  \Omega\right)  .\label{160}%
\end{equation}
This dynamical system $\mu_{t}$ is nothing else but the non-linear Markov
process, discussed above.

Another way of obtaining the same dynamical system is to look on the behavior
of a given server. Here instead of taking the symmetrization of the initial
process $\mathfrak{M}$ on $\Theta_{M,N},$ we have to consider its projection
on the first coordinate, $\Omega_{1}$, say. To make the correspondence with
the above, we have to take for the initial state of this process a measure
$\tilde{\nu}^{j}$ on $\Theta_{M,N},$ which is $\mathcal{S}_{M}$-invariant, and
which symmetrization is the initial state $\nu^{j}$ of the preceding
paragraph. The projection of $\mathfrak{M}$ on $\Omega_{1}$ would not be, of
course, a Markov process. However, it becomes the very same non-linear Markov
process $\mu_{t}$ in the above limit $N_{j},M_{j}\rightarrow\infty$.

We can generalize further, and study the projection of $\mathfrak{M}$ to a
finite product, $\prod_{j=1}^{r}\Omega_{j}.$ Then in the limit $N_{j}%
,M_{j}\rightarrow\infty$ this projection converges to a process on
$\prod_{j=1}^{r}\Omega_{j},$ which factors into the product of $r$ independent
copies of the same non-linear Markov process $\mu_{t}.$ This statement is
known as the \textquotedblleft propagation of chaos\textquotedblright\ property.

The main result of the present paper is in particular to give conditions under
which for every $\rho$ the dynamical system $\left(  \ref{160}\right)  $ has
exactly one fixed point $\nu_{c},$ $c=c\left(  \rho\right)  $, and that it is
globally attractive. That would imply that $\pi_{N_{j},M_{j}}\rightarrow
\nu_{c},$ provided $\frac{N_{j}}{M_{j}}\rightarrow\rho$ as $j\rightarrow
\infty$ and $c=c\left(  \rho\right)  .$

\section{Main result \label{main}}

Our main result states that under certain restrictions PH\ holds. In view of
what was said in the beginning of the Introduction, it is sufficient to prove
the following:

\begin{theorem}
Consider the system, described in the Section \ref{notation}, and let its
service time $\eta$ has properties 1-6 of this Section. For every initial
state $\nu$ with finite expected service time and finite mean queue:
\begin{equation}
S\left(  \nu\right)  <\infty,~\ N\left(  \nu\right)  <\infty,\label{182}%
\end{equation}
the solution $\lambda_{\nu}\left(  \cdot\right)  $ of the equation
\[
A\left(  \nu,\lambda\left(  \cdot\right)  \right)  =\lambda\left(
\cdot\right)
\]
has the \textbf{relaxation }property:
\[
\lambda_{\nu}\left(  t\right)  \rightarrow c\text{ as }t\rightarrow\infty,
\]
where the constant $c$ satisfies
\[
N\left(  \nu\right)  \equiv\mathbb{E}_{\nu}\left(  N\left(  \omega\right)
\right)  =N\left(  \nu_{c}\right)  \equiv\mathbb{E}_{\nu_{c}}\left(  N\left(
\omega\right)  \right)  .
\]
Moreover, $\mu_{\nu,\lambda_{\nu}\left(  \cdot\right)  }\left(  t\right)
\rightarrow\nu_{c}$ weakly, as $t\rightarrow\infty.$
\end{theorem}

In particular, if the service time $\eta$ has the property that $R_{\eta
}\left(  \tau\right)  <\bar{C}$ for all $\tau$ (see $\left(  \ref{11}\right)
$), then the relaxation property holds for every initial state $\nu.$

A special case of the above theorem is the following

\begin{proposition}
\label{P} Let $T>0$ be some time moment, and suppose that the function
$\lambda\left(  \cdot\right)  $ satisfies
\begin{equation}
\lambda\left(  t\right)  =b\left(  t\right)  \text{ for all }t\geq
T,\label{31}%
\end{equation}
where
\begin{equation}
b\left(  \cdot\right)  =A\left(  \mathbf{0},\lambda\left(  \cdot\right)
\right)  .\label{30}%
\end{equation}
Suppose that
\[
\int_{0}^{T}\lambda\left(  t\right)  \,dt\leq C<\infty.
\]
Then for some $c\geq0$%
\begin{equation}
\lambda\left(  t\right)  \rightarrow c\text{ as }t\rightarrow\infty.\label{32}%
\end{equation}

\end{proposition}

Our theorem follows from the Proposition \ref{P} immediately in the special
case when the initial state $\nu$ is of the form $\nu=\mu_{\mathbf{0}%
,\lambda\left(  \cdot\right)  }\left(  t\right)  $ for some $\lambda$ and some
$t>0.$ These initial states are easier to handle, so we treat them separately.

The heuristics behind the Proposition \ref{P} is the following. One expects
that if
\[
b\left(  \cdot\right)  =A\left(  \nu,\lambda\left(  \cdot\right)  \right)  ,
\]
then the function $b$ for large times is \textquotedblleft closer to a
constant\textquotedblright\ than the function $\lambda.$ More precisely, if
$t\ $belongs to some segment $\left[  T_{1},T_{2}\right]  $, with $T_{1}\gg1,$
then the dependence of $b\left(  t\right)  $ on $\nu$ is very weak, so $b$ is
determined mainly by $\lambda.$ One then argues that once the segment $\left[
T_{1},T_{2}\right]  $ is large, $\sup_{t\in\left[  T_{1},T_{2}\right]
}b\left(  t\right)  $ should be strictly less than $\sup_{t\in\left[
T_{1},T_{2}\right]  }\lambda\left(  t\right)  .$ Indeed, one can visualize the
random configuration of the exit moments $y_{i}$-s as being obtained from the
input moments configuration of $x_{i}$-s by making it \textit{sparser}.
Namely, we have to consider a sequence $\eta_{i}$ of i.i.d. random variables,
having the same distribution as $\eta,$ and then to shift the particles
$x_{i}$ to the right, positioning them at locations $z_{i},$ so that in the
result
\begin{equation}
z_{i+1}-z_{i}\geq\eta_{i}\label{10}%
\end{equation}
for all $i$-s, and $y_{i}=z_{i}+\eta_{i}.$\textit{\ }Note that before the
shift it might have been that $x_{i+1}-x_{i}<\eta_{i}$ for some $i$-s, see
$\left(  \ref{07}\right)  $, $\left(  \ref{06}\right)  $ below for more
details. However this is a very rough idea, since some particles need not be
moved, due to the fact that $x_{i+1}-x_{i}\geq\eta_{i}$ may hold already
before to the sparsening step, in which case it will happen that
$z_{i+1}=x_{i+1},$ while $z_{i}>x_{i},$ and so the configuration becomes
locally denser. (And if $\lambda$ is a constant, then $b$ is this same
constant, so again the above argument is not literally true.)

To be more precise, we will show the following\textbf{\ self-averaging
property}. Let the functions $\lambda\left(  \cdot\right)  $ and $b\left(
\cdot\right)  $ are related by
\[
b\left(  \cdot\right)  =A\left(  \mathbf{0},\lambda\left(  \cdot\right)
\right)  .
\]
One of the main points of the following will be to show that for every value
$x$ one can find a probability density $q_{\lambda,x}\left(  t\right)  ,$
vanishing for $t\leq0,$ such that
\begin{equation}
b\left(  x\right)  =\left[  \lambda\ast q_{\lambda,x}\right]  \left(
x\right)  .\label{34}%
\end{equation}
We then will show that this self-averaging property of the system implies
(\ref{32}), provided we know in advance certain regularity properties of the
family $\left\{  q_{\lambda,x}\right\}  $. Note that apriori the condition
(\ref{34}) is not evident at all for our FIFO (=First-In-First-Out) system:
one has to rule out the situation that, say, the input rate function $\lambda$
is uniformly bounded from above by $1/3,$ while the output rate $b$ is
occasionally reaching the level $2/3;$ this is clearly inconsistent with
(\ref{34}).

In general case, when
\[
b\left(  \cdot\right)  =A\left(  \mu,\lambda\left(  \cdot\right)  \right)
\]
and $S\left(  \mu\right)  <\infty,$ we have
\begin{equation}
b\left(  x\right)  =\left(  1-\varepsilon_{\lambda,\mu}\left(  x\right)
\right)  \left[  \lambda\ast q_{\lambda,\mu,x}\right]  \left(  x\right)
+\varepsilon_{\lambda,\mu}\left(  x\right)  Q_{\lambda,\mu}\left(  x\right)
,\label{134}%
\end{equation}
where $\varepsilon_{\lambda,\mu}\left(  x\right)  >0,$ $\varepsilon
_{\lambda,\mu}\left(  x\right)  \rightarrow0$ as $x\rightarrow\infty,$ while
$Q_{\lambda,\mu}\left(  x\right)  $ is a bounded term,\textbf{\ }see Section
\ref{gencase} for details.

\section{The self-averaging relation \label{self-av}}

Here we will derive a formula, expressing the function $b\left(  \cdot\right)
=A\left(  \mathbf{0},\lambda\left(  \cdot\right)  \right)  $ in terms of the
functions $\lambda\left(  \cdot\right)  $ and the density $p\left(
\cdot\right)  $ of $\eta.$ This will be the needed self-averaging relation
(\ref{34}).

First, we define the kernels $q_{\lambda,x}\left(  t\right)  .$ To do it, let
$e\left(  u\right)  $ be the probability that our server is idle at the time
$u.$ (Note that the dependence of $e\left(  u\right)  $ on $\lambda$ is only
via $\left\{  \lambda\left(  s\right)  ,s\leq u\right\}  .$) Now define the
function $c\left(  u,t\right)  $ as follows. Let us condition on the event
that the server is idle just before time $u,$ while at $u$ the customer
arrives. Under this condition define
\begin{equation}
c\left(  u,t\right)  =\lim_{h\searrow0}\frac{1}{h}\mathbf{\Pr}\left\{
\begin{array}
[c]{c}%
\text{the server is never idle during }\left[  u,u+t\right]  ;\text{ }\\
\text{during }\left[  u+t,u+t+h\right]  \text{ the server gets }\\
\text{through with some client}%
\end{array}
\right\}  .\label{008}%
\end{equation}
Then
\begin{equation}
q_{\lambda,x}\left(  t\right)  =e\left(  x-t\right)  c\left(  x-t,t\right)
.\label{005}%
\end{equation}

\begin{theorem}
\label{form} Let the functions $b\left(  \cdot\right)  $ and $\lambda\left(
\cdot\right)  $ are related by
\[
b\left(  \cdot\right)  =A\left(  \mathbf{0},\lambda\left(  \cdot\right)
\right)  .
\]
Then also
\begin{equation}
b\left(  x\right)  =\int_{0}^{\infty}\lambda\left(  x-t\right)  q_{\lambda
,x}\left(  t\right)  \,dt.\label{007}%
\end{equation}
Moreover, for all $\lambda,x$%
\begin{equation}
\int_{0}^{\infty}q_{\lambda,x}\left(  t\right)  \,dt\leq1.\label{006}%
\end{equation}

\end{theorem}

\textbf{Note. }The relation $\left(  \ref{006}\right)  $ is not at all
obvious, and we see no way to derive it from $\left(  \ref{005}\right)  $
without going into the details of the serving mechanism. In fact, it does not
hold for some more general models.\smallskip

\begin{proof}
We first introduce some new notions.

Let $l_{1},...,l_{n}>0$ be a collection of positive real numbers, which we
will interpret as the lengths of hard rods ($\equiv$service times), placed in
$\mathbb{R}^{1}.$ A configuration of rods can be then given by specifying the
sequence $x_{1},x_{2},...,x_{n}$ of their left-ends: the rod $l_{i}$ occupies
the segment $\left[  x_{i},x_{i}+l_{i}\right]  .$ This configuration will be
denoted by $\sigma_{n}\left(  x_{1},x_{2},...,x_{n};l_{1},...,l_{n}\right)  .$

In case some of the rods from $\sigma_{n}\left(  x_{1},x_{2},...,x_{n}%
;l_{1},...,l_{n}\right)  $ are intersecting over a nondegenerate segments, we
say that such a configuration has conflicts. By a resolution of conflicts we
call another placement of the rods $l_{1},...,l_{n}$ on the line. To define
it, we first need to reenumerate the points of the sequence $x_{1}%
,x_{2},...,x_{n}$ so that it will become increasing. To save on notation we
suppose \textit{until the end of this paragraph only} that it is initially so.
Then the new placing of the rods have the following sequence $z_{1}%
<z_{2}<...<z_{n}$ of the left-ends:

\noindent it is defined inductively by
\[
z_{1}=x_{1},
\]
and
\begin{equation}
z_{i}=\max\left\{  z_{i-1}+l_{i-1},x_{i}\right\} \label{07}%
\end{equation}
(Lindley equation). We will denote by $y$-s the corresponding set of the
right-ends:
\begin{equation}
y_{i}=z_{i}+l_{i}.\label{06}%
\end{equation}
Any configuration with no conflicts, and in particular any configuration
obtained by resolution of the conflicting one, will be called an
\textbf{r-configuration}. The operation of resolving the conflict will be
denoted by $R,$ so
\[
\sigma_{n}\left(  z_{1},z_{2},...,z_{n};l_{1},...,l_{n}\right)  =R\sigma
_{n}\left(  x_{1},x_{2},...,x_{n};l_{1},...,l_{n}\right)  .
\]
(The $R$-operation is not well defined for $x$-sequences with coinciding
entries. However, they have zero Lebesgue measure, which makes them irrelevant
for our future needs.)

For any configuration $\sigma$ of rods we will denote by $Y\left(
\sigma\right)  $ the set of their right-ends. So, in our notations
\[
\left(  y_{1},...,y_{n}\right)  =Y\left(  R\sigma_{n}\left(  x_{1}%
,x_{2},...,x_{n};l_{1},...,l_{n}\right)  \right)  .
\]

Suppose now that the lengths $l_{1},...,l_{n},$ the set $x_{1},x_{2}%
,...,x_{n-1}$ (with $n-1$ points) and the location $y\in\mathbb{R}^{1}$ are
specified. We define the values $X\left(  y\right)  \equiv X\left(
y\Bigm|x_{1},x_{2},...,x_{n-1};l_{1},...,l_{n}\right)  \in\mathbb{R}^{1}$ as
the solutions of the equation
\begin{equation}
y\in Y\left(  R\sigma_{n}\left(  x_{1},x_{2},...,x_{n-1},X\left(  y\right)
;l_{1},...,l_{n}\right)  \right)  .\label{21}%
\end{equation}
It is clear that the function $X\left(  y\Bigm|x_{1},x_{2},...,x_{n-1}%
;l_{1},...,l_{n}\right)  $ is not defined everywhere, and on the set where it
is defined, it is multivalued, provided $n\geq2$. (The case $n=1$ is trivial:
$X\left(  y\Bigm|l_{1}\right)  =y-l_{1}.$) However, outside the set of $x$-s
and $l$-s of Lebesgue measure zero in $\mathbb{R}^{2n-1}$, which set is
irrelevant for our future purposes, its multivaluedness is reduced to ``finitely-many-valuedness''.

Now we can write the desired formula:
\begin{align}
& b\left(  y\right)  =\exp\left\{  -I_{\lambda}\left(  y\right)  \right\}
\sum_{n=1}^{\infty}\frac{1}{\left(  n-1\right)  !}\times\label{129}\\
& \times\underset{n}{\underbrace{\int_{0}^{\infty}...\int_{0}^{\infty}}%
}\left[  \underset{n-1}{\underbrace{\int_{0}^{y}...\int_{0}^{y}}}%
\lambda\left(  X\left(  y\Bigm|\left\{  x_{1},...,x_{n-1}\right\}
;l_{1},...,l_{n}\right)  \right)  \prod_{i=1}^{n-1}\lambda\left(
x_{i}\right)  \,dx_{i}\right]  \prod_{i=1}^{n}p\left(  l_{i}\right)
\,dl_{i},\nonumber
\end{align}
where
\[
I_{\lambda}\left(  y\right)  =\int_{0}^{y}\lambda\left(  x\right)  \,dx.
\]
The integral in (\ref{129}) should be understood as follows: the range of
integration coincides with the domain where the function $X\left(
y\Bigm|\left\{  x_{1},...,x_{n-1}\right\}  ;l_{1},...,l_{n}\right)  $ is
defined, while over the domains where the function $X$ is multivalued one
should integrate each branch separately and then take the sum of integrals.

In words, the meaning of the relation $\left(  \ref{129}\right)  $ is the
following: for every realization $x_{1},...,x_{n-1}$ of the Poisson random
field and every realization $l_{1},...,l_{n}$ of the sequence of the service
times, we look for time moments $X=X\left(  y\Bigm|x_{1},...,x_{n-1}%
;l_{1},...,l_{n}\right)  ,$ at which the $l_{n}$-customer has to arrive, so as
to ensure that at the moment $y$ some customer (perhaps a different one) will
exit, after being served. That is, moments $X=X\left(  y\Bigm|x_{1}%
,...,x_{n-1};l_{1},...,l_{n}\right)  $ are beginnings of busy periods, during
which there happens an exit at time $y.$ In some cases such moments might not
exist, while in other cases there might be more than one such moment. If
$X_{i}$ are these moments, we then have to add all the rate values,
$\lambda\left(  X_{i}\right)  ,$ and to integrate the sum $\sum_{i}%
\lambda\left(  X_{i}\right)  $ over all $n$ and all $x_{1},...,x_{n-1}%
;l_{1},...,l_{n},$ thus getting the exit rate $b\left(  y\right)  .$ A
one-second thought will convince the reader that the formula $\left(
\ref{129}\right)  $ contains in itself another definition of the kernel
$\left(  \ref{005}\right)  ,$ together with the proof of the relation $\left(
\ref{007}\right)  .$ So we need only to prove $\left(  \ref{006}\right)  ,$
which turns out to be quite delicate.

The first summand $\left(  n=1\right)  $ in the sum in (\ref{129}) is by
definition the convolution,
\begin{equation}
b_{1}\left(  y\right)  =\int_{0}^{y}\lambda\left(  y-l\right)  p\left(
l\right)  \,dl.\label{18}%
\end{equation}
Since $p\left(  l\right)  \geq0$ and
\begin{equation}
\int_{0}^{y}p\left(  l\right)  \,dl\leq1,\label{17}%
\end{equation}
we have indeed that $b_{1}\left(  y\right)  <\sup_{x\leq y}\lambda\left(
x\right)  $ in case when, say, the maxima of $\lambda$ are isolated, or when
$\lambda$ is not a constant and the support of the distribution $p$ is the
full semiaxis $\left\{  l>0\right\}  .$ We want to show that in some sense the
same is true for all the functions $b_{n},$ defined as
\begin{equation}
b_{n}\left(  y\right)  =\int\left[  \int\lambda\left(  X\left(  y\Bigm|x_{1}%
,...,x_{n-1};l_{1},...,l_{n}\right)  \right)  \prod_{i=1}^{n-1}\left(
\frac{\lambda\left(  x_{i}\right)  }{I_{\lambda}\left(  y\right)  }%
\,dx_{i}\right)  \right]  \prod_{i=1}^{n}p\left(  l_{i}\right)  \,dl_{i}%
.\label{20}%
\end{equation}
Since
\[
b\left(  y\right)  =\exp\left\{  -I_{\lambda}\left(  y\right)  \right\}
\sum_{n=1}^{\infty}\frac{I_{\lambda}\left(  y\right)  ^{n-1}}{\left(
n-1\right)  !}b_{n}\left(  y\right)  ,
\]
the crucial step will be the analog of (\ref{18}), (\ref{17}) for all $n>1,$
that is that
\[
b_{n}\left(  y\right)  =\int_{0}^{y}\lambda\left(  y-l\right)  p_{n}\left(
l\right)  \,dl,
\]
for some $p_{n}\left(  l\right)  \geq0,\;\int_{0}^{y}p_{n}\left(  l\right)
\,dl\nearrow1$ for $y\rightarrow\infty.$ This turns out to be quite an
involved combinatorial statement.

Note that, evidently, the measure $\prod_{i=1}^{n}p\left(  l_{i}\right)
\,dl_{i}$ is invariant under the coordinate permutations in $\mathbb{R}^{n};$
therefore we can rewrite the expression (\ref{20}) for the function
$b_{n}\left(  y\right)  $ as
\begin{equation}
b_{n}\left(  y\right)  =\int\left[  \int\frac{1}{n!}\lambda\left(  \bar
{X}\left(  y\Bigm|x_{1},...,x_{n-1};\left\{  l_{1},...,l_{n}\right\}  \right)
\right)  \prod_{i=1}^{n-1}\left(  \frac{\lambda\left(  x_{i}\right)
}{I_{\lambda}\left(  y\right)  }\,dx_{i}\right)  \right]  \prod_{i=1}%
^{n}p\left(  l_{i}\right)  \,dl_{i},\label{23}%
\end{equation}
where the following notations and conventions are used:

\begin{itemize}
\item the (multivalued) function $\bar{X}\left(  y\Bigm|x_{1},...,x_{n-1}%
;\left\{  l_{1},...,l_{n}\right\}  \right)  $ by definition assigns to every
$y$ the union of the sets of solutions $X\left(  y\right)  $ of all the
equations
\begin{equation}
y\in Y\left(  R\sigma_{n}\left(  x_{1},...,x_{n-1},X\left(  y\right)
;l_{\pi\left(  1\right)  },...,l_{\pi\left(  n\right)  }\right)  \right)
,\label{22}%
\end{equation}
with $\pi$ running over all the permutation group $\mathcal{S}_{n}$ (the
notation $\left\{  l_{1},...,l_{n}\right\}  $ stresses the fact that the
function $\bar{X}$ does not depend on the order of $l_{i}$-s);

\item the entries of the set $\bar{X}\left(  y\Bigm|x_{1},...,x_{n-1};\left\{
l_{1},...,l_{n}\right\}  \right)  $ have to be counted with multiplicities,
which for a given $x\in\bar{X}\left(  y\Bigm|x_{1},...,x_{n-1};\left\{
l_{1},...,l_{n}\right\}  \right)  $ is by definition the number of equations
(\ref{22}) with different $\pi$-s, to which $x$ is a solution;

\item the integration in (\ref{23}) of the multivalued function means that
each sheet should be integrated and the results added. Moreover, each sheet
has to be taken as many times as its multiplicity is.
\end{itemize}

Since each contribution $\lambda\left(  X\left(  y\Bigm|x_{1},...,x_{n-1}%
;l_{1},...,l_{n}\right)  \right)  $ to $\left(  \ref{20}\right)  $ appears
$n!$ times in $\left(  \ref{23}\right)  ,$ we have to divide by $n!.$

We repeat that while for some $x$-s, $\pi$-s and $l$-s the equation (\ref{22})
might have no solutions, for other data it can have more than one solution.
Clearly, the set $\bar{X}\left(  y\Bigm|x_{1},...,x_{n-1};\left\{
l_{1},...,l_{n}\right\}  \right)  $, for Lebesgue-almost every data
$x_{1},...,x_{n-1},$ can have no other entries than those of the form
\[
x_{A,y,\left\{  l_{i}\right\}  }=y-\sum_{i\in A\subset\left\{
1,2,...,n\right\}  }l_{i},
\]
where $A$ runs over all nonempty subsets of $\left\{  1,2,...,n\right\}  $
(i.e. at most $2^{n}-1$ different entries). So the function $\bar{X}\left(
y\Bigm|x_{1},...,x_{n-1};\left\{  l_{1},...,l_{n}\right\}  \right)  $, as a
function of $x_{1},...,x_{n-1},$ has to be piecewise constant. It is not ruled
out apriori that for some data the set $\bar{X}\left(  y\Bigm|x_{1}%
,...,x_{n-1};\left\{  l_{1},...,l_{n}\right\}  \right)  $ can be empty. This
is not, however, the case. Moreover, as the crucial Theorem \ref{T6} below states,

\begin{itemize}
\item the number of elements in the set $\bar{X}\left(  y\Bigm|x_{1}%
,...,x_{n-1};\left\{  l_{1},...,l_{n}\right\}  \right)  ,$ counted with
multiplicities, is \textbf{precisely }$n!$ for Lebesgue-almost every value of
the arguments.
\end{itemize}

Therefore we have for the inner integral in (\ref{23}):
\begin{align*}
& \int\frac{1}{n!}\lambda\left(  \bar{X}\left(  y\Bigm|x_{1},...,x_{n-1}%
;\left\{  l_{1},...,l_{n}\right\}  \right)  \right)  \prod_{i=1}^{n-1}\left(
\frac{\lambda\left(  x_{i}\right)  }{I_{\lambda}\left(  y\right)  }%
\,dx_{i}\right) \\
& =\int\frac{1}{n!}\sum_{\substack{A\subset\left\{  1,2,...,n\right\}  ,
\\A\neq\emptyset}}k\left(  A,y,x_{1},...,x_{n-1};\left\{  l_{1},...,l_{n}%
\right\}  \right)  \lambda\left(  x_{A,y,\left\{  l_{i}\right\}  }\right)
\prod_{i=1}^{n-1}\left(  \frac{\lambda\left(  x_{i}\right)  }{I_{\lambda
}\left(  y\right)  }\,dx_{i}\right)  ,
\end{align*}
where the integer $k\left(  A,y,x_{1},...,x_{n-1};\left\{  l_{1}%
,...,l_{n}\right\}  \right)  $ is the multiplicity of the value
$x_{A,y,\left\{  l_{i}\right\}  }$ of the function $\bar{X}$ at the point
$\left(  y,x_{1},...,x_{n-1};\left\{  l_{1},...,l_{n}\right\}  \right)  .$
Rewriting it as
\begin{align*}
& \int\frac{1}{n!}\lambda\left(  \bar{X}\left(  y\Bigm|x_{1},...,x_{n-1}%
;\left\{  l_{1},...,l_{n}\right\}  \right)  \right)  \prod_{i=1}^{n-1}\left(
\frac{\lambda\left(  x_{i}\right)  }{I_{\lambda}\left(  y\right)  }%
\,dx_{i}\right) \\
& =\sum_{\substack{A\subset\left\{  1,2,...,n\right\}  , \\A\neq\emptyset}%
}q_{\lambda,y}\left(  A\Bigm|\left\{  l_{1},...,l_{n}\right\}  \right)
\lambda\left(  x_{A,y,\left\{  l_{i}\right\}  }\right)  ,
\end{align*}
where
\begin{align}
& q_{\lambda,y}\left(  A\Bigm|\left\{  l_{1},...,l_{n}\right\}  \right)
\label{25}\\
& =\int\frac{1}{n!}k\left(  A,y,x_{1},...,x_{n-1};\left\{  l_{1}%
,...,l_{n}\right\}  \right)  \prod_{i=1}^{n-1}\left(  \frac{\lambda\left(
x_{i}\right)  }{I_{\lambda}\left(  y\right)  }\,dx_{i}\right)  ,\nonumber
\end{align}
we have, due to the fact that Lebesgue-a.e.
\[
\sum_{\substack{A\subset\left\{  1,2,...,n\right\}  , \\A\neq\emptyset
}}k\left(  A,y,x_{1},...,x_{n-1};\left\{  l_{1},...,l_{n}\right\}  \right)
=n!,
\]
the relations
\begin{equation}
0\leq q_{\lambda,y}\left(  A\Bigm|\left\{  l_{1},...,l_{n}\right\}  \right)
\leq1,\text{ with }\sum_{\substack{A\subset\left\{  1,2,...,n\right\}  ,
\\A\neq\emptyset}}q_{\lambda,y}\left(  A\Bigm|\left\{  l_{1},...,l_{n}%
\right\}  \right)  =1,\label{26}%
\end{equation}
since the measures $\frac{\lambda\left(  x_{i}\right)  }{I_{\lambda}\left(
y\right)  }\,dx_{i}$ are probability measures on $\left[  0,y\right]  $. (Note
that the functions $k\left(  A,y,x_{1},...,x_{n-1};\left\{  l_{1}%
,...,l_{n}\right\}  \right)  $ do depend on the variables $x_{1},...,x_{n-1};$
hence the measures $q_{\lambda,y}\left(  \cdot\Bigm|\left\{  l_{1}%
,...,l_{n}\right\}  \right)  $ indeed depend on $\lambda,y. $) Therefore, for
the function $b_{n}\left(  y\right)  $ we obtain a sort of a convolution
expression:
\begin{equation}
b_{n}\left(  y\right)  =\int\sum_{\substack{A\subset\left\{
1,2,...,n\right\}  , \\A\neq\emptyset}}\left[  q_{\lambda,y}\left(
A\Bigm|\left\{  l_{1},...,l_{n}\right\}  \right)  \lambda\left(
x_{A,y,\left\{  l_{i}\right\}  }\right)  \right]  \prod_{i=1}^{n}p\left(
l_{i}\right)  \,dl_{i}.\label{24}%
\end{equation}
Be it the case that the probability measure $q_{\lambda,y}\left(
\cdot\Bigm|\left\{  l_{1},...,l_{n}\right\}  \right)  $ is concentrated on
just one subset $A=\left\{  1,2,...,n\right\}  ,$ we would obtain the usual
convolution
\[
b_{n}\left(  y\right)  =\int\lambda\left(  y-l_{1}-...-l_{n}\right)
\prod_{i=1}^{n}p\left(  l_{i}\right)  \,dl_{i}=\lambda\ast\underset
{n}{\underbrace{p\ast...\ast p}}\left(  y\right)  .
\]
Here the situation is more subtle, and in (\ref{24}) we have a stochastic
mixture of convolutions with random number of summands.

Taking into account the relations (\ref{129}), (\ref{20}), (\ref{24}), the
result can be summarized as follows. Let $\theta\equiv\theta_{\lambda,y}$ be
the integer valued random variable with the distribution
\[
\mathbf{\Pr}\left\{  \theta=n\right\}  =\exp\left\{  -I_{\lambda}\left(
y\right)  \right\}  \frac{\left[  I_{\lambda}\left(  y\right)  \right]  ^{n}%
}{n!},\,n=0,1,2,...,
\]
and $\eta_{1},\eta_{2},...$ be the i.i.d. random serving times. Consider the
random function $\xi_{\lambda,y}=\xi_{\lambda,y}\left(  \theta_{\lambda
,y};\eta_{1},\eta_{2},...\right)  ,$ such that its conditional distribution
under condition that the realization $\theta_{\lambda,y};\eta_{1},\eta
_{2},...$ is given, is supported by the finite set
\[
L\left(  \theta_{\lambda,y};\eta_{1},\eta_{2},...\right)  =\left\{  \sum_{i\in
A}\eta_{i}:A\subset\left\{  1,2,...,\theta_{\lambda,y}\mathbf{+}1\right\}
,A\neq\emptyset\right\}  \subset\mathbb{R}^{1},
\]
and is given by
\[
\mathbf{\Pr}\left\{  \xi_{\lambda,y}=\sum_{i\in A}\eta_{i}\Bigm|\theta
_{\lambda,y};\eta_{1},\eta_{2},...\right\}  =q_{\lambda,y}\left(
A\Bigm|\left\{  \eta_{1},...,\eta_{\theta_{\lambda,y}+1}\right\}  \right)
\]
(see (\ref{25})). Then the following holds:
\[
b\left(  y\right)  =\mathbb{E}\left(  \lambda\left(  y-\xi_{\lambda,y}\right)
\right)  .
\]
This is precisely the relation (\ref{34}), with $q_{\lambda,y}$ being the
distribution of $\xi_{\lambda,y}.$ The relation $\left(  \ref{006}\right)  $
follows directly from $\left(  \ref{26}\right)  .$
\end{proof}

\section{Combinatorics of the rod placements \label{comb}}

In this section we will prove the Theorem \ref{T6}, which was used in the
previous section. We will use the notation of the previous section, introduced
in the proof of the Theorem \ref{form}, up to relation $\left(  \ref{21}%
\right)  .$

By a cluster of the r-configuration $\sigma_{n}\left(  z_{1},...,z_{n}%
;l_{1},...,l_{n}\right)  $ with $z_{1}<...<z_{n}$ we call any maximal
subsequence $z_{i}<z_{i+1}<...<z_{j}$ such that $z_{j}=z_{i}+l_{i}%
+l_{i+1}+...+l_{j-1}.$ (The segment $\left[  z_{i},z_{i}+l_{i}+l_{i+1}%
+...+l_{j}\right]  \equiv\left[  x_{i},z_{j}+l_{j}\right]  $ is what is called
"busy period" for the queue.) If $z_{i}<z_{i+1}<...<z_{j}$ is a cluster of an
r-configuration, then the point $z_{i}$ will be called the root of the
cluster, while the point $z_{j}$ will be called the head of the cluster. Note
that for Lebesgue almost every configuration $\sigma_{n}\left(  x_{1}%
,...,x_{n};l_{1},...,l_{n}\right)  $ the point $z_{i}$ is a root of a cluster
of the corresponding r-configuration if and only if $z_{i}=x_{i}.$ The segment
$\left[  z_{i},z_{j}+l_{j}\right]  $ will be called the body of the cluster
$z_{i}<z_{i+1}<...<z_{j},$ and the point $z_{j}+l_{j}$ will be called the end
of the cluster.

The notation $\sigma_{n}\left(  x_{1},...,x_{n};l_{1},...,l_{n}\right)
\cup\sigma_{1}\left(  X,L\right)  $ has the obvious meaning of adding an extra
rod of the length $L$ at the location $X.$ Note though, that in general
\[
R\left[  \sigma_{n}\left(  x_{1},...,x_{n};l_{1},...,l_{n}\right)  \cup
\sigma_{1}\left(  X,L\right)  \right]  \neq R\left[  R\sigma_{n}\left(
x_{1},...,x_{n};l_{1},...,l_{n}\right)  \cup\sigma_{1}\left(  X,L\right)
\right]  .
\]
It is however the case, if the point $X$ is outside the union of all bodies of
clusters of $R\sigma_{n}\left(  x_{1},...,x_{n};l_{1},...,l_{n}\right)  .$
This will be used later.

In what follows we will need a marked point in $\mathbb{R}^{1}.$ For all our
purposes it is convenient to chose the origin, $0\in\mathbb{R}^{1},$ as such a point.

We will say that the resolution of conflicts in the configuration

\noindent$\sigma_{n}\left(  x_{1},...,x_{n};l_{1},...,l_{n}\right)  $ results
in a \textbf{hit} of the origin, iff for some $k$ we have
\begin{equation}
y_{k}\equiv z_{k}+l_{k}=0.\label{04}%
\end{equation}
Such a hit will be called an $x_{r}$-hit, iff the cluster of the point $z_{k}
$ has its root at $z_{r}=x_{r}.$ (Necessarily, we have that $r\leq k.$) An
$x_{r}$-hit will be called an $\left(  x_{r},x_{k}\right)  $-hit, if
(\ref{04}) holds.

Now we are ready to formulate our problem. Let $n$ be an integer, and
$\lambda_{1}<\lambda_{2}<...<\lambda_{n}$ be a fixed set of positive lengths
of rods. Let $x_{1}<x_{2}<...<x_{n-1}$ be a set of $\left(  n-1\right)  $
left-ends. We want to compute the number $N\left(  x_{1},x_{2},...,x_{n-1}%
;\lambda_{1},\lambda_{2},...,\lambda_{n}\right)  ,$ which is defined as
follows. For any permutation $\pi$ of $n$ elements and for any $X\in
\mathbb{R}^{1},\,X\neq x_{1},x_{2},...,x_{n-1}$ we can consider the
configuration $\sigma_{n-1}\left(  x_{1},...,x_{n-1};\lambda_{\pi\left(
1\right)  },...,\lambda_{\pi\left(  n-1\right)  }\right)  $ $\cup\sigma
_{1}\left(  X,\lambda_{\pi\left(  n\right)  }\right)  $ of rods, when the rods
$l_{i}=\lambda_{\pi\left(  i\right)  }$ are placed at $x_{i},i=1,...,n-1,$
while the free rod $l_{n}=\lambda_{\pi\left(  n\right)  } $ is placed at $X.$
Given $\pi,$ we count the number $N_{\pi}\left(  x_{1},...,x_{n-1};\lambda
_{1},...,\lambda_{n}\right)  $ of different locations $X,$ such that the
corresponding r-configuration $R\left[  \sigma_{n-1}\left(  x_{1}%
,...,x_{n-1};\lambda_{\pi\left(  1\right)  },...,\lambda_{\pi\left(
n-1\right)  }\right)  \cup\sigma_{1}\left(  X,\lambda_{\pi\left(  n\right)
}\right)  \right]  $ has a hit, and moreover this hit is an $X$-hit. (In
certain cases one cannot produce an $X$-hit by putting the rod $l_{n}%
=\lambda_{\pi\left(  n\right)  }$ anywhere on $\mathbb{R}^{1};$ then $N_{\pi
}\left(  x_{1},...,x_{n-1};\lambda_{1},...,\lambda_{n}\right)  =0.$ In certain
other cases there are more than one possibility to place the free rod so as to
produce an $X$-hit.) Then we define
\[
N\left(  x_{1},...,x_{n-1};\lambda_{1},...,\lambda_{n}\right)  =\sum_{\pi
\in\mathcal{S}_{n}}N_{\pi}\left(  x_{1},...,x_{n-1};\lambda_{1},...,\lambda
_{n}\right)  .
\]

\begin{theorem}
\label{T6} For Lebesgue-almost every $x_{1},...,x_{n-1}$ and $\lambda
_{1},...,\lambda_{n},$%
\[
N\left(  x_{1},...,x_{n-1};\lambda_{1},...,\lambda_{n}\right)  =n!
\]

\end{theorem}

\begin{proof}
Let us explain why the result is plausible. Let the set $x_{1},...,x_{n-1}$ be
given. Then we can choose the positive numbers $\lambda_{1},...,\lambda_{n}$
so small that for any $\pi$ the configuration

\noindent$\sigma_{n-1}\left(  x_{1},...,x_{n-1};\lambda_{\pi\left(  1\right)
},...,\lambda_{\pi\left(  n-1\right)  }\right)  \cup\sigma_{1}\left(
X=-\lambda_{\pi\left(  n\right)  },\lambda_{\pi\left(  n\right)  }\right)  ,$
having the $\left(  X,X\right)  $-hit, has no conflicts, while no other choice
of $X$ results in a hit. Therefore in our case \noindent$N_{\pi}\left(
x_{1},...,x_{n-1};\lambda_{1},...,\lambda_{n}\right)  =1$ for every $\pi,$ so indeed

\noindent$N\left(  x_{1},...,x_{n-1};\lambda_{1},...,\lambda_{n}\right)  =n! $.

Now we explain why our result is non-trivial. To see it, take $n=2,$
$x_{1}=-3,$ $\lambda_{1}=1,$ $\lambda_{2}=10.$ Then
\[
N_{12}\left(  x_{1};\lambda_{1},\lambda_{2}\right)  =2
\]
-- one can place the rod $10$ at $-10$ or at $-11.$ On the other hand,
\[
N_{21}\left(  x_{1};\lambda_{1},\lambda_{2}\right)  =0
\]
-- the rod $10,$ placed at $-3,$ blocks the origin from being hit. Still,
$2+0=2!.$ Note that this example is a general position one.

We will derive our theorem from its special case, explained in the first
paragraph of the present proof. The idea of computing $N\left(  x_{1}%
,...,x_{n-1};\lambda_{1},...,\lambda_{n}\right)  $ for a general data is to
decrease one by one the numbers $\lambda_{1}<\lambda_{2}<...<\lambda_{n},$
starting from the smallest one, to the values very small, keeping track on the
quantities $N_{\pi}\left(  x_{1},...,x_{n-1};\lambda_{1},...,\lambda
_{n}\right)  .$ During this evolution some of these will jump, but the total
sum $N\left(  x_{1},...,x_{n-1};\lambda_{1},...,\lambda_{n}\right)  $ would
stay unchanged, as we will show. That will prove our theorem.

We begin by presenting a simple formula for the number $N_{\pi}\left(
x_{1},...,x_{n-1};\lambda_{1},..,\lambda_{n}\right)  .$ Consider the rod
configuration $R\left[  \sigma_{n-1}\left(  x_{1},...,x_{n-1};\lambda
_{\pi\left(  1\right)  },...,\lambda_{\pi\left(  n-1\right)  }\right)
\right]  ,$ which will be abbreviated as $R_{\pi}\left(  \lambda
_{1},..,\lambda_{n}\right)  \equiv R_{\pi}\left(  \mathbf{\lambda}\right)  .$
Let us compute the quantity $S_{\pi}\left(  x_{1},...,x_{n-1};\lambda
_{1},..,\lambda_{n}\right)  ,$ which is the number of points $y_{i}\in
Y\left(  R_{\pi}\left(  \lambda_{1},..,\lambda_{n}\right)  \right)  ,$ falling
into the segment $\left[  -\lambda_{\pi\left(  n\right)  },0\right]  .$

\begin{lemma}
\label{formula}%

\begin{equation}
N_{\pi}\left(  x_{1},...,x_{n-1};\lambda_{1},..,\lambda_{n}\right)  =\left\{
\begin{array}
[c]{ll}%
S_{\pi}\left(  x_{1},...,x_{n-1};\lambda_{1},..,\lambda_{n}\right)  &
\begin{array}
[c]{l}%
\text{ if the point }-\lambda_{\pi\left(  n\right)  }\text{ }\\
\text{belongs to a cluster }\\
\text{of }R_{\pi}\left(  \lambda_{1},..,\lambda_{n}\right)  ,
\end{array}
\\
& \\
S_{\pi}\left(  x_{1},...,x_{n-1};\lambda_{1},..,\lambda_{n}\right)  +1 &
\text{ otherwice.}%
\end{array}
\right. \label{0110}%
\end{equation}

\end{lemma}

\textbf{Proof of the Lemma \ref{formula}. }Indeed, for every $y_{i},$ falling
inside $\left[  -\lambda_{\pi\left(  n\right)  },0\right]  $, there is a position

\noindent$X_{i}\left(  z_{1},...,z_{n-1},y_{1},...,y_{n-1}\right)  <0,$ such
that once the free rod $\lambda_{\pi\left(  n\right)  }$ is placed there, the
site $y_{i}$ is pushed to the right and hits the origin. In case the point
$-\lambda_{\pi\left(  n\right)  }$ is outside all clusters of $R_{\pi}\left(
\lambda_{1},..,\lambda_{n}\right)  ,$ placing the free rod $\lambda
_{\pi\left(  n\right)  }$ at $X_{0}=-\lambda_{\pi\left(  n\right)  }$ produces
an extra hit. $\blacksquare$

Now let $\Delta>0$ be such that
\[
\lambda_{1}<\lambda_{2}<...<\lambda_{i-1}<\lambda_{i}-\Delta<\lambda
_{i}+\Delta<\lambda_{i+1}<...<\lambda_{n}%
\]
for some $i=1,...,n,$ and some of the functions $N_{\pi}$ exhibit jumps in the
variable $\lambda_{i}$ as it goes down from $\lambda_{i}+\Delta$ to
$\lambda_{i}-\Delta.$ We denote by $\mathbf{\lambda}\left(  \delta\right)  $
the vector $\lambda_{1},...,\lambda_{i}+\delta,...,\lambda_{n}.$ We suppose
that $\Delta$ is small enough, so that for any $\pi$ the difference
\[
\left\vert N_{\pi}\left(  x_{1},...,x_{n-1};\mathbf{\lambda}\left(
\Delta\right)  \right)  -N_{\pi}\left(  x_{1},...,x_{n-1};\mathbf{\lambda
}\left(  -\Delta\right)  \right)  \right\vert
\]
is at most one. Moreover, we want $\Delta$ to be so small that on the segment
$\lambda\in\left[  \lambda_{i}-\Delta,\lambda_{i}+\Delta\right]  $ there is
precisely one point, say $\lambda_{i},$ at which some of the functions
$N_{\pi}\left(  x_{1},...,x_{n-1};\mathbf{\lambda}\right)  $ do jump. (In
general, there will be several permutations $\pi,$ for which such a jump will
happen at $\lambda=\lambda_{i}.$ Indeed, if we observe an $\left(
X,x_{k}\right)  $-hit in our rod configuration with $l_{i}=\lambda_{\pi\left(
i\right)  }$, while we have that $x_{1}<x_{2}<...x_{s-1}<X<x_{s}%
<...<x_{k}<...<x_{n-1},$ then in some cases we will have an $\left(
X,x_{k}\right)  $-hit for every rearrangement of the rods $l_{s},...,l_{k},$
i.e. for all permutations of the form $\pi\circ\rho,$ where $\rho$ permutes
the elements $s,...,k,$ leaving the other fixed, see Lemma \ref{formula}.)

\textbf{\ }Let us begin with the case when
\begin{equation}
N_{\pi}\left(  x_{1},...,x_{n-1};\mathbf{\lambda}\left(  \Delta\right)
\right)  -N_{\pi}\left(  x_{1},...,x_{n-1};\mathbf{\lambda}\left(
-\Delta\right)  \right)  =1.\label{010}%
\end{equation}
That means that either the intersection $Y\left(  R_{\pi}\left(
\mathbf{\lambda}\left(  \Delta\right)  \right)  \right)  \cap\left[
-\mathbf{\lambda}\left(  \Delta\right)  _{\pi\left(  n\right)  },0\right]  $
is non-empty and after the $2\delta$-evolution its cardinality decreases by
one, or else that the point $-\lambda\left(  \Delta\right)  _{\pi\left(
n\right)  }$\ is outside all clusters of $R_{\pi}\left(  \mathbf{\lambda
}\left(  \Delta\right)  \right)  ,$ while the point $-\lambda\left(
-\Delta\right)  _{\pi\left(  n\right)  }$ is inside some cluster of $R_{\pi
}\left(  \mathbf{\lambda}\left(  -\Delta\right)  \right)  .$ In the first case
let $y_{k}\left(  \mathbf{\lambda}\left(  \Delta\right)  ,\pi\right)
<...<y_{r}\left(  \mathbf{\lambda}\left(  \Delta\right)  ,\pi\right)  \ $be
all the points of the above intersection. The relation $\left(  \ref{010}%
\right)  $ implies via $\left(  \ref{0110}\right)  $ that the point
$y_{k}\left(  \mathbf{\lambda}\left(  \delta\right)  ,\pi\right)  $ leaves the
segment $\left[  -\mathbf{\lambda}\left(  \delta\right)  _{\pi\left(
n\right)  },0\right]  $ as $\delta$ passes the zero value:%

\begin{equation}
y_{k}\left(  \mathbf{\lambda}\left(  \delta\right)  ,\pi\right)
>-\mathbf{\lambda}\left(  \delta\right)  _{\pi\left(  n\right)  }\text{ for
}\delta>0,\label{012}%
\end{equation}

\begin{equation}
y_{k}\left(  \mathbf{\lambda}\left(  0\right)  ,\pi\right)  =-\mathbf{\lambda
}\left(  0\right)  _{\pi\left(  n\right)  },\label{013}%
\end{equation}

\begin{equation}
y_{k}\left(  \mathbf{\lambda}\left(  \delta\right)  ,\pi\right)
<-\mathbf{\lambda}\left(  \delta\right)  _{\pi\left(  n\right)  }\text{ for
}\delta<0.\label{014}%
\end{equation}
Moreover, the point $y_{k}\left(  \mathbf{\lambda}\left(  \delta\right)
,\pi\right)  $ is not the end of the cluster -- otherwise we would have
$N_{\pi}\left(  x_{1},...,x_{n-1};\mathbf{\lambda}\left(  \Delta\right)
\right)  =N_{\pi}\left(  x_{1},...,x_{n-1};\mathbf{\lambda}\left(
-\Delta\right)  \right)  .$ Therefore $y_{k}\left(  \mathbf{\lambda}\left(
\delta\right)  ,\pi\right)  =z_{k+1}\left(  \mathbf{\lambda}\left(
\delta\right)  ,\pi\right)  .$ We now claim that if we assign the rod
$\mathbf{\lambda}\left(  \delta\right)  _{\pi\left(  n\right)  }$ to
$x_{k+1},$ and will take for the free rod the rod $\mathbf{\lambda}\left(
\delta\right)  _{\pi\left(  k+1\right)  }$, then for the corresponding
permutation the opposite to $\left(  \ref{010}\right)  $ happens:
\begin{equation}
N_{\pi\circ\left(  n\leftrightarrow k+1\right)  }\left(  x_{1},...,x_{n-1}%
;\mathbf{\lambda}\left(  \Delta\right)  \right)  -N_{\pi\circ\left(
n\leftrightarrow k+1\right)  }\left(  x_{1},...,x_{n-1};\mathbf{\lambda
}\left(  -\Delta\right)  \right)  =-1.\label{015}%
\end{equation}
(Here we denote by $\pi\circ\left(  n\leftrightarrow k+1\right)  $ the
permutation which is the composition of the transposition $n\leftrightarrow
k+1,$ followed by $\pi.$) Indeed, after the above reassignment and the
resolution of conflicts, the rod $\mathbf{\lambda}\left(  \delta\right)
_{\pi\left(  n\right)  }$ will be positioned at the point $y_{k}\left(
\mathbf{\lambda}\left(  \delta\right)  ,\pi\right)  .$ The relations $\left(
\ref{012}\right)  -\left(  \ref{014}\right)  $ then tell us that during the
$2\delta$-evolution the right endpoint of this rod will move from the positive
semiaxis to the negative one, thus adding one unit to the value $S_{\pi
\circ\left(  n\leftrightarrow k+1\right)  }\left(  x_{1},...,x_{n-1}%
;\mathbf{\lambda}\left(  \Delta\right)  \right)  .$

In the second case\textbf{\ }we have that the point $-\lambda\left(
\Delta\right)  _{\pi\left(  n\right)  }$\textbf{\ }is outside all clusters of
$R_{\pi}\left(  \mathbf{\lambda}\left(  \Delta\right)  \right)  ,$ while the
point $-\lambda\left(  -\Delta\right)  _{\pi\left(  n\right)  }$ is inside a
cluster of $R_{\pi}\left(  \mathbf{\lambda}\left(  -\Delta\right)  \right)  $.
That however can happen only if $\pi\left(  n\right)  =i,$ so the rod
$\lambda\left(  \delta\right)  _{\pi\left(  n\right)  }$ itself is varying
with $\delta,$ and in addition to it we have for some $k\neq n$ that
\[
-\lambda\left(  \delta\right)  _{\pi\left(  n\right)  }\left\{
\begin{array}
[c]{ll}%
<x_{k} & \text{ if }\delta>0,\\
=x_{k} & \text{ if }\delta=0,\\
>x_{k} & \text{ otherwice,}%
\end{array}
\right.
\]
while $\lambda_{\pi\left(  k\right)  }>-x_{k}.$ But then for a permutation
$\pi\circ\left(  n\leftrightarrow k\right)  $ we have that the variable rod
$\lambda_{i}\left(  \delta\right)  $ is assigned to the location $x_{k},$ is
not relocated due to the resolution of conflicts, so the endpoint $y_{k}$
equals $x_{k}+\lambda_{i}\left(  \delta\right)  ,$ and during the evolution
crosses the zero threshold from right to left; therefore again
\[
N_{\pi\circ\left(  n\leftrightarrow k\right)  }\left(  x_{1},...,x_{n-1}%
;\mathbf{\lambda}\left(  \Delta\right)  \right)  -N_{\pi\circ\left(
n\leftrightarrow k\right)  }\left(  x_{1},...,x_{n-1};\mathbf{\lambda}\left(
-\Delta\right)  \right)  =-1.
\]

The above construction corresponds to every permutation $\pi,$ satisfying
$\left(  \ref{010}\right)  ,$ another permutation, $\pi^{\prime}=\Phi\left(
\pi\right)  ,$ which satisfy $\left(  \ref{015}\right)  .$ We will be done if
we show that $\Phi$ is one to one. We prove this by constructing the inverse
of $\Phi.$

So let $\pi^{\prime}$ be such that
\[
N_{\pi^{\prime}}\left(  x_{1},...,x_{n-1};\mathbf{\lambda}\left(
\Delta\right)  \right)  -N_{\pi^{\prime}}\left(  x_{1},...,x_{n-1}%
;\mathbf{\lambda}\left(  -\Delta\right)  \right)  =-1.
\]
That means that the cardinality of the intersection $Y\left(  R_{\pi^{\prime}%
}\left(  \mathbf{\lambda}\left(  \Delta\right)  \right)  \right)  \cap\left[
-\mathbf{\lambda}\left(  \Delta\right)  _{\pi^{\prime}\left(  n\right)
},0\right]  $ increases by one during the $2\delta$-evolution: indeed, the
only other option, due to the Lemma \ref{formula}, is that the point
$-\lambda\left(  \Delta\right)  _{\pi^{\prime}\left(  n\right)  }$\ is inside
some cluster of $R_{\pi^{\prime}}\left(  \mathbf{\lambda}\left(
\Delta\right)  \right)  ,$ while the point $-\lambda\left(  -\Delta\right)
_{\pi^{\prime}\left(  n\right)  }$ is outside all clusters of $R_{\pi^{\prime
}}\left(  \mathbf{\lambda}\left(  -\Delta\right)  \right)  .$ But in this case
$N_{\pi^{\prime}}\left(  x_{1},...,x_{n-1};\mathbf{\lambda}\left(
\Delta\right)  \right)  =N_{\pi^{\prime}}\left(  x_{1},...,x_{n-1}%
;\mathbf{\lambda}\left(  -\Delta\right)  \right)  .$

In our situation we necessarily have that the intersection $Y\left(
R_{\pi^{\prime}}\left(  \mathbf{\lambda}\left(  \Delta\right)  \right)
\right)  \cap\left(  0,+\infty\right)  \neq\emptyset.$ Let $y_{k^{\prime}%
}\left(  \mathbf{\lambda}\left(  \Delta\right)  ,\pi^{\prime}\right)
<...<y_{r^{\prime}}\left(  \mathbf{\lambda}\left(  \Delta\right)  ,\pi
^{\prime}\right)  \ $are all the points of this intersection. The relation
$\left(  \ref{010}\right)  $ implies via $\left(  \ref{0110}\right)  $ that
the point $y_{k^{\prime}}\left(  \mathbf{\lambda}\left(  \delta\right)
,\pi^{\prime}\right)  $ moves from the positive semiaxis to the negative one
as $\delta$ passes the zero value:
\begin{equation}
y_{k^{\prime}}\left(  \mathbf{\lambda}\left(  \delta\right)  ,\pi^{\prime
}\right)  >0\text{ for }\delta>0,\label{021}%
\end{equation}
\begin{equation}
y_{k^{\prime}}\left(  \mathbf{\lambda}\left(  0\right)  ,\pi^{\prime}\right)
=0,\label{022}%
\end{equation}
\begin{equation}
y_{k^{\prime}}\left(  \mathbf{\lambda}\left(  \delta\right)  ,\pi^{\prime
}\right)  <0\text{ for }\delta<0.\label{023}%
\end{equation}
Two subcases are possible. The first one is when the point $z_{k^{\prime}}$ is
the head of the cluster, so $z_{k^{\prime}}=x_{k^{\prime}},$ and the rod
$\lambda_{\pi^{\prime}\left(  k^{\prime}\right)  }$ is the variable one. Then
for the permutation $\pi^{\prime\prime}=\pi^{\prime}\circ\left(
n\leftrightarrow k^{\prime}\right)  $ the following happen: the free rod
$\lambda_{\pi^{\prime\prime}\left(  n\right)  }$ is the variable one, and as
$\delta$ varies, the point $-\lambda_{\pi^{\prime\prime}\left(  n\right)  }$
moves from the location to the left of $x_{k^{\prime}}$ to the location to the
right to $x_{k^{\prime}}.$ Since the point $x_{k^{\prime}}$ is the head of the
cluster, we have therefore that
\begin{equation}
N_{\pi^{\prime\prime}}\left(  x_{1},...,x_{n-1};\mathbf{\lambda}\left(
\Delta\right)  \right)  -N_{\pi^{\prime\prime}}\left(  x_{1},...,x_{n-1}%
;\mathbf{\lambda}\left(  -\Delta\right)  \right)  =1.\label{024}%
\end{equation}
In the second subcase $z_{k^{\prime}}=y_{k^{\prime}-1}.$ Then the relations
$\left(  \ref{021}\right)  -\left(  \ref{023}\right)  $ mean that the point
$y_{k^{\prime}-1}\left(  \mathbf{\lambda}\left(  \delta\right)  ,\pi^{\prime
}\right)  $ is inside the segment $\left[  -\mathbf{\lambda}\left(
\delta\right)  _{\pi^{\prime}\left(  k^{\prime}\right)  },0\right]  $ for
$\delta=\Delta,$ and outside it for $\delta=-\Delta.$ So if we again assign
the free rod $\mathbf{\lambda}\left(  \delta\right)  _{\pi^{\prime}\left(
n\right)  }$ to the point $x_{k^{\prime}},$ making instead the rod
$\mathbf{\lambda}\left(  \delta\right)  _{\pi^{\prime}\left(  k^{\prime
}\right)  }$ to be free, then for the permutation $\pi^{\prime\prime}%
=\pi^{\prime}\circ\left(  n\leftrightarrow k^{\prime}\right)  =\Phi^{\prime
}\left(  \pi^{\prime}\right)  $ the relation $\left(  \ref{024}\right)  $
again holds.

The statement that $\Phi^{\prime}$ is inverse to $\Phi$ is straightforward.
\end{proof}

Below we will need a version of the above theorem, which follows. Let $T,L$ be
positive real, $L<T.$ Let again $n$ be an integer, and $\lambda_{1}%
<\lambda_{2}<...<\lambda_{n}$ be a fixed set of positive lengths of rods. Let
$-T<x_{1}<x_{2}<...<x_{n-1}<0$ be a set of $\left(  n-1\right)  $ left-ends.
We want to compute the number $\tilde{N}\left(  -T,x_{1},x_{2},...,x_{n-1}%
;L,\lambda_{1},\lambda_{2},...,\lambda_{n}\right)  ,$ which is defined as
follows. For any permutation $\pi$ of $n$ elements and for any $X\in\left(
\mathbb{-}T,0\right)  ,\,X\neq x_{1},x_{2},...,x_{n-1}$ we can consider the
configuration $\sigma_{n}\left(  -T,x_{1},...,x_{n-1};L,\lambda_{\pi\left(
1\right)  },...,\lambda_{\pi\left(  n-1\right)  }\right)  $ $\cup\sigma
_{1}\left(  X,\lambda_{\pi\left(  n\right)  }\right)  $ of rods, when the rod
$L$ is placed at $-T,$ the rods $l_{i}=\lambda_{\pi\left(  i\right)  }$ are
placed at $x_{i},i=1,...,n-1,$ while the free rod $l_{n}=\lambda_{\pi\left(
n\right)  }$ is placed at $X,$ $-T<X<0.$ Given $\pi,$ we count the number
$\tilde{N}_{\pi}\left(  -T,x_{1},x_{2},...,x_{n-1};L,\lambda_{1},\lambda
_{2},...,\lambda_{n}\right)  $ of different locations $X,$ such that the
corresponding r-configuration

\noindent$R\left[  \sigma_{n}\left(  -T,x_{1},...,x_{n-1};L,\lambda
_{\pi\left(  1\right)  },...,\lambda_{\pi\left(  n-1\right)  }\right)
\cup\sigma_{1}\left(  X,\lambda_{\pi\left(  n\right)  }\right)  \right]  $ has
a hit at zero, and moreover this hit is an $X$-hit. Then we define
\[
\tilde{N}\left(  -T,x_{1},x_{2},...,x_{n-1};L,\lambda_{1},\lambda
_{2},...,\lambda_{n}\right)  =\sum_{\pi\in\mathcal{S}_{n}}\tilde{N}_{\pi
}\left(  -T,x_{1},x_{2},...,x_{n-1};L,\lambda_{1},\lambda_{2},...,\lambda
_{n}\right)  .
\]

\begin{theorem}
\label{R} Suppose that
\begin{equation}
L+\lambda_{1}+\lambda_{2}+...+\lambda_{n}<T.\label{080}%
\end{equation}
Then $\tilde{N}\left(  -T,x_{1},x_{2},...,x_{n-1};L,\lambda_{1},\lambda
_{2},...,\lambda_{n}\right)  =n!$ for Lebesgue-almost every $x_{1}%
,...,x_{n-1}$ and $\lambda_{1},...,\lambda_{n}.$
\end{theorem}

The Theorem \ref{R} differs from the Theorem \ref{T6} by the presence of the
additional rod $L,$ which is placed at $-T,$ and by the restriction that all
points $X,x_{1},x_{2},...,x_{n-1}$ has to be within the segment $\left(
-T,0\right)  .$ In particular, the rod $L$ does not have to move under the
resolution of conflicts.

Note that without the restriction (\ref{080}) the statement of the theorem is
not valid, as it is easy to see.

\begin{proof}
Let the numbers $0<\varepsilon_{1}<...<\varepsilon_{n-1}$ be so small that the
sum $\varepsilon_{1}+...+\varepsilon_{n-1}$ is less than any of the numbers
$\left\vert \delta_{0}\left(  T-L\right)  +\delta_{1}\lambda_{1}+\delta
_{2}\lambda_{2}+...+\delta_{n}\lambda_{n}\right\vert ,$ where $\delta_{i}$ are
taking any of three values $-1,0,1,$ with the only restriction that not all of
them vanish simultaneously. Let us replace the configuration $x_{1}%
,x_{2},...,x_{n-1}$ by the configuration $x_{1}^{\prime},x_{2}^{\prime
},...,x_{n-1}^{\prime},$ where
\[
x_{i}^{\prime}=\left\{
\begin{array}
[c]{ll}%
L-T+\varepsilon_{i} & \text{ if }x_{i}<L-T,\\
x_{i} & \text{ otherwice.}%
\end{array}
\right.
\]
Let $k$ be the largest integer for which $x_{i}^{\prime}>x_{i}.$ (The meaning
of the configuration $x_{1}^{\prime},x_{2}^{\prime},...,x_{n-1}^{\prime}$ is
the following: were all $\varepsilon_{i}$ zeroes, it is the result of
resolving the first conflict, between the first rod $L$ and the rods
intersecting it, which rods have to be pushed to the right-hand end of $L.$ We
use positive $\varepsilon$-s in order to have all the point $x_{i}^{\prime}$
different.) By the previous theorem we know that $N\left(  x_{1}^{\prime
},x_{2}^{\prime},...,x_{n-1}^{\prime};\lambda_{1},...,\lambda_{n}\right)
=n!\,.$ Let the location $X$ and the permutation $\pi$ are such that the
corresponding r-configuration

\noindent$R\left[  \sigma_{n-1}\left(  x_{1}^{\prime},...,x_{n-1}^{\prime
};\lambda_{\pi\left(  1\right)  },...,\lambda_{\pi\left(  n-1\right)
}\right)  \cup\sigma_{1}\left(  X,\lambda_{\pi\left(  n\right)  }\right)
\right]  $ has an $X$-hit. The condition (\ref{080}) implies that the cluster
of the r-configuration

\noindent$R\left[  \sigma_{n-1}\left(  x_{1}^{\prime},...,x_{n-1}^{\prime
};\lambda_{\pi\left(  1\right)  },...,\lambda_{\pi\left(  n-1\right)
}\right)  \cup\sigma_{1}\left(  X,\lambda_{\pi\left(  n\right)  }\right)
\right]  ,$ rooted at $X,$ does not contain any of the points $z_{1}^{\prime
}=x_{1}^{\prime},z_{2}^{\prime},...,z_{k}^{\prime}$ (see (\ref{07}) for the
notation), so $X>L-T,$ and the r-configuration

\noindent$R\left[  \sigma_{n}\left(  -T,x_{1},...,x_{n-1};L,\lambda
_{\pi\left(  1\right)  },...,\lambda_{\pi\left(  n-1\right)  }\right)
\cup\sigma_{1}\left(  X,\lambda_{\pi\left(  n\right)  }\right)  \right]  $ has
an $X$-hit as well. Therefore

\noindent$\tilde{N}\left(  -T,x_{1},x_{2},...,x_{n-1};L,\lambda_{1}%
,\lambda_{2},...,\lambda_{n}\right)  \geq n!\,.$ On the other hand, if the r-configuration

\noindent$R\left[  \sigma_{n}\left(  -T,x_{1},...,x_{n-1};L,\lambda
_{\pi\left(  1\right)  },...,\lambda_{\pi\left(  n-1\right)  }\right)
\cup\sigma_{1}\left(  X,\lambda_{\pi\left(  n\right)  }\right)  \right]  $ has
an $X$-hit, then by the same reasoning $X$ has to be to the right of the
location $L-T,$ and moreover the cluster of this configuration, rooted at $X,$
does not contain any of the points $z_{1}=-T,z_{2}=-T+L,...,z_{k+1};$
therefore the r-configuration

\noindent$R\left[  \sigma_{n-1}\left(  x_{1}^{\prime},...,x_{n-1}^{\prime
};\lambda_{\pi\left(  1\right)  },...,\lambda_{\pi\left(  n-1\right)
}\right)  \cup\sigma_{1}\left(  X,\lambda_{\pi\left(  n\right)  }\right)
\right]  $ has an $X$-hit. Hence

\noindent$\tilde{N}\left(  -T,x_{1},x_{2},...,x_{n-1};L,\lambda_{1}%
,\lambda_{2},...,\lambda_{n}\right)  \leq n!\,,$ and the proof follows.
\end{proof}

\section{$10$ (ten) technical statements \label{ten}}

In this section we present several technical statements needed for the proof
of our main result. The first subsection deals with the regularity properties
of the NMP-s, while the second -- with the estimates on the averaging kernels
$q_{\ast,\ast}.$

\subsection{Regularity properties of Non-Linear Markov Process}

The first fact we will establish concerns the integral behavior of the output
rate (=input rate) of the NMP.

\begin{lemma}
\label{la} Let $\mu_{\nu,\lambda_{\nu}\left(  \cdot\right)  }\left(
\cdot\right)  $ be NMP, with $N\left(  \mu_{\nu,\lambda_{\nu}\left(
\cdot\right)  }\left(  t\right)  \right)  =N\left(  \nu\right)  =q.$ Then
there exists a time duration $T^{\prime}=T^{\prime}\left(  q\right)  $ and
$\varepsilon^{\prime}=\varepsilon^{\prime}\left(  q\right)  >0,$ such that for
all $T\geq T^{\prime}$ and all $s\geq0$
\begin{equation}
\int_{s}^{s+T}\lambda_{\nu}\left(  t\right)  \,dt<T\left(  1-\varepsilon
^{\prime}\right)  .\label{01}%
\end{equation}

\end{lemma}

Clearly, in order to see it, it is sufficient to prove the following statement:

\begin{proposition}
\label{damba} Let $\mu_{t}$ be GFP (i.e. time inhomogeneous Markov process
M(t)/GI/1, see Sect. 2) with arbitrary initial state $\mu_{0},$ corresponding
to the Poisson input with continuous rate $\lambda\left(  t\right)  ,$
$0\leq\lambda\left(  t\right)  \leq L.$ For any $N>0$ one can find
$\varepsilon\left(  N\right)  $ and $T\left(  N,L\right)  ,$ such that if for
some $T\geq T\left(  N,L\right)  $%
\begin{equation}
\int_{0}^{T}\lambda\left(  t\right)  dt\geq\left(  1-\varepsilon\left(
N\right)  \right)  T,\label{186}%
\end{equation}
then for some $t\leq T$%
\begin{equation}
N\left(  \mu_{t}\right)  \geq N.\label{183}%
\end{equation}

\end{proposition}

Indeed, if for any $T$ and $\varepsilon^{\prime}$ one can find $s$ such that
$\left(  \ref{01}\right)  $ is violated, then Proposition \ref{damba}, used
with $L=C\left(  \eta\right)  $ (see $\left(  \ref{13}\right)  $)$,$ would
imply that the function $N\left(  \mu_{\nu,\lambda_{\nu}\left(  \cdot\right)
}\left(  t\right)  \right)  $ is unbounded. This, however, contradicts to the
fact that it stays constant. The continuity of $\lambda_{\nu}\left(  t\right)
,$ which is a prerequisite needed to apply Proposition \ref{damba}, is a
consequence of Lemma \ref{Lip} below.

The proof of the Proposition \ref{damba} (which proposition is of course valid
even without continuity assumption) follows from the next two lemmas, preceded
by two definitions.

Let $\chi_{1},\chi_{2}$ be two measures on a segment $\left[  A,B\right]
\subset\mathbb{R}^{1}.$

\begin{definition}
We say that $\chi_{1}\prec\chi_{2},$ if for any monotone increasing function
$f$ on $\left[  A,B\right]  $ we have
\[
\int_{A}^{B}f~d\chi_{1}\leq\int_{A}^{B}f~d\chi_{2}.
\]

\end{definition}

\noindent(This is equivalent to saying that $\chi_{1}\left(  \left[
a,B\right]  \right)  \leq\chi_{2}\left(  \left[  a,B\right]  \right)  $ for
every $b\in\left[  A,B\right]  .$)

\begin{definition}
We say that $\chi_{1}\leq\chi_{2}$ if for every positive function $g$ on
$\left[  A,B\right]  $ we have
\[
\int_{A}^{B}g~d\chi_{1}\leq\int_{A}^{B}g~d\chi_{2}.
\]

\end{definition}

\noindent(This is equivalent to saying that $\chi_{1}\left(  \left[
a,b\right]  \right)  \leq\chi_{2}\left(  \left[  a,b\right]  \right)  $ for
every $a,b\in\left[  A,B\right]  .$) Note that the second relation holds for
probability measures only in case when $\chi_{1}=\chi_{2};$ in this paragraph
we are concerned, however, with arbitrary measures.

In what follows we will assume that the measures $\chi_{i}$ have densities
$\lambda_{i}$ with respect to the Lebesgue measure, which densities are
continuous and satisfy $0\leq\lambda_{i}\leq C,$ though much of what will be
said is true in more general situation.

\begin{lemma}
\label{coupling} Let us consider two GFP (non-stationary Markov processes)
$\mu_{t}^{1},\mu_{t}^{2},$ with the same distribution of the service time
$\eta,$ and with the Poisson input flows, defined by the two measures
$\chi_{1},\chi_{2}$ on the time interval $\left[  A,B\right]  .$ Suppose that
$\mu_{A}^{1}\left(  \mathbf{0}\right)  =1$ (that is, the first server is
initially idle), and that $\chi_{1}\prec\chi_{2}.$ Then
\[
N\left(  \mu_{B}^{1}\right)  \leq N\left(  \mu_{B}^{2}\right)  .
\]

\end{lemma}

\begin{proof}
The rough idea of the proof is the following: we will argue that the condition
$\chi_{1}\prec\chi_{2}$ enables us to represent the process $\mu_{t}^{2}$ as a
certain transformation of the process $\mu_{t}^{1},$ when the customers of the
second process are the same as these in the first process (i.e. they require
the same service times), but just come later, in addition to extra customers
which were not present in the first process. Clearly, in that case the queue
at the final moment has to be longer for the second process.

To make the above rigorous, we will use the coupling technique. Note first,
that if $\mathcal{P}$ is a Poisson random field on $\left[  A,B\right]  $ with
the rate function $\lambda\left(  t\right)  ,$ $t\in\left[  A,B\right]  ,$
$\omega\subset\left[  A,B\right]  $ is its realization, and $f:\left[
A,B\right]  \rightarrow\left[  A,B\right]  $ is a strictly increasing
continuous map, then the image set $f\left(  \omega\right)  $ is also a
realization of a Poisson random field $\mathcal{P}^{f},$ defined by the rate
function $\lambda^{f}\left(  t\right)  ,$ given by
\[
\lambda^{f}\left(  t\right)  =\left\{
\begin{array}
[c]{ll}%
\frac{\lambda\left(  f^{-1}\left(  t\right)  \right)  }{f^{\prime}\left(
f^{-1}\left(  t\right)  \right)  } & \text{ if }t\in f\left(  \left[
A,B\right]  \right)  ,\\
0 & \text{ otherwice.}%
\end{array}
\right.
\]
We claim now that if $\chi_{1}\prec\chi_{2},$ then there exists a map
$f:\left[  A,B\right]  \rightarrow\left[  A,B\right]  ,$ such that
\[
f\left(  x\right)  \geq x\text{ for all }x\in\left[  A,B\right]  ,
\]
\[
\lambda_{1}^{f}\left(  t\right)  \leq\lambda_{2}\left(  t\right)  .
\]

One way to construct such a function is the following. Let us extend the
function $\lambda_{2}\left(  t\right)  $ to the region $t\geq B$ by putting it
there to be equal to $\max\lambda_{1}.$ Consider now the family $\Phi$ of all
continuous functions $\varphi:\left[  A,B\right]  \rightarrow\lbrack
A,\infty),$ satisfying the properties:
\[
\varphi\left(  x\right)  \geq x;
\]
\[
\varphi\left(  y\right)  \geq\varphi\left(  x\right)  \text{ for }y>x;
\]
\[
\frac{\lambda_{1}\left(  x\right)  }{\varphi^{\prime}\left(  x\right)  }%
\leq\lambda_{2}\left(  \varphi\left(  x\right)  \right)  .
\]
For example, any shift $\varphi_{c}\left(  x\right)  =x+c$ is in $\Phi,$ once
$c>B-A,$ so $\Phi$ is non-empty. Note that if $\varphi_{1},\varphi_{2}\in
\Phi,$ then the function $\tilde{\varphi}\left(  x\right)  =\min\left\{
\varphi_{1}\left(  x\right)  ,\varphi_{2}\left(  x\right)  \right\}  $ is also
in $\Phi.$ Therefore the function
\[
f\left(  x\right)  =\inf_{\varphi\in\Phi}\left\{  \varphi\left(  x\right)
\right\}
\]
is in $\Phi$ as well, and it is easy to see that $f\left(  \left[  A,B\right]
\right)  \subset\left[  A,B\right]  .$ (In fact, if $\chi_{2}\left(  \left[
A,B\right]  \right)  =\chi_{1}\left(  \left[  A,B\right]  \right)  ,$ then $f$
is the only function in $\Phi$ with this property.)

Another, more intuitive way of defining $f$, is by constructing a coupling
between the measure $\chi_{1}$ and a suitable measure $\chi_{3}\succ\chi_{1},$
which is a \textquotedblleft part\textquotedblright\ of $\chi_{2},$ in the
sense that $\chi_{3}\leq\chi_{2}.$ {\small (The measure }$K${\small \ on
}$\left[  A,B\right]  \times\left[  A,B\right]  ${\small \ is called a
coupling between the measures }$\chi^{\prime}${\small \ and }$\chi
^{\prime\prime}${\small \ on }$\left[  A,B\right]  ${\small \ iff for every
subset }$C\subset\left[  A,B\right]  ${\small \ we have }$K\left(
C\times\left[  A,B\right]  \right)  =\chi^{\prime}\left(  C\right)
${\small \ and }$K\left(  \left[  A,B\right]  \times C\right)  =\chi
^{\prime\prime}\left(  C\right)  .${\small )}

The special coupling we need can be most easily constructed via discrete
approximations of the measures $\chi_{i}$ and subsequent limit procedure.
Since this construction is well-known due to extensive use in probability
theory of the Monge-Kantorovich-Rubinstein-Ornstein-Vasserstein distance, we
will give only a sketch of it. We replace the segment $\left[  A,B\right]  $
by the set $\left\{  1,2,...,n\right\}  \subset\mathbb{R}^{1}.$ The coupling
sought is then just a matrix $K\left(  i,j\right)  .$ We construct it by induction.

We define the first row $K\left(  1,\cdot\right)  =\left\{  K\left(
1,1\right)  ,K\left(  1,2\right)  ,...,K\left(  1,n\right)  \right\}  $ to be
the measure on $\left\{  1,2,...,n\right\}  $ with the following properties:

\textrm{i) }$K\left(  1,1\right)  +K\left(  1,2\right)  +...+K\left(
1,n\right)  =\chi_{1}\left(  1\right)  ;$

\textrm{ii) }$K\left(  1,\cdot\right)  \leq\chi_{2}\left(  \cdot\right)  ;$

\textrm{iii) }$K\left(  1,\cdot\right)  $ is the minimal (in the $\prec
$-sense) measure on $\left\{  1,2,...,n\right\}  ,$ satisfying \textrm{i)
}and\textrm{\ ii). }

\noindent(One can give more explicit definition of $K:$ we start by putting
$K\left(  1,1\right)  =\min\left\{  \chi_{1}\left(  1\right)  ,\chi_{2}\left(
1\right)  \right\}  .$ If it turns out that $\chi_{1}\left(  1\right)
=K\left(  1,1\right)  ,$ then we put $K\left(  1,j\right)  =0$ for $j>1.$
Otherwise we define $K\left(  1,2\right)  =\min\left\{  \chi_{1}\left(
1\right)  -K\left(  1,1\right)  ,\chi_{2}\left(  2\right)  \right\}  .$ If it
turns out that $\chi_{1}\left(  1\right)  =K\left(  1,1\right)  +K\left(
1,2\right)  ,$ then we put $K\left(  1,j\right)  =0$ for $j>2,$ otherwise
defining $K\left(  1,3\right)  =\min\left\{  \chi_{1}\left(  1\right)
-K\left(  1,1\right)  -K\left(  1,2\right)  ,\chi_{2}\left(  3\right)
\right\}  ,$ etc.)

Because of \textrm{ii), }the difference $\chi_{2}^{\left(  2\right)  }\left(
\cdot\right)  \equiv\chi_{2}\left(  \cdot\right)  -K\left(  1,\cdot\right)  $
is still a (positive) measure on $\left\{  1,2,...,n\right\}  ,$ and, because
of \textrm{iii),} $\chi_{2}^{\left(  2\right)  }\left(  \cdot\right)
\succ\chi_{1}\left(  \cdot\right)  $, \textit{considered as the measures on}
$\left\{  2,3,...,n\right\}  .$ Therefore we can repeat the preceding
construction, defining the second row, $K\left(  2,\cdot\right)  ,$ as the
measure on $\left\{  2,3,...,n\right\}  ,$ corresponding to $\chi_{2}^{\left(
2\right)  }\left(  \cdot\right)  $ and $\chi_{1}\Bigm|_{\left\{
2,3,...,n\right\}  }.$ Proceeding inductively, and considering the measures
$\chi_{2}^{\left(  k\right)  }\left(  \cdot\right)  $ on $\left\{
k,k+1,...,n\right\}  ,$ we obtain the coupling $K\left(  i,j\right)  $ sought.
By construction, $K\left(  i,j\right)  =0$ for $j<i.$ Finally, we define the
measure $\chi_{3}$ by
\[
\chi_{3}\left(  j\right)  =K\left(  1,j\right)  +...+K\left(  j,j\right)  .
\]
By construction, $\chi_{3}\succ\chi_{1},$ while $\chi_{3}\leq\chi_{2}.$

To perform the limiting procedure, we now define for every $n$ the atomic
measures $\chi_{i,n}$ with atoms at points $\left\{  A+k\frac{B-A}%
{n},~k=0,1,...,n-1\right\}  ,$ by
\[
\chi_{i,n}\left(  A+k\frac{B-A}{n}\right)  =\chi_{i}\left(  \left[
A+k\frac{B-A}{n},A+\left(  k+1\right)  \frac{B-A}{n}\right]  \right)  .
\]
We then construct the couplings $K_{n}$ between $\chi_{1,n}$ and $\chi_{3,n}$
(with $\chi_{1,n}\prec\chi_{3,n},\ \chi_{3,n}\leq\chi_{2,n}$) in the manner
described above. As $n\rightarrow\infty,$ the supports $S_{n}=\mathrm{supp}%
\left(  K_{n}\right)  \subset\left[  A,B\right]  \times\left[  A,B\right]  $
of the measures $K_{n}$ converge to the limiting curve $S $ in the square
$\left[  A,B\right]  \times\left[  A,B\right]  .$ This curve is the graph of
the above function $f:$%
\[
S=\left\{  \left(  x,y\right)  :y=f\left(  x\right)  \right\}  .
\]

Now we can return to the proof of the statement of the Lemma. Let $\omega
\in\left[  A,B\right]  $ be a configuration of the Poisson random field
$\mathcal{P}_{1},$ corresponding to the rate measure $\chi_{1}.$ Let $f$ be
the function defined above, $\chi_{3}=f_{\ast}\left(  \chi_{1}\right)  ,$ and
the measure $\zeta=\chi_{2}-\chi_{3}$ (that $\zeta$ is indeed a positive
measure follows from our construction). Let $\mathcal{\tilde{P}}$ be the
Poisson random field with the rate $\zeta,$ independent of $\mathcal{P}_{1},$
and $\tilde{\omega}$ be its configuration. Then the random set $\bar{\omega
}=f\left(  \omega\right)  \cup\tilde{\omega}$ has distribution of the
configuration of the Poisson random field $\mathcal{P}_{2},$ corresponding to
the rate measure $\chi_{2}.$ To specify the input flow and thus the coupling
sought, we have to specify the service times. So we assign to every point $x$
in $\omega$ the service time $\eta_{x},$ drawn independently from the
distribution of the random variable $\eta.$ To every point $y\in f\left(
\omega\right)  $ we assign the service time $\eta_{f^{-1}\left(  y\right)  },$
while to points $z\in\tilde{\omega}$ we assign independent realizations
$\eta_{z}$ of $\eta.$

Our statement now becomes almost evident. Consider a $\mathcal{P}_{1}%
$-customer, represented by a point $x\in\omega,$ whose service is not yet over
at the moment $B.$ This happens due to the service time $\eta_{x},$ needed for
him, as well as due to the service of the customers $x_{1},x_{2},...,$ who
came before $x.$ But then the corresponding $\mathcal{P}_{2}$-customer
$f\left(  x\right)  $ has the same service time $\eta_{x},$ but arrives later
than $x,$ as well as all the customers $f\left(  x_{1}\right)  ,f\left(
x_{2}\right)  ,...,$ who came before. In addition, there can be extra
$\mathcal{\tilde{P}}$-customers arrived before $f\left(  x\right)  .$ It is
evident therefore, that the service of the $\mathcal{P}_{2} $-customer
$f\left(  x\right)  $ will not be over at the moment $B.$ So the queue in the
second case is not shorter.
\end{proof}

In what follows we denote by $\kappa^{\left(  a\right)  }$ the measure on
$\left[  A,B\right]  $, having the constant density $\lambda\left(  t\right)
=a. $

\begin{lemma}
\label{calcul} Let the measure $\chi$ on $\left[  A,B\right]  $ has continuous
density $\lambda\left(  t\right)  ,$ $0\leq\lambda\left(  t\right)  \leq L,$
and satisfies the property:
\[
\chi\left(  \left[  A,B\right]  \right)  \geq\left(  1-\varepsilon\right)
\left(  B-A\right)  .
\]
Then there exists a segment $\left[  A,C\right]  \subset\left[  A,B\right]  $
of the length
\begin{equation}
C-A>\frac{\varepsilon}{L}\left(  B-A\right)  ,\label{305}%
\end{equation}
such that
\[
\chi\Bigm|_{\left[  A,C\right]  }\succ\kappa^{\left(  1-2\varepsilon\right)
}\Bigm|_{\left[  A,C\right]  }.
\]

\end{lemma}

\begin{proof}
Consider all the segments $\left[  c,d\right]  \subset\left[  A,B\right]  ,$
which have the property that
\begin{equation}
\chi\Bigm|_{\left[  c,d\right]  }\succ\kappa^{\left(  1-2\varepsilon\right)
}\Bigm|_{\left[  c,d\right]  }.\label{309}%
\end{equation}
Note that if $\left[  c_{1},d_{1}\right]  $ and $\left[  c_{2},d_{2}\right]  $
are two such segments, and $\left[  c_{1},d_{1}\right]  \cap\left[
c_{2},d_{2}\right]  \neq\varnothing,$ then the same holds for their union,
i.e.
\[
\chi\Bigm|_{\left[  c_{3},d_{3}\right]  }\succ\kappa^{\left(  1-2\varepsilon
\right)  }\Bigm|_{\left[  c_{3},d_{3}\right]  }%
\]
for $\left[  c_{3},d_{3}\right]  =\left[  c_{1},d_{1}\right]  \cup\left[
c_{2},d_{2}\right]  .$ Therefore the union of all segments with the property
$\left(  \ref{309}\right)  $ splits into the family of non-intersecting
maximal ones, $\left[  C_{1},D_{1}\right]  ,\left[  C_{2},D_{2}\right]  ,...,$
where we enumerate the segments according to their length, say, so $\left\vert
D_{1}-C_{1}\right\vert \geq$ $\left\vert D_{2}-C_{2}\right\vert \geq...$ .
Denote by $\left[  A,a\right]  $ the segment $\left[  C_{i},D_{i}\right]  $
among these maximal, which contains the point $A$. If in fact $A$ does not
belong to any maximal segment, we put $a=A$.

Let now $\left[  C_{k},D_{k}\right]  $ be one of maximal segments, which is
different from the segment $\left[  A,a\right]  .$ Then
\begin{equation}
\chi\left(  \left[  C_{k},D_{k}\right]  \right)  =\left(  1-2\varepsilon
\right)  \left(  D_{k}-C_{k}\right)  .\label{027}%
\end{equation}
Indeed, let us couple the measure $\kappa^{\left(  1-2\varepsilon\right)
}\Bigm|_{\left[  C_{k},D_{k}\right]  }$\ with (part of) the measure
$\chi\Bigm|_{\left[  C_{k},D_{k}\right]  }$\ -- that is with some measure
$\tilde{\chi}\Bigm|_{\left[  C_{k},D_{k}\right]  }\leq\chi\Bigm|_{\left[
C_{k},D_{k}\right]  },$\ see the proof of the Lemma \ref{coupling} above. This
is possible since $\kappa^{\left(  1-2\varepsilon\right)  }\Bigm|_{\left[
C_{k},D_{k}\right]  }\prec\chi\Bigm|_{\left[  C_{k},D_{k}\right]  }.$\ If
$\chi\left(  \left[  C_{k},D_{k}\right]  \right)  >\left(  1-2\varepsilon
\right)  \left(  D_{k}-C_{k}\right)  ,$\ then the difference $\Delta
\chi\Bigm|_{\left[  C_{k},D_{k}\right]  }=\chi\Bigm|_{\left[  C_{k}%
,D_{k}\right]  }-\tilde{\chi}\Bigm|_{\left[  C_{k},D_{k}\right]  }$\ is a
positive measure on $\left[  C_{k},D_{k}\right]  .$\ Therefore for some small
$\delta$\ the measure $\kappa^{\left(  1-2\varepsilon\right)  }\Bigm|_{\left[
C_{k}-\delta,C_{k}\right]  }$\ can be coupled with (part of) the measure
$\Delta\chi\Bigm|_{\left[  C_{k},D_{k}\right]  }.$\ That implies that
$\chi\Bigm|_{\left[  C_{k}-\delta,D_{k}\right]  }\succ\kappa^{\left(
1-2\varepsilon\right)  }\Bigm|_{\left[  C_{k}-\delta,D_{k}\right]  },$\ which
contradicts to the maximality of $\left[  C_{k},D_{k}\right]  $.

Note that for any point $x\in\left[  A,B\right]  ,$ which is outside all of
the segments $\left[  C_{i},D_{i}\right]  ,$ we have $\lambda\left(  x\right)
\leq1-2\varepsilon.$ Together with $\left(  \ref{027}\right)  $ it implies
that
\[
\chi\left(  \left[  a,B\right]  \right)  \leq\left(  1-2\varepsilon\right)
\left\vert B-a\right\vert .
\]
On the other hand,
\[
\chi\left(  \left[  A,B\right]  \right)  =\chi\left(  \left[  A,a\right]
\right)  +\chi\left(  \left[  a,B\right]  \right)  \geq\left(  1-\varepsilon
\right)  \left\vert B-A\right\vert ,
\]
so
\[
\left(  a-A\right)  L\geq\chi\left(  \left[  A,a\right]  \right)
\geq\varepsilon\left\vert B-A\right\vert ,
\]
and the proof follows.
\end{proof}

\medskip

\textbf{Proof of the Proposition \ref{damba}. }Let $N$ be fixed. As we know
(see relation $\left(  \ref{03}\right)  $ and the statement after it), there
is a value $c=c\left(  N\right)  <1,$ such that the homogeneous process with
the rate function $\lambda\equiv c$ has the invariant measure $\nu_{c}$ with
$N\left(  \nu_{c}\right)  =2N.$ We define $\varepsilon\left(  N\right)  $ by
\[
\varepsilon\left(  N\right)  =\frac{1-c\left(  N\right)  }{2}.
\]
If we start the process with $\lambda\equiv c$ in the state $\mathbf{0},$ then
it weakly converges with time to $\nu_{c}.$ In particular, $N\left(
\mu_{\mathbf{0},c}\left(  t\right)  \right)  \rightarrow2N,$ as $t\rightarrow
\infty.$ Define $\bar{T}=\bar{T}\left(  N\right)  $ to be the time duration,
after which the (monotone) function $N\left(  \mu_{\mathbf{0},c}\left(
t\right)  \right)  $ satisfies $N\left(  \mu_{\mathbf{0},c}\left(  t\right)
\right)  \geq N$ for all $t\geq\bar{T}.$ We want this value $\bar{T}$ to
appear as the lower bound of the length of the segment $\left[  A.C\right]  $,
obtained in the Lemma \ref{calcul}. That will be the case if the whole segment
$\left[  A,B\right]  $ will be of length
\[
T\left(  N,L\right)  =\frac{L\bar{T}\left(  N\right)  }{\varepsilon\left(
N\right)  },
\]
see $\left(  \ref{305}\right)  .$ Now, if $\left(  \ref{186}\right)  $ is
satisfied with the $T\left(  N,L\right)  $ and $\varepsilon\left(  N\right)  $
just chosen, then by Lemma \ref{calcul} the measure $\lambda\left(  t\right)
dt$ is bigger (in the $\succ$ sense) than the measure $c\left(  N\right)  dt$
on a segment $\left[  0,\bar{t}\right]  $ with $\bar{t}\geq\bar{T}.$ By Lemma
\ref{coupling} we have the relation $\left(  \ref{183}\right)  $ at $\bar{t}.$
$\blacksquare$

\medskip

Next we show that for every initial state $\nu$ of NMP there exists a time
moment after which the probability of observing the system $\mu_{\nu
,\lambda_{\nu}\left(  \cdot\right)  }\left(  t\right)  $ to be in the idle
state is uniformly positive.

\begin{lemma}
\label{pusto} Let $\mu_{\nu,\lambda_{\nu}\left(  \cdot\right)  }\left(
\cdot\right)  $ be NMP, with $N\left(  \mu_{\nu,\lambda_{\nu}\left(
\cdot\right)  }\left(  t\right)  \right)  =N\left(  \nu\right)  =q.$ Then
there exists a time moment $T=T\left(  \nu\right)  $ and $\varepsilon
=\varepsilon\left(  \nu\right)  >0,$ such that for all $t>T$%
\begin{equation}
\left\langle \omega=\mathbf{0}\right\rangle _{\mu_{\nu,\lambda_{\nu}\left(
\cdot\right)  }\left(  t\right)  }>\varepsilon.\label{302}%
\end{equation}

\end{lemma}

\begin{proof}
We first construct an auxiliary stationary ergodic Markov process, which in a
certain sense dominates our NMP from above. Namely, let $T=T^{\prime}\left(
q\right)  $ be the time duration proven to exist in Lemma \ref{la}, and
$\varepsilon$ be the corresponding quantity $\varepsilon^{\prime}\left(
q\right)  $. Our Markov process $M_{\nu,\varepsilon}\left(  t\right)  $ will
consist of the states of some auxiliary server $\mathfrak{S}$ at moments $t,$
with $M_{\nu,\varepsilon}\left(  0\right)  =\nu.$ The customers are arriving
to $\mathfrak{S}$ only at moments $kT,$ $k=1,2,...$ . Their numbers
$\mathfrak{N}_{k}$ are i.i.d., distributed as the number of Poisson flow of
customers with constant rate $\left(  1-\varepsilon\right)  $, arriving during
the time intervals $\left[  \left(  k-1\right)  T,kT\right]  .$ The service
times are described by our random variable $\eta.$ Since $\mathbb{E}\left(
\eta\right)  =1,$ while $\varepsilon>0,$ $M_{\nu,\varepsilon}$ is indeed
ergodic. In particular,
\begin{equation}
L\left(  \varepsilon,\delta\right)  \equiv\lim_{k\rightarrow\infty
}\left\langle \omega=\mathbf{0}\right\rangle _{M_{\nu,\varepsilon}\left(
kT-\delta\right)  }>0\label{184}%
\end{equation}
for any $\delta$ small enough.

Now we note that
\begin{equation}
\left\langle \omega=\mathbf{0}\right\rangle _{\mu_{\nu,\lambda_{\nu}\left(
\cdot\right)  }\left(  kT-\delta\right)  }\geq\mathbf{\Pr}\left(
\mathcal{E}_{k}\left(  T\right)  \right)  \left\langle \omega=\mathbf{0}%
\Bigm|\mathcal{E}_{k}\left(  T\right)  \right\rangle _{\mu_{\nu,\lambda_{\nu
}\left(  \cdot\right)  }\left(  kT-\delta\right)  },\label{142}%
\end{equation}
where the event
\[
\mathcal{E}_{k}\left(  T\right)  =\left\{
\begin{array}
[c]{l}%
\text{in the Poisson random flow, defined by the rate}\\
\lambda_{\nu}\left(  \cdot\right)  ,\text{ no customer arrives during the time
}\left[  \left(  k-1\right)  T,kT\right]  .
\end{array}
\right\}  .\text{ }%
\]
Since the rate $\lambda_{\nu}\left(  \cdot\right)  $ is bounded from above
uniformly in $\nu,$
\begin{equation}
\mathbf{\Pr}\left(  \mathcal{E}_{k}\left(  T\right)  \right)  \geq
\alpha\left(  T\right)  >0\label{143}%
\end{equation}
for some positive function $\alpha.$ On the other hand,
\begin{equation}
\left\langle \omega=\mathbf{0}\Bigm|\mathcal{E}_{k}\left(  T\right)
\right\rangle _{\mu_{\nu,\lambda_{\nu}\left(  \cdot\right)  }\left(
kT-\delta\right)  }\geq\left\langle \omega=\mathbf{0}\right\rangle
_{M_{\nu,\varepsilon}\left(  kT-\delta\right)  }.\label{144}%
\end{equation}
Indeed, since $\int_{s}^{s+T}\lambda_{\nu}\left(  t\right)  \,dt<T\left(
1-\varepsilon\right)  $ for all $s,$ the two processes $\mu_{\nu,\lambda_{\nu
}\left(  \cdot\right)  }$ and $M_{\nu,\varepsilon}$ can be coupled in such a
way that to each customer $C\left(  t,\bar{\eta}\right)  $ of the process
$\mu_{\nu,\lambda_{\nu}\left(  \cdot\right)  },$ who arrives at the moment
$t,$ $\left(  l-1\right)  T<t\leq lT$ and uses the server for time $\bar{\eta
},$ it corresponds a customer $C^{\prime}\left(  \mathfrak{t}\left(  t\right)
,\bar{\eta}\right)  $ of the process $M_{\nu,\varepsilon},$ who arrives at a
later moment
\[
\mathfrak{t}\left(  t\right)  =\left(  \left\lfloor \frac{t}{T}\right\rfloor
+1\right)  T=lT
\]
and needs the server for the same time duration $\bar{\eta}.$ Hence the queue
at every moment $t<kT$ of the process $M_{\nu,\varepsilon}$ is not shorter
than the one for the process $\mu_{\nu,\lambda_{\nu}^{\left(  k-1\right)
T}\left(  \cdot\right)  },$ where
\[
\lambda_{\nu}^{\left(  k-1\right)  T}\left(  t\right)  =\left\{
\begin{array}
[c]{ll}%
\lambda_{\nu}\left(  t\right)  & \text{ if }t\leq\left(  k-1\right)  T,\\
0 & \text{ otherwice.}%
\end{array}
\right.  .
\]
From $\left(  \ref{143}\right)  ,$ $\left(  \ref{144}\right)  $ we infer that
for all $k\geq k_{0}$%
\[
\left\langle \omega=\mathbf{0}\right\rangle _{\mu_{\nu,\lambda_{\nu}\left(
\cdot\right)  }\left(  kT-\delta\right)  }\geq\frac{1}{2}\alpha\left(
T\right)  L\left(  \varepsilon.\delta\right)  ,
\]
where $k_{0}$ is the smallest index $k,$ for which $\left\langle
\omega=\mathbf{0}\right\rangle _{M_{\nu,\varepsilon}\left(  kT-\delta\right)
}\geq\frac{1}{2}L\left(  \varepsilon,\delta\right)  .$ (Of course, the value
$k_{0}=k_{0}\left(  \nu\right)  $ does depend on the initial state $\nu.)$

Thus far we got the desired result $\left(  \ref{302}\right)  $ only for
values $t\in\left[  kT-\delta,kT\right]  ,$ $k\geq k_{0}\left(  \nu\right)  .$
To take care of the other values of $t$-s we should just make a change of
variables and to start our process from all the different states $\mu
_{\nu,\lambda_{\nu}\left(  \cdot\right)  }\left(  t\right)  ,$ $t\in\left[
0,T\right]  .$ Be the convergence in $\left(  \ref{184}\right)  $ uniform in
$\nu, $ that would be the end of the proof. However this is not the case. Yet,
this does not create any problem, since the family of states $\left\{
\mu_{\nu,\lambda_{\nu}\left(  \cdot\right)  }\left(  t\right)  ,t\in\left[
0,T\right]  \right\}  ,$ which will be taken for initial states of the process
$M_{\ast,\varepsilon},$ is compact, and on it the convergence in $\left(
\ref{184}\right)  $ is uniform indeed.
\end{proof}

We finish this subsection with a statement about the regularity of the exit flow.

\begin{lemma}
\label{Lip} Let the function $p\left(  t\right)  $ satisfies the strong
Lipschitz condition $\left(  \ref{02}\right)  $: for some $C$%
\begin{equation}
\left\vert p\left(  t+\Delta t\right)  -p\left(  t\right)  \right\vert \leq
Cp\left(  t\right)  \Delta t.\label{028}%
\end{equation}
Let $\nu$ be some initial state of GFP, and $\lambda\left(  \cdot\right)  $ be
arbitrary rate function of arriving customers. (In particular one can take
$\lambda=\lambda_{\nu},$ thus getting NMP.) Then the rate function $b\left(
t\right)  $ of the corresponding exit flow is Lipschitz, with Lipschitz
constant independent of $\nu,$ $\lambda$.
\end{lemma}

\begin{proof}
Let $t$ be fixed. The informal idea of the proof is the following: consider
the elementary event $\varpi$, which contributes to the output rate $b\left(
t\right)  ,$ and let $\mathfrak{c}$ be the customer, who leaves the server at
the moment $t,$ being under the service for the time duration $\mathfrak{t.}$
Then the elementary event $\varpi^{\prime},$ obtained from $\varpi$ by
enlarging the service time of $\mathfrak{c}$ from $\mathfrak{t}$ to
$\mathfrak{t}+\Delta t$, contributes to $b\left(  t+\Delta t\right)  .$ So we
can get the desired result by comparing the probabilities of $\varpi$ and
$\varpi^{\prime}.$

This correspondence, however, does not \textquotedblleft
cover\textquotedblright\ all the events, contributing to $b\left(  t+\Delta
t\right)  .$ Namely, the elementary events not covered by the above
correspondence, are precisely those events $\varpi^{\prime\prime}$, for which
the customer $\mathfrak{c,}$ which is the last customer who has started
his/her service before the moment $t$, is different from the customer
$\mathfrak{c}^{\prime\prime},$ whose service terminates at $t+\Delta t.$ This,
though, means that the service time of $\mathfrak{c}^{\prime\prime}$ was less
than $\Delta t.$

To implement this idea, let us write the measure of our process
\[
\exp\left\{  -I_{\lambda}\left(  T\right)  \right\}  \sum_{n=0}^{\infty}%
\frac{1}{n!}\left[  \prod_{i=1}^{n}\lambda\left(  x_{i}\right)  \,dx_{i}%
\prod_{i=1}^{n}p\left(  l_{i}\right)  \,dl_{i}\right]
\]
on a segment $\left[  0,T\right]  ,$ with any $T>t+\Delta t$ (compare with
$\left(  \ref{129}\right)  $), as
\[
d\Pi\left(  \varpi\right)  =\frac{1}{Z_{T}}\pi\left(  \varpi\right)  d\varpi,
\]
with $Z_{T}=\exp\left\{  I_{\lambda}\left(  T\right)  \right\}  $ and
$\pi\left(  \varpi\right)  =\prod_{i=1}^{n}\left[  \lambda\left(
x_{i}\right)  p\left(  l_{i}\right)  \right]  ,$ so
\[
b\left(  t\right)  =\int_{B\left(  t\right)  }d\Pi\left(  \varpi\right)  .
\]
Here $B\left(  t\right)  $ is the manifold of all elementary events $\varpi,$
which have a moment of service termination to happen at $t.$ (We are treating
here the case when the initial state $\nu$ is concentrated on the
configuration $\omega=\mathbf{0};$ the general case is totally similar.)

Let us split the rate
\[
b\left(  t+\Delta t\right)  =\int_{B\left(  t+\Delta t\right)  }d\Pi\left(
\varpi\right)
\]
into two parts. The first one, $b^{\prime}\left(  t+\Delta t\right)  ,$ is
given by
\[
b^{\prime}\left(  t+\Delta t\right)  =\int_{B^{\prime}\left(  t+\Delta
t\right)  }d\Pi\left(  \varpi\right)  ,
\]
where $B^{\prime}\left(  t+\Delta t\right)  \subset B\left(  t+\Delta
t\right)  $ is the image of the manifold $B\left(  t\right)  $ under the map
$\varpi\rightsquigarrow\varpi^{\prime},$ defined in the first paragraph of the
present proof. Therefore
\[
\left\vert b^{\prime}\left(  t+\Delta t\right)  -b\left(  t\right)
\right\vert =\left\vert \int_{B\left(  t\right)  }\left(  \frac{p\left(
\mathfrak{t}+\Delta t\right)  }{p\left(  \mathfrak{t}\right)  }-1\right)
d\Pi\left(  \varpi\right)  \right\vert ,
\]
where $\mathfrak{t}=\mathfrak{t}\left(  \varpi\right)  $ is the service time
of the customer $\mathfrak{c}$ in $\varpi$, whose service terminates at the
moment $t.$ From $\left(  \ref{028}\right)  $ it follows that
\[
\left\vert b^{\prime}\left(  t+\Delta t\right)  -b\left(  t\right)
\right\vert \leq C\Delta t.
\]
The remaining part $b^{\prime\prime}\left(  t+\Delta t\right)  $ of the rate
$b\left(  t+\Delta t\right)  ,$%
\[
b^{\prime\prime}\left(  t+\Delta t\right)  =\int_{B^{\prime\prime}\left(
t+\Delta t\right)  }d\Pi\left(  \varpi\right)  ,
\]
corresponds to the event when some customer $\mathfrak{c}$ finishes his
service during the time period $\left[  t,t+\Delta t\right]  ,$ while the next
customer needs time at most $\Delta t$ to be served. Such an event has
probability below $b\left(  t\right)  \left(  \Delta t\right)  ^{2},$ so the
proof follows.
\end{proof}

\subsection{Estimates on the averaging kernels}

Here we will estimate the densities $q_{\lambda,x}\left(  t\right)  ,$
entering into the relation $b\left(  x\right)  =\left[  \lambda\ast
q_{\lambda,x}\right]  \left(  x\right)  .$

\begin{lemma}
\label{lowerb} The family $q_{\lambda,y}\left(  t\right)  $ is weakly
continuous in $y,$ for every $\lambda.$ Also
\begin{equation}
q_{\lambda,y}\left(  t\right)  \geq p\left(  t\right)  \mathbf{\Pr}\left\{
\text{server is idle at the moment }y-t\right\}  .\label{301}%
\end{equation}

\end{lemma}

\begin{lemma}
\label{oc}
\begin{equation}
q_{\lambda,y}\left(  t\right)  \leq\sum_{n=1}^{\infty}p^{\ast n}\left(
t\right)  \mathbf{\Pr}\left\{  N_{t}^{\lambda,y}\geq n-1\right\}
\equiv\mathcal{Q}_{\lambda,y}\left(  t\right)  ,\label{152}%
\end{equation}
where $N_{t}^{\lambda,y}$ is the random number of $\lambda$-Poisson points in
the segment $\left[  y-t,y\right]  .$

In particular, there exists a constant $\tilde{C}=\tilde{C}\left(  p\right)
,$ such that for $t$ and all $\lambda,y$%
\begin{equation}
q_{\lambda,y}\left(  t\right)  \leq\tilde{C}.\label{157}%
\end{equation}

\end{lemma}

\begin{proof}
Both the relations $\left(  \ref{301}\right)  $ and $\left(  \ref{152}\right)
$ follow easily from the definition $\left(  \ref{005}\right)  .$ To see
$\left(  \ref{301}\right)  $ we note that, evidently, $c\left(  u,t\right)
\geq p\left(  t\right)  .$ To get $\left(  \ref{152}\right)  ,$ we split the
event $\mathcal{C}\left(  u,t\right)  ,$ entering $\left(  \ref{008}\right)
$, into the sum of events $\mathcal{C}_{n}\left(  u,t\right)  ,$ $n\geq1,$
where $\mathcal{C}_{n}\left(  u,t\right)  $ consists of all these outcomes
when between the moment $u$ of arrival of the first customer and the arrival
of the customer who terminates during $\left[  u+t,u+t+h\right]  $ precisely
$n-2$ other customers came. (The event $\mathcal{C}_{1}\left(  u,t\right)  $
consists from the outcomes when the first customer himself terminates during
$\left[  u+t,u+t+h\right]  .$) But if $\mathcal{C}_{n}\left(  u,t\right)  $
holds, then two independent events have to happen:
\[
\eta_{1}+...+\eta_{n}\in\left[  t,t+h\right]  ,
\]
and
\[
N_{t}^{\lambda,u}\geq n-1,
\]
which imply $\left(  \ref{152}\right)  .$

To see $\left(  \ref{157}\right)  ,$ we use a rough form of $\left(
\ref{152}\right)  :$%
\begin{equation}
q_{\lambda,y}\left(  t\right)  \leq\sum_{n=1}^{\infty}p^{\ast n}\left(
t\right)  .\label{8158}%
\end{equation}
Let $A=\sup_{t}p\left(  t\right)  .$ Then it is immediate from $\left(
\ref{8158}\right)  $ that for all $t\leq C$%
\[
q_{\lambda,y}\left(  t\right)  \leq A\left(  1+\sum_{n=1}^{\infty}\mathbf{\Pr
}\left\{  \eta_{1}+...+\eta_{n}\leq C\right\}  \right)  ,
\]
where $\eta_{i}$ are i.i.d. random variables, distributed as $\eta.$ But the
probabilities $\mathbf{\Pr}\left\{  \eta_{1}+...+\eta_{n}\leq C\right\}  $
decay exponentially in $n,$ so the series converges for every $t$. It is a
classical result of the renewal theory, that the sum $\left(  \ref{8158}%
\right)  $ goes to a finite limit as $t\rightarrow\infty,$ see e.g. the
relation (1.17) of Chapter XI in \cite{F}. That proves $\left(  \ref{157}%
\right)  .$

Continuity of $q_{\lambda,y}$ in $y$ follows from the definition in a
straightforward way.
\end{proof}

We will need the compactness estimate on the distributions $q_{\lambda
,y}\left(  t\right)  .$ We will obtain them using the estimate $\left(
\ref{152}\right)  .$ As the following statement shows, the estimate $\left(
\ref{152}\right)  $ is rather rough; we believe that all the moments of the
distribution $q_{\lambda,y}\left(  t\right)  $ of order less than $1+\delta$
are finite.

\begin{lemma}
\label{Qu} Suppose that $\lambda$ is such that for some $T^{\prime}$ and
$\varepsilon^{\prime}>0$ and for all $T\geq T^{\prime}$ and $s\geq0$
\begin{equation}
\int_{s}^{s+T}\lambda\left(  t\right)  \,dt<T\left(  1-\varepsilon^{\prime
}\right) \label{51}%
\end{equation}
(see $\left(  \ref{01}\right)  $). Then for any $b<\frac{\delta}{2}$
\begin{equation}
\int_{0}^{\infty}t^{b}q_{\lambda,y}\left(  t\right)  \,dt<C\left(
\lambda,b\right)  <\infty,\label{303}%
\end{equation}
where $C\left(  \lambda,b\right)  $ depends on $\lambda$ only via $T^{\prime}$
and $\varepsilon^{\prime}.$
\end{lemma}

\textbf{Proof of Lemma \ref{Qu}. }We are going to use the simple estimate: for
every random variable $\zeta$ and every $\varkappa>0$
\begin{equation}
\tilde{Q}\left(  T\right)  \equiv\mathbf{\Pr}\left\{  \zeta>T\right\}  \leq
T^{-\varkappa}\mathbb{E}\left(  \left|  \zeta\right|  ^{\varkappa}\right)
.\label{151}%
\end{equation}

We also will need an estimate on $\int_{A}^{\infty}t^{a}\tilde{q}\left(
t\right)  \,dt,$ $a<\varkappa,$ where $\tilde{q}$ is the density of $\zeta.$
We have:
\begin{align}
\int_{A}^{\infty}t^{a}\tilde{q}\left(  t\right)  \,dt  & =-\int_{A}^{\infty
}t^{a}\,d\left(  \tilde{Q}\left(  t\right)  \right) \label{168}\\
& =A^{a}\tilde{Q}\left(  A\right)  +a\int_{A}^{\infty}t^{a-1}\tilde{Q}\left(
t\right)  \,dt.\nonumber
\end{align}

To apply $\left(  \ref{151}\right)  $ to $\left(  \ref{152}\right)  $ we will
use the Dharmadhikari-Yogdeo estimate (see, e.g. \cite{P}, p.79): if $\xi_{i}$
are independent centered random variables, then
\begin{equation}
\mathbb{E}\left(  \left\vert \xi_{1}+...+\xi_{n}\right\vert ^{2+\delta
}\right)  \leq Rn^{\delta/2}\sum_{1}^{n}\mathbb{E}\left(  \left\vert \xi
_{i}\right\vert ^{2+\delta}\right)  .\label{159}%
\end{equation}
Here $R=R\left(  \delta\right)  $ is some universal constant.

Introducing $\xi_{i}=\eta_{i}-1$ (see $\left(  \ref{153}\right)  $), and using
$\left(  \ref{151}\right)  $ with $\varkappa=2+\delta$ and $\left(
\ref{159}\right)  ,$ we have (see $\left(  \ref{009}\right)  $) \textbf{\ }
\begin{align}
Q_{n}\left(  t\right)   & \equiv\mathbf{\Pr}\left\{  \eta_{1}+...+\eta
_{n}>t\right\}  =\mathbf{\Pr}\left\{  \xi_{1}+...+\xi_{n}>t-n\right\}
\label{156}\\
& \leq RM_{\delta}\left(  t-n\right)  ^{-\left(  2+\delta\right)  }%
n^{1+\delta/2}.\nonumber
\end{align}

To proceed, we use $\left(  \ref{152}\right)  $ to write
\begin{align}
\int_{0}^{\infty}t^{b}q_{\lambda,y}\left(  t\right)  \,dt  & \leq\int
_{0}^{\infty}t^{b}\mathcal{Q}_{\lambda,y}\left(  t\right)  \,dt\label{155}\\
& =\sum_{n=1}^{\infty}\left[  \int_{0}^{\infty}t^{b}p^{\ast n}\left(
t\right)  \mathbf{\Pr}\left\{  N_{t}^{\lambda,y}\geq n-1\right\}  \,dt\right]
.\nonumber
\end{align}

Note that due to $\left(  \ref{51}\right)  $ there exists an $\alpha
>0\mathbb{\ }$such that $\mathbb{E}\left(  N_{t}^{\lambda,y}\right)
\leq\left(  1-\alpha\right)  t$ once $t$ is large enough, uniformly in $y$.
The first step is to estimate every summand by
\begin{align}
& \int_{0}^{\infty}t^{b}p^{\ast n}\left(  t\right)  \mathbf{\Pr}\left\{
N_{t}^{\lambda,y}\geq n-1\right\}  \,dt\label{300}\\
& \leq\int_{0}^{n\left(  1+\frac{\alpha}{2}\right)  }t^{b}p^{\ast n}\left(
t\right)  \mathbf{\Pr}\left\{  N_{t}^{\lambda,y}\geq n-1\right\}
\,dt+\int_{n\left(  1+\frac{\alpha}{2}\right)  }^{\infty}t^{b}p^{\ast
n}\left(  t\right)  \,dt.\nonumber
\end{align}
Now, using $\left(  \ref{168}\right)  $ and $\left(  \ref{156}\right)  ,$ we
have for the second term in $\left(  \ref{300}\right)  :$
\begin{align*}
\int_{n\left(  1+\frac{\alpha}{2}\right)  }^{\infty}t^{b}p^{\ast n}\left(
t\right)  \,dt  & \leq\left[  n\left(  1+\frac{\alpha}{2}\right)  \right]
^{b}RM_{\delta}\left(  \frac{\alpha}{2}n\right)  ^{-\left(  2+\delta\right)
}n^{1+\delta/2}\\
& +bRM_{\delta}n^{1+\delta/2}\int_{n\left(  1+\frac{\alpha}{2}\right)
}^{\infty}t^{b-1}\left(  t-n\right)  ^{-\left(  2+\delta\right)  }\,dt\\
& \leq Cn^{b-1-\delta/2},
\end{align*}
where $C=C\left(  \alpha,\delta,M_{\delta}\right)  .$

The first term in $\left(  \ref{300}\right)  $ is negligible. To see that, we
first observe:

\begin{lemma}
Let $0<\nu<1,$ and $N_{n}^{\nu}$ be a Poisson random variable:
\[
\mathbf{\Pr}\left\{  N_{n}^{\nu}=k\right\}  =e^{-\nu n}\frac{\left(  \nu
n\right)  ^{k}}{k!}.
\]
Then
\[
\mathbf{\Pr}\left\{  N_{n}^{\nu}\geq n\right\}  \leq\frac{1}{1-\nu}%
e^{-\frac{\left(  1-\nu\right)  ^{2}}{2}n},
\]
provided $n$ is large enough.
\end{lemma}

\begin{proof}
Note first of all, that if $\chi>0$ and $n>\chi,$ then
\[
e^{-\chi}\sum_{k\geq n}\frac{\chi^{k}}{k!}\leq e^{-\chi}\frac{\chi^{n}}%
{n!}\sum_{k\geq0}\left(  \frac{\chi}{n+1}\right)  ^{k}=e^{-\chi}\frac{\chi
^{n}}{n!}\frac{1}{1-\frac{\chi}{n+1}}.
\]
In our case we thus have
\[
\sum_{k\geq n}\mathbf{\Pr}\left\{  N_{n}^{\nu}=k\right\}  \leq e^{-\nu n}%
\frac{\left(  \nu n\right)  ^{n}}{n!}\frac{1}{1-\nu}.
\]
By Stirling, for $n$ large
\begin{align*}
\sum_{k\geq n}\mathbf{\Pr}\left\{  N_{n}^{\nu}=k\right\}   & \leq\frac
{1}{1-\nu}e^{-\nu n}\frac{\nu^{n}n^{n}}{n^{n}e^{-n}}\\
& =\frac{1}{1-\nu}e^{\left(  1-\nu+\ln\nu\right)  n}\\
& \leq\frac{1}{1-\nu}e^{-\frac{\left(  1-\nu\right)  ^{2}}{2}n}.
\end{align*}

\end{proof}

To estimate the integral $\int_{0}^{n\left(  1+\frac{\alpha}{2}\right)  }%
t^{b}p^{\ast n}\left(  t\right)  \mathbf{\Pr}\left\{  N_{t}^{\lambda,y}\geq
n-1\right\}  \,dt$ we note, that in the range of $t\in\left[  0,n\left(
1+\frac{\alpha}{2}\right)  \right]  $ we have
\[
\mathbf{\Pr}\left\{  N_{t}^{\lambda,y}\geq n-1\right\}  \leq\mathbf{\Pr
}\left\{  N_{n\left(  1+\frac{\alpha}{2}\right)  }^{\lambda,y}\geq
n-1\right\}  .
\]
We apply to the r.h.s. the last lemma, with $\nu=\frac{n-1}{n\left(
1+\frac{\alpha}{2}\right)  }.$ Therefore, for all $n$ large enough and
uniformly in $y$
\begin{align*}
& \int_{0}^{n\left(  1+\frac{\alpha}{2}\right)  }t^{b}p^{\ast n}\left(
t\right)  \mathbf{\Pr}\left\{  N_{t}^{\lambda,y}\geq n-1\right\}  \,dt\\
& \leq\frac{2}{\alpha}e^{-\frac{\alpha^{2}}{8}n}\int_{0}^{n\left(
1+\frac{\alpha}{2}\right)  }t^{b}p^{\ast n}\left(  t\right)  \,dt\\
& \leq\frac{2}{\alpha}e^{-\frac{\alpha^{2}}{8}n}\left[  n\left(
1+\frac{\alpha}{2}\right)  \right]  ^{b}.
\end{align*}
Hence, the moment $\int_{0}^{\infty}t^{b}q_{\lambda,y}\left(  t\right)  \,dt$
is finite as soon as the series $\sum_{n}n^{b-1-\delta/2}$ converges, which
happens when $b<\frac{\delta}{2}.$ That proves Lemma \ref{Qu}. $\blacksquare$

\section{The self-averaging relation: general case \label{gencase}}

Here we derive a formula, expressing the function $b\left(  \cdot\right)
=A\left(  \mu,\lambda\left(  \cdot\right)  \right)  $ in terms of the
functions $\lambda\left(  \cdot\right)  $, $p\left(  \cdot\right)  $ and the
initial state $\mu$\textbf{\ }of our non-stationary (GFP) Markov process. This
will be the needed self-averaging relation (\ref{134}). We remind the reader
that $\mu$ is a probability measure on the set of pairs $\left\{  \left(
n,\tau\right)  \right\}  \cup\mathbf{0.}$

\begin{theorem}
Let $N\left(  \mu\right)  =q,$ and the rate function $\lambda\left(
\cdot\right)  $ satisfies the conclusions of the Lemma \ref{la}:
\begin{equation}
\int_{s}^{s+T}\lambda\left(  t\right)  \,dt<T\left(  1-\varepsilon^{\prime
}\right)  \text{ }\,\text{for all }T\geq T^{\prime}>0,\text{ all }s\geq0\text{
and some }\varepsilon^{\prime}>0.\label{185}%
\end{equation}
Then there exists the family of probability densities $q_{\lambda,\mu
,x}\left(  \cdot\right)  ,$ $x>0,$ and the functionals $\varepsilon
_{\lambda,\mu}\left(  x\right)  $ and $Q_{\lambda,\mu}\left(  x\right)  ,$
such that
\begin{equation}
b\left(  x\right)  =\left(  1-\varepsilon_{\lambda,\mu}\left(  x\right)
\right)  \left[  \lambda\ast q_{\lambda,\mu,x}\right]  \left(  x\right)
+\varepsilon_{\lambda,\mu}\left(  x\right)  Q_{\lambda,\mu}\left(  x\right)
.\label{103}%
\end{equation}

Moreover,
\begin{equation}
\varepsilon_{\lambda,\mu}\left(  x\right)  \rightarrow0\text{ as }%
x\rightarrow\infty,\label{102}%
\end{equation}
thought not necessarily uniformly in $\mu,$ while $Q_{\lambda,\mu}\left(
x\right)  \leq C,$ uniformly\textbf{\ }in $\lambda,\mu$ and $x.$
\end{theorem}

\begin{proof}
We start by defining the functional $\varepsilon_{\lambda,\mu}\left(
x\right)  .$ Note that the description of the realization of our process up to
the moment $x$ consists of the following data:

$i)$ the initial configuration $\left(  n,\tau\right)  ,$ drawn from the
distribution $\mu$ (with $n$ to be the number of customers in the system
before time zero, $\tau$ being the time the first one already spent in the server);

$ii)$ the random set $0<x_{1}<...<x_{m}<x$ (with random number $m$ of points),
which is a realization of the Poisson random field defined by the rate
function $\lambda$ (restricted to the segment $\left[  0,x\right]  $),
independent of $\left(  n,\tau\right)  $ (the arrival moments of the
customers, which come after time moment zero);

$iii)$ one realization $\eta_{1}$ of the conditional random variable
$\eta_{\tau}\equiv\left(  \eta-\tau\Bigm|\eta>\tau\right)  $ and $n+m-1$
independent realizations $\eta_{k},k=2,...,n+m$ of the random variable $\eta$
(service times for the customers).

\noindent We denote by $\mathbb{P}_{\mu\otimes\lambda\otimes\eta}$ the
corresponding (product) distribution. The difference $1-\varepsilon
_{\lambda,\mu}\left(  x\right)  $ is by definition just the $\mathbb{P}%
_{\mu\otimes\lambda\otimes\eta}$-probability of the event
\begin{equation}
\sum_{1}^{n+m}\eta_{k}<x.\label{140}%
\end{equation}
(If $n=0,$ then by definition we put $\tau=0;$ we put also $\sum_{1}^{0}%
\equiv0.$)

The meaning of the decomposition (\ref{103}) can be explained now: the first
term corresponds to the exit flow computed over those realizations where the
relation (\ref{140}) holds, while the second term represents the rest of the flow.

Let us prove $\left(  \ref{102}\right)  ,$ that is that
\[
\mathbf{\Pr}\left\{  \sum_{1}^{n+m}\eta_{k}>x\right\}  \rightarrow0\,\text{as
}x\rightarrow\infty.
\]
To do this, we introduce two independent random variables:
\[
S_{\mu}=\sum_{1}^{n}\eta_{k},\;S_{\lambda}=\sum_{n+1}^{n+m}\eta_{k}.
\]
Then for every $\alpha\in\left(  0,1\right)  $ we have
\begin{align*}
\mathbf{\Pr}\left\{  \sum_{1}^{n+m}\eta_{k}>x\right\}   & =\mathbf{\Pr
}\left\{  S_{\mu}+S_{\lambda}>x\right\} \\
& \leq\mathbf{\Pr}\left\{  S_{\mu}>\alpha x\right\}  +\mathbf{\Pr}\left\{
S_{\lambda}>\left(  1-\alpha\right)  x\right\}  .
\end{align*}
Indeed, if $S_{\mu}+S_{\lambda}>x,$ then either $S_{\mu}>\alpha x,$ or else
$S_{\lambda}>\left(  1-\alpha\right)  x.$ Since $S_{\mu}$ is a random
variable, the probability $\mathbf{\Pr}\left\{  S_{\mu}>\alpha x\right\}  $
goes to zero for every $\alpha$ positive, as $x\rightarrow\infty,$ though not
necessarily uniformly in $\mu.$ For the second term we have
\begin{align*}
& \mathbf{\Pr}\left\{  S_{\lambda}>\left(  1-\alpha\right)  x\right\} \\
& =\sum_{m=1}^{\infty}\left(  \int_{\left(  1-\alpha\right)  x}^{\infty
}p^{\ast m}\left(  t\right)  \,dt\right)  \mathbf{\Pr}\left\{  N^{\lambda
,x}=m\right\}  .
\end{align*}
Here $N^{\lambda,x}$ is the random number of points of the $\lambda$-Poisson
field in $\left[  0,x\right]  .$ Note that $\mathbb{E}\left(  N^{\lambda
,x}\right)  <x\left(  1-\varepsilon^{\prime}\right)  $ once $x>T^{\prime}.$
Therefore we can apply the same argument which was used in the proof of Lemma
\ref{Qu} when showing that the integral $\int_{T}^{\infty}\mathcal{Q}%
_{\lambda,y}\left(  t\right)  \rightarrow0$ as $T\rightarrow\infty,$ see
$\left(  \ref{152}\right)  .$ It implies that $\mathbf{\Pr}\left\{
S_{\lambda}>\left(  1-\alpha\right)  x\right\}  \rightarrow0$ once $\alpha$ is
small enough, uniformly in $\lambda,$ satisfying $\left(  \ref{185}\right)  $.
That establishes $\left(  \ref{102}\right)  .$

Next we define the distributions $q_{\lambda,\mu,x}.$ They are constructed
from the random field of the rods $\left\{  \eta_{k},k=1,...,n+m\right\}  ,$
defined above, placed at locations $\left\{  \underset{n}{\underbrace
{0,...,0}},x_{1},...,x_{m}\right\}  ,$ via the procedure of resolution of
conflicts, defined in the previous section. To do it we first introduce the
rate $b_{L}\left(  x\right)  $ to be the exit rate of the conditional service
process under the conditions that
\begin{equation}
\sum_{1}^{n}\eta_{k}=L,\;\;\sum_{n+1}^{n+m}\eta_{k}<x-L.\label{141}%
\end{equation}
We claim that for some probability distributions $q_{\lambda,L,x}$ we have
\[
b_{L}\left(  x\right)  =\left[  \lambda\ast q_{\lambda,L,x}\right]  \left(
x\right)  .
\]
The distribution $q_{\lambda,\mu,x}$ is then obtained by integration:
\[
q_{\lambda,\mu,x}=\int q_{\lambda,L,x}\mathbb{P}_{\mu\otimes\lambda\otimes
\eta}\left(  \sum_{1}^{n}\eta_{k}\in dL\right)  .
\]
(The random variable $\sum_{1}^{n}\eta_{k}$ is of course independent of the
Poisson $\lambda$-field.) The output rate $b_{L}\left(  x\right)  $
corresponds to the situation when we have customers arriving at the moments
$0,x_{1},...,x_{m},$ which have serving times $L,\eta_{n+1},...,\eta_{n+m},$
and which satisfy the relation
\[
L+\sum_{n+1}^{n+m}\eta_{k}<x.
\]
So we have to repeat the construction of the Section 5 in the present
situation. Few steps require some comments. The transition from the relation
(\ref{20}) to (\ref{23}) uses the fact that for any $s$ the measure
$\prod_{i=1}^{s}p\left(  l_{i}\right)  \,dl_{i}$ is invariant under the
coordinate permutations $S_{s}$ in $\mathbb{R}^{s}.$ But the same $S_{m}$
symmetry evidently holds for the conditional distribution of the random vector
$\left\{  \left(  \eta_{k},k=n+1,...,n+m\right)  \Bigm|\sum_{n+1}^{n+m}%
\eta_{k}<x-L\right\}  ,$ since both the unconditional distribution and the
distribution of the condition are $S_{m}$-invariant. The next crucial step was
the relation (\ref{26}), stating that the functions $q_{\lambda,y}$ are
probability distributions. It was based on the Theorem \ref{T6}. The situation
at hand is somewhat more delicate, since the rods we are dealing now with, are
of two kinds: the first one has a non-random length $L,$ produced by the
initial state $\mu,$ while others are situated at the Poissonian locations
$\left\{  x_{i}\right\}  ,$ defined by the rate function $\lambda.$ However,
under condition $\sum_{n+1}^{n+m}\eta_{k}<x-L$ the needed combinatorial
statement (about the quantity $m!$) still holds, and is the content of the
Theorem \ref{R}. These remarks allow one to carry over the construction of the
Section \ref{self-av}, and so to establish the existence of the probability
densities $q_{\lambda,L,x},$ and thus also $q_{\lambda,\mu,x}.$ The upper and
lower estimates on $q_{\lambda,\mu,x}$ are obtained in the same way as were
the estimates for $q_{\lambda,x}$ in the preceding section.

The function $Q_{\lambda,\mu}\left(  x\right)  $ is the rate of exit flow of
our process, conditioned by the event
\[
\sum_{1}^{n+m}\eta_{k}\geq x.
\]
The boundedness of the $Q_{\lambda,\mu}\left(  x\right)  $ follows from the
following property of the service time distribution $p\left(  x\right)  $: for
every $x,\tau,\,x>\tau>0,$ $1>t>0$
\begin{equation}
\frac{p\left(  x\right)  }{p\left(  x+t\right)  }\leq C^{\prime}%
,\;\;\frac{p\left(  x-\tau\Bigm|\eta>\tau\right)  }{p\left(  x-\tau
+t\Bigm|\eta>\tau\right)  }\leq C^{\prime}.\label{104}%
\end{equation}
The relation (\ref{104}) follows easily from the condition $\left(
\ref{02}\right)  ,$ with $C^{\prime}=C^{\prime}\left(  C\right)  .$ To explain
the boundedness, consider the elementary event
\[
\left(  n,\tau\right)  \times\left\{  x_{1},...,x_{m}:0<x_{1}<...<x_{m}%
<x\right\}  \times\left\{  \eta_{1},...,\eta_{n+m}\right\}  ,
\]
which contributes to the output flow inside the segment $\left[  x,x+\Delta
x\right]  ,$ which flow is accounted by the second term of (\ref{103}). That
means that our rod configuration produces after resolution of conflicts a hit
inside $\left[  x,x+\Delta x\right]  ,$ and also that
\begin{equation}
\sum_{1}^{n+m}\eta_{k}>x.\label{105}%
\end{equation}
In the notation of the Section 6 it means that after resolution of conflicts
the endpoint $y_{k}$ of some (shifted) rod fits within $\left[  x,x+\Delta
x\right]  ,$ for some $k\in\left\{  1,...,n+m\right\}  .$ Let $\bar{k}$ be the
smallest such index. But then the elementary events
\[
\left(  n,\tau\right)  \times\left\{  x_{1},...,x_{m}:0<x_{1}<...<x_{m}%
<x\right\}  \times\left\{  \eta_{1},...,\eta_{\bar{k}-1},\eta_{\bar{k}}%
+t,\eta_{\bar{k}+1},...,\eta_{n+m}\right\}  ,
\]
with any $t\in\left(  \Delta x,1\right)  ,$ do not contribute to the output
flow inside the segment $\left[  x,x+\Delta x\right]  ,$ while still
satisfying (\ref{105}). Therefore, due to (\ref{104}), the probability that
the customer would finish his service during the period $\left[  x,x+\Delta
x\right]  ,$ is of the order of $\Delta x,$ and, moreover,
\[
Q_{\lambda,\mu}\left(  x\right)  \leq\frac{1}{C^{\prime}}\text{.}%
\]

\end{proof}

\medskip

Let now $\mathfrak{M\in}\mathcal{M}\left(  \mathcal{M}_{q}\left(
\Omega\right)  \right)  $ be some invariant measure of the dynamical system
$\left(  \ref{160}\right)  .$ Then $\mathfrak{M}$-almost every state
$\tilde{\mu}_{0}\in\mathcal{M}_{q}\left(  \Omega\right)  $ belongs to the
family $\left\{  \tilde{\mu}_{t}:-\infty<t<+\infty\right\}  ,$ such that for
all $\tau>0,$ all $t$
\[
\mathcal{T}_{\tau}\left(  \tilde{\mu}_{t}\right)  =\tilde{\mu}_{t+\tau}.
\]
Let us fix one such family $\left\{  \tilde{\mu}_{t}\right\}  $. Then the
function $\lambda\left(  t\right)  $, $-\infty<t<+\infty,$ which for every
$-\infty<\tau<+\infty$ satisfies on $[\tau,+\infty)$ the equation
\[
\lambda\left(  \cdot\right)  =A\left(  \tilde{\mu}_{\tau},\lambda\left(
\cdot\right)  ,\tau\right)  ,
\]
is well defined. Then, according to the equation $\left(  \ref{103}\right)  ,$
for every $\tau,$ $-\infty<\tau<+\infty,$ and for all $x\geq\tau$%

\begin{equation}
\lambda\left(  x\right)  =\left(  1-\varepsilon_{\lambda,\tilde{\mu}_{\tau}%
}\left(  x\right)  \right)  \left[  \lambda\ast q_{\lambda,\tilde{\mu}_{\tau
},x}\right]  \left(  x\right)  +\varepsilon_{\lambda,\tilde{\mu}_{\tau}%
}\left(  x\right)  Q_{\lambda,\tilde{\mu}_{\tau}}\left(  x\right)
.\label{122}%
\end{equation}
One would like to pass here to the limit $\tau\rightarrow-\infty.$ According
to $\left(  \ref{102}\right)  ,$ for every $x$ we have $\varepsilon
_{\lambda,\tilde{\mu}_{\tau}}\left(  x\right)  \rightarrow0$ as $\tau
\rightarrow-\infty.$ Moreover, it is not difficult to show that in the same
limit $q_{\lambda,\tilde{\mu}_{\tau},x}\left(  \cdot\right)  \rightarrow
q_{\lambda,x}\left(  \cdot\right)  .$ So the following equation holds for
$\lambda:$%
\begin{equation}
\lambda\left(  x\right)  =\left[  \lambda\ast q_{\lambda,x}\right]  \left(
x\right)  ,\;\;-\infty<x<+\infty.\label{123}%
\end{equation}
By the methods developed below one can show that every bounded solution of
$\left(  \ref{123}\right)  $ is a constant. Since, however, we are proving a
stronger statement, that the dynamical system $\mathcal{T}_{\tau}$ has one
fixed point on each $\mathcal{M}_{q}\left(  \Omega\right)  ,$ which is,
moreover, globally attractive, we will not provide the details.

\section{Self-averaging $\Longrightarrow$ relaxation: a warm-up
\label{warm-up}}

Before presenting the general proof that self-averaging implies relaxation, we
consider the following simpler system: we have infinitely many servers, with
service time $\eta,$ distributed according to the probability density $p.$ As
the customer comes, he chooses any free server, and is served, leaving the
system afterwards. The inflow is Poissonian, given by the rate function
$f\left(  x\right)  .$ If we impose the condition that the customers are
coming at the rate they are living the system, we get the non-linear Markov
process. The self-averaging relation (\ref{34}) in such a case simplifies to
\[
b\left(  x\right)  =\left[  f\ast p\right]  \left(  x\right)  .
\]

\begin{lemma}
Let $p\left(  x\right)  $ be some probability density with support, belonging
to $\mathbb{R}^{+},$ satisfying the conditions of Section \ref{notation}. Let
$f$ be a positive bounded function on $\mathbb{R}^{1}.$ Suppose that
\begin{equation}
f\ast p\left(  x\right)  =f\left(  x\right)  \text{ for all }x\geq
0.\label{100}%
\end{equation}
Then $f\left(  x\right)  \rightarrow c$ as $x\rightarrow\infty,$ for some
$c>0.$
\end{lemma}

\begin{proof}
Our statement follows easily from the well known results of the renewal
theory. Let us introduce the function
\[
\varphi\left(  x\right)  =\left\{
\begin{array}
[c]{ll}%
f\left(  x\right)  & \text{ if }x<0,\\
0 & \text{ if }x>0.
\end{array}
\right.
\]
Then the values of the function $f$ for $x\geq0$ can be recovered from the
function $\varphi$ and the fact that it satisfies $\left(  \ref{100}\right)
$. Indeed, iterating $\left(  \ref{100}\right)  $, we find:
\begin{equation}
f\left(  x\right)  =\left\{
\begin{array}
[c]{ll}%
\left[  \varphi\ast p\right]  \left(  x\right)  +\left[  \;\left(  \varphi\ast
p\Bigm|_{\left\{  x\geq0\right\}  }\right)  \ast\left(  \sum_{n=1}^{\infty
}p^{\ast n}\right)  \right]  \left(  x\right)  & \text{ for }x>0,\\
\varphi\left(  x\right)  & \text{ for }x<0.
\end{array}
\right. \label{101}%
\end{equation}

The function
\[
s\left(  x\right)  =\sum_{n=1}^{\infty}p^{\ast n}\left(  x\right)
\]
is a key object of the renewal theory; in particular, it is known that it goes
to a positive limit as $x\rightarrow\infty;$ see for example the relation
(1.17) of Chapter XI in \cite{F}. From that and the relation $\left(
\ref{101}\right)  $ our claim follows.
\end{proof}

In fact, the constant $c$ can be computed: if $m$ is the mean value of $\eta,$
then
\[
\lim_{x\rightarrow\infty}f\left(  x\right)  =\frac{1}{m}\int_{0}^{+\infty
}\left[  \varphi\ast p\right]  \left(  x\right)  \;dx.
\]
This fact, as well as the renewal relation itself -- $s\left(  x\right)
\rightarrow const$ as $x\rightarrow\infty$ -- is a consequence of the
statement that the probability density $p^{\ast n}\left(  x\right)  $ of the
sum of i.i.d. random variables, $\eta_{1}+...+\eta_{n},$ is well approximated,
due to the local limit theorem, by the Gaussian distribution
\[
\frac{1}{\sqrt{2\pi nv}}e^{-\left(  x-nm\right)  ^{2}/2nv},
\]
where $v$ is the variance of $\eta.$ It becomes very flat as $n$ increases.

\section{Self-averaging $\Longrightarrow$ relaxation: probabilistic proof?
\label{no-go}}

As was already said in Section 2, any function $\lambda$, defined for $x<0,$
and vanishing for $x<-T,$ can be uniquely extended to $x\geq0$ in such a way
that the relation
\[
A\left(  \mathbf{0},\lambda\left(  \cdot\right)  ,-T\right)  =b\left(
\cdot\right)
\]
holds with $b\left(  x\right)  =\lambda\left(  x\right)  $ for $x\geq0.$
Therefore for every $x\geq0$ we have
\begin{equation}
\lambda\left(  x\right)  =\left[  \lambda\ast q_{\lambda,x}\right]  \left(
x\right)  ,\label{35}%
\end{equation}
where $q_{\lambda,x}\left(  \cdot\right)  $ is a probability density supported
by the semiaxis $\left\{  y\geq0\right\}  ,$ and depending on $\lambda$ via
its restriction $\lambda\Bigm|_{\left\{  y\leq x\right\}  }.$ Our goal is to
show that (\ref{35}) implies that $\lambda\left(  x\right)  $ relaxes to some
constant $c$ as $x\rightarrow\infty.$

Since the distributions $q_{\lambda,x}$ depend on $\lambda\left(
\cdot\right)  $ in a very complicated way, we have to treat a more general
statement. Suppose a family of probability densities $q_{x}\left(
\cdot\right)  ,$ supported by the semiaxis $\left\{  y\geq0\right\}  ,$ is
given, where $x\geq0.$ Let $f\left(  x\right)  $ be a non-negative function,
defined on $\mathbb{R}^{1},$ such that
\begin{align}
f\left(  x\right)   & \leq C\text{ for }x<0,\nonumber\\
f\left(  x\right)   & =\left[  f\ast q_{x}\right]  \left(  x\right)  \text{
for }x\geq0.\label{36}%
\end{align}
One would like to show that
\begin{equation}
\lim_{x\rightarrow\infty}f\left(  x\right)  =c,\label{43}%
\end{equation}
for some $c\geq0.$ That will imply the relaxation needed.

Motivated by the analysis of the previous section, we will study the equation
(\ref{36}) by considering the corresponding non-stationary renewal process and
the resulting inhomogeneous Markov random walk. Unfortunately, the relation
(\ref{43}) does not follow from (\ref{36}) in general, and the reasons are
probabilistic! Before explaining it let us \textquotedblleft
solve\textquotedblright\ (\ref{36}).

So, let the family $\left\{  q_{x},x\geq0\right\}  $ be given; we solve
(\ref{36}) for $f,$ given its restriction $f\Bigm|_{\left\{  x<0\right\}  }.$
We do this in close analogy with the previous section, see (\ref{101}). We
put
\[
f_{0}\left(  x\right)  =\left\{
\begin{array}
[c]{ll}%
f\left(  x\right)  & \text{ for }x<0\\
0 & \text{ for }x\geq0.
\end{array}
\right.
\]
We define
\[
f_{n+1}\left(  x\right)  =\left\{
\begin{array}
[c]{ll}%
f\left(  x\right)  & \text{ for }x<0,\\
\left[  f_{n}\ast q_{x}\right]  \left(  x\right)  & \text{ for }x\geq0.
\end{array}
\right.
\]
Then for every $x$ the sequence $f_{n}\left(  x\right)  $ is increasing, and
the function $f\left(  x\right)  =\lim_{n\rightarrow\infty}f_{n}\left(
x\right)  $ solves (\ref{36}).

We can rewrite the function $f$ in a different way. We define
\[
g_{1}\left(  x\right)  =\left\{
\begin{array}
[c]{ll}%
\left[  f_{0}\ast q_{x}\right]  \left(  x\right)  & \text{ for }x\geq0,\\
0 & \text{ for }x<0,
\end{array}
\right.
\]
\begin{equation}
g_{n+1}\left(  x\right)  =\left[  g_{n}\ast q_{x}\right]  \left(  x\right)
.\label{37}%
\end{equation}
Then for $x\geq0$ we have
\[
f\left(  x\right)  =\sum_{n\geq1}g_{n}\left(  x\right)  .
\]
Now we will write the formula for $g_{n}$ in terms of convolution. Here
instead of the density $p^{\ast n}$ of the sum $S_{n}$ of i.i.d. random
variables $\eta_{1}+...+\eta_{n}$ of the previous Section we have to consider
the distribution $p_{x}^{\left(  n\right)  }$ of the inhomogeneous Markov
walker $\bar{S}_{n,x},$ defined as follows. Remember that at each point
$x\in\mathbb{R}^{+}$ we have the probability density $q_{x}.$ So when our
walker after some steps happen to arrive to the location $x,$ then the next
move is to the location $x+y,$ with the increment $y>0$ distributed with the
density $q_{x}\left(  y\right)  .$ The random variable $\bar{S}_{n,x}$ is now
defined as the position of the above described Markov walker after $n$ steps,
with the initial position $x.$ With these notation we have, by (\ref{37}):
\[
g_{n+1}\left(  x\right)  =\left[  g_{1}\ast p_{x}^{\left(  n\right)  }\right]
\left(  x\right)  .
\]

Summarizing, we have for $x>0:$
\[
f\left(  x\right)  =g_{1}\left(  x\right)  +\sum_{n\geq2}\left[  g_{1}\ast
p_{x}^{\left(  n\right)  }\right]  \left(  x\right)  ,
\]
in analogy with $\left(  \ref{101}\right)  .$ As in the previous section, the
Local Limit Theorem for the Markov chain $\bar{S}_{n,x}$ would imply the
relaxation $\left(  \ref{43}\right)  .$

We have to note, however, that the relation between the validity of the LLT
for our Markov chain and the validity of the relation (\ref{43}) is more
complicated. Namely, LLT for $\bar{S}_{n,x}$ might fail, while the relaxation
(\ref{43}) might still remain valid -- or else fail as well! First of all, let
us explain that even the Central Limit Theorem for $\bar{S}_{n,x}$ might not
hold, notwithstanding the family $q_{x}\left(  \cdot\right)  $ to have very
nice compactness properties. To give one example, consider the family of
probability densities $u_{x}\left(  t\right)  ,$ $x\in\mathbb{R}^{1},$ where
all $u_{x}\left(  \cdot\right)  $ have for their support the segment $\left[
0,1\right]  ,$ and satisfy there $0<c<u_{x}\left(  t\right)  <C<\infty,$
uniformly in $x$ and $t.$ We define now
\[
q_{x}\left(  t\right)  =u_{x}\left(  t-\left(  1-\left\{  x\right\}  \right)
\right)  ,
\]
where $\left\{  \cdot\right\}  $ stays for the fractional part. Then all
$q_{x}\left(  \cdot\right)  $-s have their supports within the segment
$\left[  0,2\right]  .$ But the random variables $\bar{S}_{n,x}$ do not have
CLT behavior! Indeed, the random variable $\bar{S}_{n,x}$ is localized in the
segment $\left[  \left\lfloor x\right\rfloor +n,\left\lfloor x\right\rfloor
+n+1\right]  ,$ where $\left\lfloor \cdot\right\rfloor $ denotes the integer
part.\textbf{\ }So the variance of $\bar{S}_{n,x}$ remains bounded in $n~$!
Nevertheless, for this example it can be shown that the relation (\ref{43})
still holds, and that involves certain statement of the type of
Perron-Frobenius theorem for our Markov chain. Further modification of this
example, when
\[
q_{x}\left(  t\right)  =u_{x}\left(  t-\left(  2-\left\{  x\right\}  \right)
\right)  ,
\]
results in the Markov chain with two classes, and in this case both the CLT
and the relation (\ref{43}) might fail.

We conjecture here that the LLT theorem for the sums $\bar{S}_{n,x}$ holds, if
the family $q_{x}\left(  \cdot\right)  $ of transition densities has the
following additional property:

\begin{itemize}
\item For some $k,K,$ $0<k<K<\infty,$
\begin{equation}
k\leq\frac{q_{x_{1}}\left(  t\right)  }{q_{x_{2}}\left(  t\right)  }\leq
K,\label{45}%
\end{equation}
provided at least one of the values $q_{x_{i}}\left(  t\right)  $ is positive.
\end{itemize}

The condition (\ref{45}) is reminiscent of the \textit{positivity of
ergodicity coefficient}\textbf{\ }condition, introduced by Dobrushin \cite{D1}
in his study of the limit theorems for the non-stationary Markov chains.

In what follows we will take another road, and we get the relaxation property
by analytic methods, which seems in our case to be simpler. But we still use
probability theory, though not the CLT. It would be interesting to obtain the
desired result by proving the corresponding limit theorem.

\section{Self-averaging $\Longrightarrow$ relaxation: finite range case
\label{frc}}

In this section we prove the relaxation for the solution of the equation
(\ref{36}) in the finite range case.

The reader of the paper of course understands that for our initial problem we
have to consider the infinite range case. We think nevertheless that the
finite range case is of independent interest, and moreover, it holds in a much
more general setting than the infinite-range case. This is why we devote to it
the present Section. Its content is not used in what follows.

\begin{theorem}
\label{finite range} Suppose that
\begin{align*}
0  & \leq f\left(  x\right)  \leq C\text{ and is continuous for }x<0,\\
f\left(  x\right)   & =\left[  f\ast q_{x}\right]  \left(  x\right)  \text{
for }x\geq0,
\end{align*}
while the following three conditions on the family $q_{x}$ hold:

$i)$ for some $T$ and all $x$%
\[
\int_{0}^{T}q_{x}\left(  t\right)  \,dt=1,\ \ q_{x}\left(  t\right)  =0\text{
for }t>T,
\]

$ii)$ the family $q_{x}\left(  \cdot\right)  $ depends on $x$ continuously,

$iii)$ for $0\leq t\leq T$%
\begin{equation}
C\geq q_{x}\left(  t\right)  \geq\kappa\left(  t\right)  >0,\label{016}%
\end{equation}
with continuous positive $\kappa\left(  t\right)  .$

Then the limit exists:
\[
\lim_{x\rightarrow\infty}f\left(  x\right)  =c\geq0.
\]

\end{theorem}

The property $\left(  \ref{016}\right)  $ holds for the NMP, as follows from
the relations $\left(  \ref{301}\right)  ,$ $\left(  \ref{302}\right)  $ and
$\left(  \ref{157}\right)  .$

\begin{proof}
$i)$ We know that the function $f$ is continuous and bounded, $0\leq f\leq C.
$ So if there exists a value $X$ such that $f$ is monotone for $x\geq X,$ then
the function $f$ has to be constant for $x\geq X+T,$ and we are done. So we
are left with the case when the function $f$ has infinitely many points of
local maxima and local minima, which go to $\infty.$

$ii)$ Given a local maximum, $x_{0},$ we will construct now a sequence $x_{i}
$ of local maximums, $i=0,-1,-2,...,-n=-n\left(  f,x_{0}\right)  $ such that

\begin{itemize}
\item $x_{0}>x_{-1}>x_{-2}>...,$

\item $x_{i}-x_{i-1}<2T,\;\;x_{i}-x_{i-2}\geq T$ for all $i,$

\item $0<x_{-n}<2T,$

\item $f\left(  x_{i-1}\right)  \geq f\left(  x\right)  $ for any $x_{i-1}\leq
x, $ and $f\left(  x_{i-1}\right)  >f\left(  x_{i}\right)  ,$

\item for every $x\in\left[  x_{i-1},x_{i}-T\right]  $ we have $f\left(
x\right)  \geq f\left(  x_{i}\right)  $ (of course if the segment is
non-empty, i.e. $x_{i}-x_{i-1}>T$).
\end{itemize}

The construction is the following. Let $x_{0}$ be some point of local maxima.
Since
\[
f\left(  x_{0}\right)  =\int_{0}^{T}f\left(  x_{0}-t\right)  q_{x_{0}}\left(
t\right)  \,dt,
\]
we have $f\left(  x_{0}\right)  <F\left(  x_{0}\right)  \equiv\sup\left\{
f\left(  x\right)  :x\in\left[  x_{0}-T,x_{0}\right]  \right\}  ,$ unless $f$
is a constant on $\left[  x_{0}-T,x_{0}\right]  ,$ in which case we are done. Let

\noindent$y=\inf\left\{  x\in\left[  x_{0}-T,x_{0}\right]  :f\left(  x\right)
=F\left(  x_{0}\right)  \right\}  .$ If $y>x_{0}-T,$ or if $y=x_{0}-T$ and is
a local maximum, we define $x_{-1}=y.$ In the opposite case we have that the
point $x_{0}-T$ is not a local maximum of the function $f$ on the segment
$\left[  x_{0}-2T,x_{0}-T\right]  .$ We then consider two cases.

In the first one we suppose that the function $f$ on the segment $\left[
x_{0}-2T,x_{0}-T\right]  $ takes values below $\bar{F}=\frac{f\left(
x_{0}\right)  +f\left(  x_{0}-T\right)  }{2}.$ Let $\left[  y,x_{0}-T\right]
\subset\left[  x_{0}-2T,x_{0}-T\right]  $ be the largest segment for which the
inequality $f\left(  x\right)  \geq\bar{F}$ holds for every $x\in\left[
y,x_{0}-T\right]  .$ We define $x_{-1}$ to be the leftmost point of maximum of
$f$ in $\left[  y,x_{0}-T\right]  .$

In the opposite case we consider the set $S=\left\{  x\in\left[
x_{0}-2T,x_{0}-T\right]  :f\left(  x\right)  \geq f\left(  x_{0}-T\right)
\right\}  .$ It contains other points besides $x_{0}-T.$ However, it can not
contain all the segment $\left[  x_{0}-2T,x_{0}-T\right]  .$ Since $f\ $is not
a constant on $\left[  x_{0}-2T,x_{0}-T\right]  ,$ $\,\sup_{S}f>f\left(
x_{0}-T\right)  .$ Let $z\in\left(  x_{0}-2T,x_{0}-T\right)  $ be such that
$f\left(  z\right)  <f\left(  x_{0}-T\right)  .$ Let $S_{1}=S\cap\left[
z,x_{0}-T\right]  .$ We necessarily have that $\sup_{S_{1}}f>f\left(
x_{0}-T\right)  $ as well. We define $x_{-1}$ to be any point in $S_{1}$ where
$f\left(  x_{-1}\right)  =\sup_{S_{1}}f.$ Clearly, $x_{-1}$ is a local maxima
of $f,$ while $x_{0}-x_{-1}<2T.$

We proceed to define the sequences $x_{i}$ by induction, $i=0,-1,-2,...~.$ It
is not excluded that $x_{i}\in\left[  x_{i+1}-T,x_{i+1}\right]  $ for some
$i.$ However that means in particular that the point $x_{i}$ is the first
maximum point of the function $f$ on the segment $\left[  x_{i+1}%
-T,x_{i+1}\right]  ,$ and since $f\left(  x_{i-1}\right)  >f\left(
x_{i}\right)  ,$ we have that $x_{i-1}\notin\left[  x_{i+1}-T,x_{i+1}\right]
$ for all $i.$ We stop when we arrive to a first value below $2T$.

$iii)$ In the same way, starting from a local minima $y_{0},$ we can construct
a sequence $y_{i}$ of local minima, such that

\begin{itemize}
\item $y_{0}>y_{-1}>y_{-2}>...,$

\item $y_{i}-y_{i-1}<2T,\;\;y_{i}-y_{i-2}\geq T$ for all $i,$

\item $0<y_{-n}<2T,$

\item $f\left(  y_{i-1}\right)  \leq f\left(  x\right)  $ for any $y_{i-1}\leq
x, $ and $f\left(  y_{i-1}\right)  <f\left(  y_{i}\right)  ,$

\item for every $x\in\left[  y_{i-1},y_{i}-T\right]  $ we have $f\left(
x\right)  \leq f\left(  y_{i}\right)  $ (if the segment is non-empty).
\end{itemize}

We can suppose additionally that $x_{0}\geq y_{0}\geq x_{-1}.$

$iv)$ Note that the (finite) sequence $x_{i}$ do depend on the initial local
minima $x_{0},$ which was used for the starter. The bigger $x_{0}$ is, the
longer the sequence $x_{i}$ is. So let us introduce the sequence
$x_{0}^{\left(  N\right)  }$ of such starters, and we suppose that
$x_{0}^{\left(  N\right)  }\geq N.$ In that way we will obtain the family
$x_{i}^{\left(  N\right)  }$ of sequences of local maximums of $f,$
$i=0,-1,...,-n\left(  f,x_{0}^{\left(  N\right)  }\right)  ,$ with $n\left(
f,x_{0}^{\left(  N\right)  }\right)  \rightarrow\infty$ as $N\rightarrow
\infty.$ (Of course, in well may happen that for different $N$-s the
corresponding sequences share many common terms.)

Denote by $M$ the limit $\liminf_{N\rightarrow\infty}f\left(  x_{0}^{\left(
N\right)  }\right)  .$ In the same way we can introduce the limit
$m=\limsup_{N\rightarrow\infty}f\left(  y_{0}^{\left(  N\right)  }\right)  .$
Clearly, $M\geq m,$ and if we can show that $M=m,$ then we are done. So we
suppose that $M-m>0,$ and we will bring that to contradiction.

$v)$ Let us fix $\varepsilon>0,$ $\varepsilon<\frac{M-m}{10},$ which is
possible if $M-m>0.$ Then one can choose $N$ so large, that at least 99\% of
terms of the sequence $f\left(  x_{i-1}^{\left(  N\right)  }\right)  -f\left(
x_{i}^{\left(  N\right)  }\right)  $ are less than $\frac{\varepsilon^{2}}%
{2}.$ We will fix that value of $N,$ and we will omit $N$ from our notation.
Therefore without loss of generality we can assume that for some $i$ (in fact,
for many) $\,$we have $f\left(  x\right)  <f\left(  x_{i}\right)
+\varepsilon^{2}$ for all $x\in\left[  x_{i}-T,x_{i}\right]  .$ Therefore for
the set $K\equiv\left\{  x\in\left[  x_{i}-T,x_{i}\right]  :f\left(  x\right)
>f\left(  x_{i}\right)  -\varepsilon\right\}  $ we have:
\begin{equation}
\int_{K-\left(  x_{i}-T\right)  }q_{x_{i}}\left(  t\right)  \,dt>1-\varepsilon
.\label{46}%
\end{equation}
Hence, for its Lebesgue measure we have
\[
\mathrm{mes}\left\{  K\right\}  \geq\frac{1-\varepsilon}{C},
\]
with $C=\sup f.$

Consider now the corresponding sequence of minima, $\left\{  y_{k}\right\}  ,$
and the segments $\left[  y_{k}-T,y_{k}\right]  .$ We claim that the set $K$
has to belong to the union of these segments. That would be evident if the
segments in question were covering the corresponding region without any holes.
However, that is not necessarily the case, and there can be holes between the
segments, since in general the differences $y_{i}-y_{i-1}$ can be bigger than
$T.$ Yet, this does not present a problem, since by construction the function
$f$ is smaller than $m$ outside the union of the segments $\left[
y_{k}-T,y_{k}\right]  ,$ which implies that the set $K$ indeed is covered by
these segments. Since $\mathrm{diam\,}\left(  K\right)  \leq T,$ there exists
$k=k\left(  K\right)  ,$ such that $K\subset\left[  y_{k-1}-T,y_{k-1}\right]
\cup\left[  y_{k}-T,y_{k}\right]  \cup\left[  y_{k+1}-T,y_{k+1}\right]  .$
Without loss of generality we can assume the set $K$ ``fits into $\left[
y_{k}-T,y_{k}\right]  $'', in the sense that
\[
\mathrm{mes}\left\{  K\cap\left[  y_{k}-T,y_{k}\right]  \right\}  \geq
\frac{\mathrm{mes}\left\{  K\right\}  }{3}\geq\frac{1-\varepsilon}{3C},
\]
while we have $f\left(  x\right)  >f\left(  y_{k}\right)  -\varepsilon^{2}$
for all $x\in\left[  y_{k}-T,y_{k}\right]  .$ So we have
\begin{equation}
\int_{\left\{  K\cap\left[  y_{k}-T,y_{k}\right]  \right\}  -\left(
y_{k}-T\right)  }q_{y_{k}}\left(  t\right)  \,dt\geq\bar{\kappa}\left(
\frac{1-\varepsilon}{3C}\right)  ,\label{47}%
\end{equation}
where we define the function $\bar{\kappa}\left(  \alpha\right)  $ by
\[
\bar{\kappa}\left(  \alpha\right)  =\inf_{A\subset\left[  0,T\right]
:\,\mathrm{mes}\left\{  A\right\}  \geq\alpha}\int_{A}\kappa\left(  t\right)
\,dt.
\]
By construction, the set $K\cap\left[  y_{k}-T,y_{k}\right]  $ is disjoint
from the set $L\subset\left[  y_{k}-T,y_{k}\right]  ,$ which is defined by
$L=\left\{  x\in\left[  y_{k}-T,y_{k}\right]  :f\left(  x\right)  <f\left(
y_{k}\right)  +\varepsilon\right\}  .$ Since
\[
f\left(  y_{k}\right)  =\int_{0}^{T}f\left(  y_{k}-t\right)  q_{y_{k}}\left(
t\right)  \,dt,
\]
we have similar to (\ref{46}) that
\begin{equation}
\int_{L-\left(  y_{k}-T\right)  }q_{y_{k}}\left(  t\right)  \,dt>1-\varepsilon
.\label{48}%
\end{equation}
But since $q_{y_{k}}\left(  t\right)  \,dt$ is a probability measure, we
should have that
\[
\bar{\kappa}\left(  \frac{1-\varepsilon}{3C}\right)  +1-\varepsilon\leq1,
\]
because of (\ref{47}), (\ref{48}). This, however, fails once $\varepsilon$ is
small enough.
\end{proof}

\section{Self-averaging $\Longrightarrow$ relaxation: infinite range case
\label{infrange}}

We return to the equation (\ref{36}), $f\left(  x\right)  =\left[  f\ast
q_{x}\right]  \left(  x\right)  .$ Now we will not suppose that the measures
$q_{x}$ have finite support. Instead we suppose that

\begin{enumerate}
\item The family $q_{x}$ has the following compactness property: for every
$\varepsilon>0$ there exists a value $K\left(  \varepsilon\right)  ,$ such
that
\begin{equation}
\int_{0}^{K\left(  \varepsilon\right)  }q_{x}\left(  t\right)  \,dt\geq
1-\varepsilon\label{111}%
\end{equation}
uniformly in $x.$

\item For every $T$ the (monotone continuous) function
\begin{equation}
F_{T}\left(  \delta\right)  =\inf_{x\geq X\left(  T\right)  }\inf
_{\substack{D\subset\left[  0,T\right]  : \\\mathrm{mes}D\geq\delta}}\int
_{D}q_{x}\left(  t\right)  \,dt\label{112}%
\end{equation}
is positive once $\delta>0,$ for some choice of the function $X\left(
T\right)  <\infty.$

\item The family $q_{x}$ is such that the function $f,$ with solves $\left(
\ref{36}\right)  ,$ is Lipschitz, with Lipschitz constant $\mathcal{L}%
=\mathcal{L}\left(  \left\{  q_{\cdot}\right\}  \right)  .$
\end{enumerate}

As we know from the Section \ref{ten}, these conditions are indeed satisfied
in the specific case of the non-linear Markov process and the equation
$\left(  \ref{35}\right)  .$ Indeed, $\left(  \ref{111}\right)  $ follows from
Lemma \ref{la} and Lemma \ref{Qu}, $\left(  \ref{112}\right)  $ -- from Lemma
\ref{lowerb} and Lemma \ref{pusto}, while Lipschitz property follows from
Lemma \ref{Lip}.

\begin{theorem}
\label{infinite range} Let $f$ satisfies $f\left(  x\right)  =\left[  f\ast
q_{x}\right]  \left(  x\right)  $ for $x\geq0,$ with the kernels $q_{x}$
having three properties listed above. Then $f$ relaxes to a constant value as
$x\rightarrow\infty.$
\end{theorem}

\subsection{Approaching stationary point \label{appp}}

\begin{lemma}
\label{T} $i)$ Let $M=\limsup_{x\rightarrow\infty}f\left(  x\right)  .$ Then
for every $T$ and every $\varepsilon$ given there exists some value $K_{1},$
such that
\[
\inf_{x\in\left[  K_{1},K_{1}+T\right]  }f\left(  x\right)  \geq
M-\varepsilon.
\]
$ii)$ Let $m=\liminf_{x\rightarrow\infty}f\left(  x\right)  .$ Then for every
$T$ and every $\varepsilon$ given there exists some value $K_{2},$ such that
\[
\sup_{x\in\left[  K_{2},K_{2}+T\right]  }f\left(  x\right)  \leq
m+\varepsilon.
\]
Moreover, the conclusions of the lemma remains valid if the function $f$
satisfies a weaker equation (see $\left(  \ref{103}\right)  $)
\begin{equation}
f\left(  x\right)  =\left(  1-\varepsilon\left(  x\right)  \right)  \left[
f\ast q_{x}\right]  \left(  x\right)  +\varepsilon\left(  x\right)  Q\left(
x\right)  ,\label{211}%
\end{equation}
with $\varepsilon\left(  x\right)  \rightarrow0$ as $x\rightarrow\infty$ and
$Q\left(  \cdot\right)  \leq C.$
\end{lemma}

\begin{proof}
$i)$ Let $\delta>0.$ Then there exists a value $S=S\left(  \delta\right)  >0,
$ such that for all $x>S$ we have $f\left(  x\right)  <M+\delta,$ and
$\varepsilon\left(  x\right)  Q\left(  x\right)  <\frac{\delta}{2}.$ Further,
there exists a value $R=R\left(  \delta\right)  >S,$ such that for all $y\geq
R$
\[
\int_{R-S}^{\infty}q_{y}\left(  t\right)  dt<\delta,
\]
see (\ref{111}). Finally, there exists a point $y=y\left(  \delta\right)
>R+T,$ such that $f\left(  y\right)  >M-\frac{\delta}{2}.$ Due to the equation
$\left(  \ref{211}\right)  $ we have
\[
f\left(  y\right)  =\left(  1-\varepsilon\left(  y\right)  \right)  \left[
\int_{0}^{y-S}f\left(  y-t\right)  q_{y}\left(  t\right)  \,dt+\int
_{y-S}^{\infty}f\left(  y-t\right)  q_{y}\left(  t\right)  \,dt\right]
+\varepsilon\left(  y\right)  Q\left(  y\right)  .
\]
Let $\Delta>0,$ and $A=\left\{  x\in\left[  y-T,y\right]  :f\left(  x\right)
<M-\Delta\right\}  ,$ while $a=\int_{A}q_{y}\left(  t\right)  \,dt.$ We want
to show that the measure $a$ has to be small for small $\delta$. Splitting the
first integral into two, according to whether the point $y-t$ is in $A$ or
not, we have
\[
M-\delta<a\left(  M-\Delta\right)  +\left(  1-a-\delta\right)  \left(
M+\delta\right)  +\delta C,
\]
so
\[
a<\delta\frac{C+2-M}{\Delta},
\]
which goes to zero with $\delta,$ provided $\Delta$ is fixed. Therefore
\[
\mathrm{mes}\left\{  A\right\}  \leq F_{T}^{-1}\left(  \delta\frac
{C+2-M}{\Delta}\right)  ,
\]
(see (\ref{112})). Since $F_{T}^{-1}\left(  u\right)  \rightarrow0$ as
$u\rightarrow0,$ that proves that for any given $\Delta$ the Lebesgue measure
$\mathrm{mes}\left\{  A\right\}  \rightarrow0$ as $\delta\rightarrow0$. Since
the function $f$ is Lipschitz, we conclude that $\inf_{x\in\left[
y-T,y\right]  }f\left(  x\right)  \geq M-\Delta-\mathcal{L}\mathrm{mes}%
\left\{  A\right\}  \geq M-2\Delta,$ provided $\delta$ is small enough. Taking
$\Delta=\varepsilon/2$ finishes the proof.

$ii)$ Let $\delta>0.$ Then there exists a value $S>0,$ such that for all $x>S$
we have $f\left(  x\right)  >m-\delta.$ Again, take $R>S,$ such that for all
$y\geq R$
\[
\int_{0}^{R-S}q_{y}\left(  t\right)  dt>1-\delta.
\]
Finally, there exists a point $y>R+T,$ such that $f\left(  y\right)
<m+\delta.$ Due to the equation $\left(  \ref{211}\right)  $ we have
\begin{equation}
f\left(  y\right)  >\left(  1-\varkappa\right)  \left[  \int_{0}^{y-S}f\left(
y-t\right)  q_{y}\left(  t\right)  \,dt+\int_{y-S}^{\infty}f\left(
y-t\right)  q_{y}\left(  t\right)  \,dt\right]  ,\label{113}%
\end{equation}
where $\varkappa$ can be supposed arbitrarily small. Let $\Delta>0,$ and

\noindent$A=\left\{  t\in\left[  0,T\right]  :f\left(  y-t\right)
>m+\Delta\right\}  ,$ while $a=\int_{A}q_{y}\left(  t\right)  \,dt.$ We want
to show that the measure $a$ has to be small for small $\delta$. Splitting the
first integral into two, according to whether the point $y-t$ is in $A$ or
not, and disregarding the second one, we have
\[
m+\delta>\left(  1-\varkappa\right)  \left[  a\left(  m+\Delta\right)
+\left(  1-\delta-a\right)  \left(  m-\delta\right)  \right]  .
\]
For $\varkappa$ so small that $\varkappa\left[  a\left(  m+\Delta\right)
+\left(  1-a-\delta\right)  \left(  m-\delta\right)  \right]  <\delta,$ we
have
\[
m+2\delta>a\left(  m+\Delta\right)  +\left(  1-a-\delta\right)  \left(
m-\delta\right)  ,
\]
so
\[
a<\delta\frac{m+3}{\Delta},
\]
which goes to zero with $\delta,$ provided $\Delta$ is fixed. Therefore
\[
\mathrm{mes}\left\{  A\right\}  \leq F_{T}^{-1}\left(  \delta\frac{m+3}%
{\Delta}\right)  ,
\]
and the rest of the argument coincides with that of the part $i).$
\end{proof}

\subsection{Absorbing by stationary point \label{absorbing}}

We now will show that if $f$ satisfies $\left(  \ref{36}\right)  ,$ then the
property $\inf_{x\in\left[  K,K+T\right]  }f\left(  x\right)  \geq
M-\varepsilon$ implies that for all $x>K+T$
\begin{equation}
f\left(  x\right)  >M-\varepsilon-c\left(  T\right)  ,\label{203}%
\end{equation}
with $c\left(  T\right)  \rightarrow0$ as $T\rightarrow\infty.$ That clearly
implies relaxation. (Note that we do not claim that $\left(  \ref{203}\right)
$ holds for the solutions of $\left(  \ref{211}\right)  $). We will show it
under the extra assumption that the distribution $p\left(  \cdot\right)  $ has
finite moment of some order above $4.$ This assumption will be used only
throughout the rest of the present subsection.

Using the linearity of $\left(  \ref{36}\right)  ,$ we will rewrite our
problem slightly, in order to simplify the notation.

Let the function $f\geq0$ satisfies $f\left(  x\right)  =\left[  f\ast
q_{x}\right]  \left(  x\right)  $ for $x>0,$ and

$i)$ $f\left(  x\right)  >1$ for $x\in\left[  -T,0\right]  ,$

$ii)$ for some $\beta>1$ and $B<\infty$ and for every $x$ we have
\begin{equation}
\int_{0}^{\infty}t^{\beta}q_{x}\left(  t\right)  \,dx\leq B,\label{110}%
\end{equation}
compare with $\left(  \ref{303}\right)  .$ We want to derive from that data
that for some $c\left(  T\right)  >0,$ $c\left(  T\right)  \rightarrow0$ as
$T\rightarrow\infty$
\[
f\left(  x\right)  >1-c\left(  T\right)  \text{ for all }x>0.
\]

Denote by
\[
g_{0}\left(  x\right)  =\left\{
\begin{array}
[c]{cc}%
1 & x\in\left[  -T,0\right] \\
0 & x\notin\left[  -T,0\right]
\end{array}
\right.  .
\]
Since $f\geq g,$ we have $f\left(  x\right)  \geq g_{1}\left(  x\right)
=\left[  g_{0}\ast q_{x}\right]  \left(  x\right)  $ for $x\geq0.$ We define
$g_{1}\left(  x\right)  =g_{0}\left(  x\right)  $ for $x<0.$ Iterating, we
have $f\left(  x\right)  \geq g_{n}\left(  x\right)  ,$ where
\[
g_{n}\left(  x\right)  =\left\{
\begin{array}
[c]{cc}%
g_{0}\left(  x\right)  & x<0\\
\left[  g_{n-1}\ast q_{x}\right]  \left(  x\right)  & x\geq0
\end{array}
\right.  .
\]

Hence, $f\left(  x\right)  \geq g_{\infty}\left(  x\right)  .$ The function
$g_{\infty}\left(  x\right)  $ has the following probabilistic interpretation:
we have a Markov chain on $\mathbb{R}^{1},$ where transition from the point
$x$ is governed by transition densities $q_{x}$ to make the step (to the
left), (and which steps to the left are defined in an arbitrary way for
$x\leq0$); then the value $g_{\infty}\left(  x\right)  $ for $x>0$ is the
probability that starting from $x$ we will visit the interval $\left[
-T,0\right]  .$ The question now is about the lower bound on $g_{\infty
}\left(  x\right)  $ over all possible $q_{x}$ from our class.

So let us take $x>0,$ and let start the Markov chain $X_{n}$ from $x,$ (i.e.
$X_{0}=x$), which goes to the left, and which makes a transition from $y$ to
$y-t$ with the probability $q_{y}\left(  t\right)  dt.$ We need to know the
probability of the event
\[
\mathbb{P}_{x}\left\{  \text{there exists }n\text{ such that }X_{n}\in\left[
-T,0\right]  \right\}  .
\]
In other words, we want to know the probability of $X_{\left\{  \cdot\right\}
}$ visiting $\left[  -T,0\right]  .$ We would like to show that
\begin{equation}
\mathbb{P}_{x}\left\{  X_{\left\{  \cdot\right\}  \text{ }}\text{ visits
}\left[  -T,0\right]  \right\}  \geq\gamma\left(  \beta,B,T\right) \label{106}%
\end{equation}
with
\[
\gamma\left(  \beta,B,T\right)  \rightarrow1\text{ as }T\rightarrow\infty
\]
uniformly over the families $q_{x}$ from our class.

Note, however, that in general such an estimate does not hold. For example,
the process $X_{\left\{  \cdot\right\}  \text{ }}$ can well stay positive for
all times. The more interesting example where the process goes to $-\infty$,
follows, so we will need further restrictions on the family $q_{x}$.

\textbf{Example. }Let $T$ be given. We will construct the family $q_{x}^{T}$
from our class (\ref{110}), such that for the corresponding process
$X_{\left\{  \cdot\right\}  }^{T}$
\[
\mathbb{P}_{x}\left\{  X_{\left\{  \cdot\right\}  }^{T}\,\text{visits }\left[
0,T\right]  \right\}  =0.
\]
We define $q_{x}^{T}\left(  t\right)  $ for $x\in(k,k+1]$ with integer
$k\neq0$ to be any distribution localized in the segment $\left[
k-1,k\right]  $ (the uniform distribution on $\left[  k-1,k\right]  $ is OK).
For $x\in(\frac{1}{2^{k}},\frac{1}{2^{k-1}}],$ $k=1,2,...,$ it is defined by
\[
q_{x}^{T}\left(  t\right)  =\left\{
\begin{array}
[c]{ll}%
e^{-t} & \text{ if }t>T+1\\
2^{k+1}\left(  1-\int_{T+1}^{\infty}e^{-t}dt\right)  & \text{ if }t\in\left[
x-\frac{1}{2^{k}},x-\frac{1}{2^{k+1}}\right]  \text{ }\\
0 & \text{ otherwice.}%
\end{array}
\right.
\]
For $x\leq0$ it is defined in an arbitrary way. $\blacksquare$

The mechanism of violating the relation (\ref{106}) is that the time the
process $X_{\left\{  \cdot\right\}  }^{T}$ can spend in the segment $\left[
0,1\right]  $ is unbounded in $T$. As the following theorem shows, this
feature is the only obstruction for the statement desired to hold.

\begin{theorem}
Consider the Markov chain $X_{\left\{  \cdot\right\}  }$ defined above via the
transition densities $q_{x}\left(  t\right)  .$ Suppose that condition
(\ref{110}) holds, and that in addition these densities are uniformly bounded
in the vicinity of the origin: for all real $x$ and all $t$ in the segment
$\left[  0,1\right]  $, say,
\begin{equation}
q_{x}\left(  t\right)  \leq C.\label{107}%
\end{equation}
Then for some $\gamma\left(  \beta,B,C,T\right)  \rightarrow1$ as
$T\rightarrow\infty$ we have:
\[
\mathbb{P}_{x}\left\{  X_{\left\{  \cdot\right\}  }\text{visits }\left[
-T,0\right]  \right\}  \geq\gamma\left(  \beta,B,C,T\right)  .
\]

\end{theorem}

The condition $\left(  \ref{107}\right)  $ holds in the case of NMP, see
estimate $\left(  \ref{157}\right)  \ $from Lemma \ref{oc}.

\begin{proof}
We will estimate the probability of the complementary event:
\begin{align*}
& \mathbb{P}_{x}\left\{  X_{\left\{  \cdot\right\}  }\text{misses }\left[
-T,0\right]  \right\} \\
& =\sum_{k=0}^{\infty}\int_{0}^{x}\left[  \int_{y+T}^{\infty}q_{y}\left(
t\right)  \,dt\right]  P_{k}\left(  x,dy\right)  .
\end{align*}
Here $P_{k}\left(  x,dy\right)  $ is the probability distribution of the chain
$X_{\left\{  \cdot\right\}  \text{ }}$ after $k$ steps, and the expression
$\left[  \int_{y+T}^{\infty}q_{y}\left(  t\right)  \,dt\right]  P_{k}\left(
x,dy\right)  $ is the probability that the chain $X_{\left\{  \cdot\right\}
\text{ }}$ arrives after $k$ steps to the location $y,$ and then makes a jump
over the segment $\left[  -T,0\right]  .$ (In our case the probability of the
event that $X_{\cdot}$ never becomes negative equals zero.) So we have
\begin{align*}
& \mathbb{P}_{x}\left\{  X_{\left\{  \cdot\right\}  \text{ }}\text{misses
}\left[  -T,0\right]  \right\} \\
& \leq\int_{0}^{x}B\left(  y+T\right)  ^{-\beta}\sum_{k=0}^{\infty}%
P_{k}\left(  x,dy\right) \\
& \leq\sum_{n=0}^{\left[  x\right]  +1}B\left(  n+T\right)  ^{-\beta}%
\sum_{k=0}^{\infty}P_{k}\left(  x,\left[  n,n+1\right]  \right)  ,
\end{align*}
where $P_{k}\left(  x,\left[  n,n+1\right]  \right)  $ is the probability of
the event $X_{k}\in\left[  n,n+1\right]  ,$ and where in the second line we
are using the following simple estimate:
\[
\int_{r}^{\infty}q_{y}\left(  t\right)  \,dt=r^{-\beta}\int_{r}^{\infty
}r^{\beta}q_{y}\left(  t\right)  \,dt\leq r^{-\beta}\int_{0}^{\infty}t^{\beta
}q_{y}\left(  t\right)  \,dt.
\]
Now,
\begin{align}
& \sum_{k=0}^{\infty}\mathbb{P}_{x}\left\{  X_{k}\in\left[  n,n+1\right]
\right\} \label{108}\\
& =\sum_{k=0}^{\infty}\sum_{l<k}\mathbb{P}_{x}\left\{  X_{k}\in\left[
n,n+1\right]  ,X_{l}>n+1,X_{l+1}\in\left[  n,n+1\right]  \right\} \nonumber\\
& =\sum_{l=0}^{\infty}\mathbb{P}_{x}\left\{  X_{l}>n+1,X_{l+1}\in\left[
n,n+1\right]  \right\} \nonumber\\
& \times\sum_{k>0}\mathbb{P}_{x}\left\{  X_{l+k}\in\left[  n,n+1\right]
\Bigm|X_{l}>n+1,X_{l+1}\in\left[  n,n+1\right]  \right\}  .\nonumber
\end{align}
Let now the random variables $\zeta_{i}$ be i.i.d., uniformly distributed in
the segment $\left[  0,\frac{1}{C}\right]  ,$ where $C$ is the same as in
(\ref{106}). Then is easy to see that
\[
\mathbb{P}_{x}\left\{  X_{l+k}\in\left[  n,n+1\right]  \Bigm|X_{l}%
>n+1,X_{l+1}\in\left[  n,n+1\right]  \right\}  \leq\mathbf{\Pr}\left\{
\zeta_{1}+...+\zeta_{k}\leq1\right\}  .
\]
Since the last probability decays exponentially in $k,$ while

\noindent$\sum_{l=0}^{\infty}\mathbb{P}_{x}\left\{  X_{l}>n+1,X_{l+1}%
\in\left[  n,n+1\right]  \right\}  =1,$ we conclude that
\[
\sum_{k=0}^{\infty}\mathbb{P}_{x}\left\{  X_{k}\in\left[  n,n+1\right]
\right\}  \leq K\left(  C\right)  .
\]
Since the series $\sum n^{-\beta}$ converges for $\beta>1,$ the proof follows.
\end{proof}

\section{Self-averaging $\Longrightarrow$ relaxation: noisy case \label{nc}}

In this section we establish the relaxation for the NMP with general initial
condition. Using the fact that the Poisson rate $\lambda\left(  x\right)
=\lambda_{\mu}\left(  x\right)  $ of the NMP with initial state $\mu$
satisfies the equation
\[
\lambda\left(  x\right)  =\left(  1-\varepsilon_{\lambda,\mu}\left(  x\right)
\right)  \left[  \lambda\ast q_{\lambda,\mu,x}\right]  \left(  x\right)
+\varepsilon_{\lambda,\mu}\left(  x\right)  Q_{\lambda,\mu}\left(  x\right)  ,
\]
we will prove the following

\begin{theorem}
\label{noisy case} Let the initial state $\mu$ of the NMP $\mu_{t}$ is such
that both the expected service time $S\left(  \mu\right)  $ and the mean queue
length $N\left(  \mu\right)  $ are finite. Then the limit $c=\lim
_{x\rightarrow\infty}\lambda\left(  x\right)  $ exists; moreover, $\mu
_{t}\rightarrow\nu_{c}$ as $t\rightarrow\infty,$ where $\nu_{c}$ is the
invariant measure of NMP, such that $\lambda_{\nu_{c}}\left(  x\right)  \equiv
c.$ Also $N\left(  \nu_{c}\right)  =N\left(  \mu\right)  .$
\end{theorem}

We are not able to prove this theorem in the generality of the previous
Sections. Below we will use all the specific features of the NMP-s, and in
particular we will use the comparison between different NMP-s and GFP-s,
corresponding to various initial states and input rates. The comparison
mentioned is based on the coupling arguments.

\subsection{Coupling}

\begin{definition}
Let $\mu_{1},\mu_{2}$ be two states on $\Omega.$ We call the state $\mu_{1}$
to be \textbf{higher} than $\mu_{2},$ $\mu_{1}\succcurlyeq\mu_{2},$ if there
exists a coupling $P\left[  d\omega_{1},d\omega_{2}\right]  $ between the
states $\mu_{1},\mu_{2},$ with the property:
\[
P\left[  \left(  \Omega\times\Omega\right)  ^{>}\right]  =1,
\]
where
\[
\left(  \Omega\times\Omega\right)  ^{>}=\left\{  \left[  \left(  n_{1}%
,\tau_{1}\right)  ,\left(  n_{2},\tau_{2}\right)  \right]  \in\Omega
\times\Omega:n_{1}\geq n_{2}\right\}  .
\]

\end{definition}

Clearly, if $\mu_{1}\succcurlyeq\mu_{2},$ then $N\left(  \mu_{1}\right)  \geq
N\left(  \mu_{2}\right)  .$

Next, we introduce the stronger relation.

\begin{definition}
Let $\mu_{1},\mu_{2}$ be two states on $\Omega.$ We call the state $\mu_{1}$
to be \textbf{taller} than $\mu_{2},$ $\mu_{1}\curlyeqsucc\mu_{2},$ if there
exists a coupling $P\left[  d\omega_{1},d\omega_{2}\right]  $ between the
states $\mu_{1},\mu_{2},$ with the property:
\[
P\left[  \left(  \Omega\times\Omega\right)  ^{\gg}\right]  =1,
\]
where
\[
\left(  \Omega\times\Omega\right)  ^{\gg}=\left\{  \left[  \left(  n_{1}%
,\tau_{1}\right)  ,\left(  n_{2},\tau_{2}\right)  \right]  \in\Omega
\times\Omega:\tau_{1}=\tau_{2},n_{1}\geq n_{2}\text{ or }\left(  n_{2}%
,\tau_{2}\right)  =\mathbf{0}\right\}  .\text{ }%
\]

\end{definition}

\begin{lemma}
Let $\mu_{1}\left(  0\right)  ,\,\mu_{2}\left(  0\right)  $ be two initial
states on $\Omega$ at $t=0,$ and $\lambda_{1}\left(  t\right)  ,$
$\,\lambda_{2}\left(  t\right)  ,$ $t\geq0$ be two Poisson densities of the
input flows. The service time distribution is the same $\eta$ as before. Let
$\mu_{i}\left(  t\right)  $ be two corresponding GFP-s. Suppose that $\mu
_{1}\left(  0\right)  \curlyeqsucc\mu_{2}\left(  0\right)  ,$ and that
$\lambda_{1}\left(  t\right)  \geq\lambda_{2}\left(  t\right)  .$ Then
$\mu_{1}\left(  t\right)  \succcurlyeq\mu_{2}\left(  t\right)  ,$ so in
particular
\[
N\left(  \mu_{1}\left(  t\right)  \right)  \geq N\left(  \mu_{2}\left(
t\right)  \right)  .
\]
Also, there exists a coupling between the processes such that for almost every
trajectory $\left(  \omega_{1}\left(  t\right)  ,\omega_{2}\left(  t\right)
\right)  $%
\begin{equation}
S\left(  \omega_{1}\left(  t\right)  \right)  \geq S\left(  \omega_{2}\left(
t\right)  \right)  .\label{304}%
\end{equation}

\end{lemma}

\begin{proof}
To see this let us construct the coupling between the processes $\mu
_{i}\left(  t\right)  .$ Let us color the customers, arriving according to the
$\lambda_{2}\left(  t\right)  $ flow, as red. We also assign the red color to
the customers which were present at time $t=0$ from the initial state $\mu
_{2}\left(  0\right)  $. Let $\gamma\left(  t\right)  =\lambda_{1}\left(
t\right)  -\lambda_{2}\left(  t\right)  ,$ and consider $\gamma\left(
t\right)  $ as the extra input flow of blue customers (with independent
service times). We also add blue customers at time $t=0,$ which are needed to
complete the state $\mu_{2}\left(  0\right)  $ up to $\mu_{1}\left(  0\right)
.$ Then the total (color blind) flow coincides with $\lambda_{1}$ flow, while
the total (color blind) process coincides with $\mu_{1}\left(  t\right)  .$

The service rule for the two-colored process is color blind: all the customers
are served in order of their arrival time. We claim now that along every
coupled trajectory $\left(  \omega_{1}\left(  t\right)  ,\omega_{2}\left(
t\right)  \right)  $ we have $r\left(  \omega_{1}\left(  t\right)  \right)
\geq r\left(  \omega_{2}\left(  t\right)  \right)  ,$ where $r\left(
\cdot\right)  $ is the number of red customers at the moment $t,$ waiting to
be served. That evidently will prove our statement.

Clearly, the number $r\left(  \omega\left(  t\right)  \right)  $ is the
difference,
\[
r\left(  \omega\left(  t\right)  \right)  =\mathcal{A}\left(  \omega\left(
t\right)  \right)  -\mathcal{B}\left(  \omega\left(  t\right)  \right)  ,
\]
where $\mathcal{A}\left(  \omega\left(  t\right)  \right)  $ is the total
number of red customers, having arrived before $t,$ while $\mathcal{B}\left(
\omega\left(  t\right)  \right)  $ is the total number of red customers, who
left the system before $t.$ Clearly, $\mathcal{A}\left(  \omega_{1}\left(
t\right)  \right)  =\mathcal{A}\left(  \omega_{2}\left(  t\right)  \right)  .$
Let us show that $\mathcal{B}\left(  \omega_{1}\left(  t\right)  \right)
\leq\mathcal{B}\left(  \omega_{2}\left(  t\right)  \right)  .$

This is easy to see once one visualizes the procedure of resolving the rod
conflicts, which corresponds to our service rule, for the two-colored rod
case. Namely, one has first to put all the red rods, and resolve all their
conflicts by shifting some of them to the right accordingly. The number of
thus obtained rods to the left of the point $t$ is the number $\mathcal{B}%
\left(  \omega_{2}\left(  t\right)  \right)  .$ Clearly, if one adds some blue
rods to the red ones, then each red rod would be shifted to the right by at
least the same amount as without the blue rods. As a result, every red rod
would either stay where it was, or move to the right, so indeed $\mathcal{B}%
\left(  \omega_{1}\left(  t\right)  \right)  \leq\mathcal{B}\left(  \omega
_{2}\left(  t\right)  \right)  .$

The relation $\left(  \ref{304}\right)  $ is evident.
\end{proof}

\subsection{Convergence}

Consider a General Flow Process $\mu\left(  t\right)  $ with initial state
$\mu\left(  0\right)  =\bar{\nu}$ at $T=0$ and the input rate $\lambda\left(
t\right)  \equiv c<1$ (i.e. the usual queueing system $M|GI|1$). This is an
ergodic process, so the weak limit
\[
\lim_{t\rightarrow\infty}\mu_{\bar{\nu},c}\left(  t\right)  =\nu_{c}%
\]
exists and does not depend on the initial state $\bar{\nu}.$ We would like to
show that if $N\left(  \bar{\nu}\right)  <\infty,$ then also
\begin{equation}
\lim_{t\rightarrow\infty}N\left(  \mu_{\bar{\nu},c}\left(  t\right)  \right)
=N\left(  \nu_{c}\right) \label{50}%
\end{equation}
(see $\left(  \ref{03}\right)  $). This however is not true in general, and
extra assumptions are needed in order to have such convergence.

\begin{lemma}
\label{Ryba} Suppose additionally that $S\left(  \bar{\nu}\right)  <\infty.$
Then
\[
\lim_{t\rightarrow\infty}N\left(  \mu_{\bar{\nu},c}\left(  t\right)  \right)
=N\left(  \nu_{c}\right)  .
\]
Moreover, for every $s<\infty$ the convergence $N\left(  \mu_{\bar{\nu}%
,c}\left(  t\right)  \right)  \rightarrow N\left(  \nu_{c}\right)  $ is
uniform on the set of all initial states $\bar{\nu}$ satisfying $S\left(
\bar{\nu}\right)  \leq s.$
\end{lemma}

\textbf{Proof of Lemma \ref{Ryba}. }Since $\lim_{t\rightarrow\infty}\mu
_{\bar{\nu},c}\left(  t\right)  =\nu_{c},$ $\lim_{t\rightarrow\infty}N\left(
\mu_{\bar{\nu},c}\left(  t\right)  \right)  \geq N\left(  \nu_{c}\right)  .$
To prove the equality we need to show the uniform integrability for the family
of random variables $N_{\mu_{\bar{\nu},c}\left(  t\right)  }$, which means the
following property: for every $\varkappa>0$ there exists a value
$s_{\varkappa}$ such that for all $t$%
\[
\mathbb{E}_{\mu_{\bar{\nu},c}\left(  t\right)  }\left(  N\left(
\omega\right)  \mathbf{I}_{N\left(  \omega\right)  \geq s_{\varkappa}}\right)
<\varkappa,
\]
where $\mathbf{I}$ stands for the indicator.

Note first that it is enough to show the uniform integrability of the family
$S_{\mu_{\bar{\nu},c}\left(  t\right)  }$ of random variables. Indeed,
consider the event $N\left(  \omega\right)  \geq N.$ Then the conditional
$\mu_{\bar{\nu},c}\left(  t\right)  $-probability of the event $S\left(
\omega\right)  \geq\frac{1}{2}N$ under the condition $N\left(  \omega\right)
\geq N $ goes to $1$ as $N\rightarrow\infty,$ since $\mathbb{E}\left(
\eta\right)  =1.$ Therefore if the family $S_{\mu_{\bar{\nu},c}\left(
t\right)  }$ is not uniformly integrable, so is $N_{\mu_{\bar{\nu},c}\left(
t\right)  }.$

To get uniform integrability we prove the following stochastic domination:

\begin{lemma}
\label{domin}
\begin{equation}
S_{\mu_{\bar{\nu},c}\left(  t\right)  }\preccurlyeq S_{\bar{\nu}}+S_{\nu_{c}%
}.\label{52}%
\end{equation}
(Here in the rhs we mean the sum of two independent random variables.)
\end{lemma}

Since the expectation $S\left(  \bar{\nu}\right)  $ of the first random
variable is finite by assumption, while the expectation of the second equals
$N\left(  \nu_{c}\right)  $ and so is also finite, the uniform integrability
of the family $S_{\mu_{\bar{\nu},c}\left(  t\right)  }$ follows.

\textbf{Proof of Lemma \ref{domin}. }To prove $\left(  \ref{52}\right)  $ we
will use the following construction. Let $x_{1},\eta_{1};x_{2},\eta_{2};... $
be a realization of the flow of customers. It means that at the moment $x_{1}$
the first customer comes, which needs the time $\eta_{1}$ to be served, at the
moment $x_{2}$ the second comes, etc. To every such realization we can assign
the function $W\left(  x\right)  ,$ which is the remaining time duration
needed for the server to serve all the customers who came before the moment
$x.$ That is,
\[
W\left(  x\right)  =\left\{
\begin{array}
[c]{ll}%
0 & \text{for }x<x_{1},\\
\max\left\{  \eta_{1}-\left(  x-x_{1}\right)  ,0\right\}  & \text{for }%
x_{1}\leq x<x_{2},\\
\max\left\{  W\left(  x_{2}-0\right)  +\eta_{2}-\left(  x-x_{2}\right)
,0\right\}  & \text{for }x_{2}\leq x<x_{3},\text{ etc.}%
\end{array}
\right.
\]
If $\mu_{t}$ is a process of states of our server with $\mu_{0}=\delta
_{\mathbf{0}}$, and $\left\{  x_{1},\eta_{1};x_{2},\eta_{2};...\right\}  $ is
its realization, then the random variable $W\left(  t\right)  $ is the same as
$S_{\mu_{t}}.$ With obvious modification the $W$ function is defined for a
process with non-empty initial state $\mu_{0}.$

Consider now two processes: $\mu_{\nu_{c},c}\left(  t\right)  $ and
$\mu_{\delta_{\mathbf{0}},c}\left(  t\right)  .$ The first one is stationary.
Since evidently $\nu_{c}\curlyeqsucc\delta_{\mathbf{0}},$ we have $\mu
_{\nu_{c},c}\left(  t\right)  \succcurlyeq\mu_{\delta_{\mathbf{0}},c}\left(
t\right)  ,$ i.e. $\nu_{c}\succcurlyeq\mu_{\delta_{\mathbf{0}},c}\left(
t\right)  $ for all $t.$ It means also that we can couple the two processes in
such a way that $W^{\mu_{\nu_{c},c}}\left(  t\right)  \geq W^{\mu
_{\delta_{\mathbf{0}},c}}\left(  t\right)  $ with probability one.

Let us see now how the two processes -- $W^{\mu_{\delta_{\mathbf{0}},c}%
}\left(  t\right)  $ and $W^{\mu_{\bar{\nu},c}}\left(  t\right)  $ -- are
related. We consider the natural coupling between $\mu_{\delta_{\mathbf{0}}%
,c}$ and $\mu_{\bar{\nu},c},$ where the latter process is obtained by adding
to a general configuration $\left\{  x_{1},\eta_{1};x_{2},\eta_{2}%
;...\right\}  $ of the former one extra customer $\left\{  x_{0},\Pi
_{0}\right\}  ,$ with $x_{0}=0$ and the independent random variable $\Pi_{0}$
distributed according to $S_{\bar{\nu}}.$ The trajectory $W^{\mu_{\bar{\nu}%
,c}}\left(  t\right)  $ is obtained in the following way: one considers first
the function
\[
\tilde{W}\left(  x\right)  =\Pi_{0}-x+\sum_{i\geq1}\eta_{i}\chi_{x_{i}}\left(
x\right)  ,
\]
where
\[
\chi_{a}\left(  x\right)  =\left\{
\begin{array}
[c]{ll}%
1 & \text{ if }x\geq a,\\
0 & \text{ otherwice.}%
\end{array}
\right.
\]
Let $x_{0}$ be the first moment when $\tilde{W}\left(  x\right)  $ vanishes.
Then
\[
W^{\mu_{\bar{\nu},c}}\left(  x\right)  =\left\{
\begin{array}
[c]{ll}%
\tilde{W}\left(  x\right)  & \text{ if }x\leq x_{0},\\
W^{\mu_{\delta_{\mathbf{0}},c}}\left(  x\right)  & \text{ otherwice.}%
\end{array}
\right.
\]
From this the relation $\left(  \ref{52}\right)  $ follows immediately.

The uniformity of convergence follows from the fact that the function $S$ is a
compact function on $\Omega,$ once the function $R_{\eta}\left(  \tau\right)
$ is unbounded (see $\left(  \ref{181}\right)  $). (This compactness means
that for every $s$ the set $\left\{  \omega\in\Omega:S\left(  \omega\right)
\leq s\right\}  $ is compact.) As a result, the family of initial states
$\bar{\nu}$ satisfying $S\left(  \bar{\nu}\right)  \leq s$ is weakly compact
as well, which together with continuity of the function $N\left(  \mu
_{\bar{\nu},c}\left(  t\right)  \right)  $ in $\bar{\nu}$ and $t$ provides the
claim needed. If the function $R_{\eta}\left(  \tau\right)  $ is uniformly
bounded in $\tau,$ then for some $\xi$ and $C$ the exponential moment
$\mathbb{E}\left(  \exp\left\{  \xi\eta_{\tau}\right\}  \right)  $ $\leq C$
for all $\tau.$ Therefore the family of all possible probability distributions
$F_{\theta}$ of the form
\[
F_{\theta}\left(  x\right)  =\int F_{\eta_{\tau}}\left(  x\right)
~d\theta\left(  \tau\right)  ,
\]
where $\theta$ runs over all probability measures on the semiaxis $\left\{
\tau\geq0\right\}  ,$ is compact. That again implies the uniformity.
$\blacksquare$

$\blacksquare$

\subsection{End of the proof in noisy case}

Let $\mu_{\nu,\lambda_{\nu}\left(  \cdot\right)  }\left(  t\right)  $ be the
non-linear Markov process\textit{\ }with the initial state $\nu,$ having
finite mean queue, $N\left(  \nu\right)  <\infty,$ and finite expected service
time $S\left(  \nu\right)  .$ We will show that the function $\lambda\left(
t\right)  \equiv\lambda_{\nu}\left(  t\right)  $ goes to a limit as
$t\rightarrow\infty.$ The idea is the following:

Suppose $m=\liminf_{t\rightarrow\infty}\lambda\left(  t\right)  <\limsup
_{t\rightarrow\infty}\lambda\left(  t\right)  =M.$ As we already know from
Lemma \ref{T}, for every $T,K$ and every $\varepsilon>0$ there exist values
$K_{1},K_{2}>K$ such that
\begin{equation}
\sup_{x\in\left[  K_{1},K_{1}+T\right]  }\lambda\left(  x\right)  \leq
m+\varepsilon,\label{161}%
\end{equation}
while
\begin{equation}
\inf_{x\in\left[  K_{2},K_{2}+T\right]  }\lambda\left(  x\right)  \geq
M-\varepsilon.\label{162}%
\end{equation}
We want to bring this to contradiction, arguing as follows:

\begin{itemize}
\item First of all, we note that the mean queue, $N\left(  \mu_{\nu
,\lambda_{\nu}\left(  \cdot\right)  }\left(  t\right)  \right)  $ does not
change in time, staying equal to the initial value $N\left(  \nu\right)  .$

On the other hand:

\item We can compare the state $\mu_{\nu,\lambda_{\nu}\left(  \cdot\right)
}\left(  K_{1}+T\right)  $ with the state $\mu_{\left[  \mu_{\nu,\lambda_{\nu
}\left(  \cdot\right)  }\left(  K_{1}\right)  \right]  ,m+\varepsilon}\left(
T\right)  .$ Due to $\left(  \ref{161}\right)  ,$ the latter is higher, so
\begin{equation}
N\left(  \mu_{\left[  \mu_{\nu,\lambda_{\nu}\left(  \cdot\right)  }\left(
K_{1}\right)  \right]  ,m+\varepsilon}\left(  T\right)  \right)  \geq N\left(
\nu\right)  .\label{163}%
\end{equation}
By the same reasoning,
\begin{equation}
N\left(  \mu_{\left[  \mu_{\nu,\lambda_{\nu}\left(  \cdot\right)  }\left(
K_{2}\right)  \right]  ,M-\varepsilon}\left(  T\right)  \right)  \leq N\left(
\nu\right)  .\label{164}%
\end{equation}

\item We then claim that once $T$ is large enough, the state $\mu_{\left[
\mu_{\nu,\lambda_{\nu}\left(  \cdot\right)  }\left(  K_{1}\right)  \right]
,m+\varepsilon}\left(  T\right)  $ is close to the equilibrium $\nu
_{m+\varepsilon},$ so in particular
\begin{equation}
N\left(  \mu_{\left[  \mu_{\nu,\lambda_{\nu}\left(  \cdot\right)  }\left(
K_{1}\right)  \right]  ,m+\varepsilon}\left(  T\right)  \right)  \leq N\left(
\nu_{m+\varepsilon}\right)  +\varepsilon^{\prime},\label{165}%
\end{equation}
with $\varepsilon^{\prime}=\varepsilon^{\prime}\left(  T\right)  \rightarrow0$
as $T\rightarrow\infty.$ By the same reasoning,
\begin{equation}
N\left(  \mu_{\left[  \mu_{\nu,\lambda_{\nu}\left(  \cdot\right)  }\left(
K_{2}\right)  \right]  ,M-\varepsilon}\left(  T\right)  \right)  \geq N\left(
\nu_{M-\varepsilon}\right)  -\varepsilon^{\prime\prime}.\label{166}%
\end{equation}

\item Since $N\left(  \nu_{M-\varepsilon}\right)  >N\left(  \nu_{m+\varepsilon
}\right)  $ once $\varepsilon$ is small, the choice of $\varepsilon^{\prime}$
and $\varepsilon^{\prime\prime}$ such that $\varepsilon^{\prime}%
+\varepsilon^{\prime\prime}<N\left(  \nu_{M-\varepsilon}\right)  -N\left(
\nu_{m+\varepsilon}\right)  $ makes it possible to deduce from the relations
$\left(  \ref{163}\right)  $-$\left(  \ref{166}\right)  $ that
\[
N\left(  \nu\right)  \geq N\left(  \nu_{M-\varepsilon}\right)  -\varepsilon
^{\prime\prime}>N\left(  \nu_{m+\varepsilon}\right)  +\varepsilon^{\prime}\geq
N\left(  \nu\right)  ,
\]
which is inconsistent with the properties of the relation $>$ between the real numbers.
\end{itemize}

We need to prove the relations $\left(  \ref{165}\right)  $ and $\left(
\ref{166}\right)  .$ It turns out that the relation $\left(  \ref{166}\right)
$ is easier. Indeed, to show it, we can compare the state $\mu_{\left[
\mu_{\nu,\lambda_{\nu}\left(  \cdot\right)  }\left(  K_{2}\right)  \right]
,M-\varepsilon}\left(  T\right)  $ with the state $\mu_{\mathbf{0}%
,M-\varepsilon}\left(  T\right)  .$ The latter is evidently lower :
\[
N\left(  \mu_{\left[  \mu_{\nu,\lambda_{\nu}\left(  \cdot\right)  }\left(
K_{2}\right)  \right]  ,M-\varepsilon}\left(  T\right)  \right)  \geq N\left(
\mu_{\mathbf{0},M-\varepsilon}\left(  T\right)  \right)  .
\]
Since $\mu_{\mathbf{0},M-\varepsilon}\left(  T\right)  $ is also lower than
$\nu_{M-\varepsilon},$
\begin{equation}
N\left(  \mu_{\mathbf{0},M-\varepsilon}\left(  T\right)  \right)  \leq
N\left(  \nu_{M-\varepsilon}\right)  .\label{167}%
\end{equation}
Since $\mu_{\mathbf{0},M-\varepsilon}\left(  T\right)  \rightarrow
\nu_{M-\varepsilon}$ as $T\rightarrow\infty,$ $\left(  \ref{167}\right)  $
implies that $N\left(  \mu_{\mathbf{0},M-\varepsilon}\left(  T\right)
\right)  \rightarrow N\left(  \nu_{M-\varepsilon}\right)  ,$ which proves
$\left(  \ref{166}\right)  .$

In the above proof the important step was to replace the state $\mu
_{\nu,\lambda_{\nu}\left(  \cdot\right)  }\left(  K_{2}\right)  $ with a lower
state $\mathbf{0,}$ which is in fact the lowest. Turning to $\left(
\ref{165}\right)  ,$ we see that this step can not be mimicked there, since
there is no highest state! So, to proceed, we need some apriori upper bound on
the state $\mu_{\nu,\lambda_{\nu}\left(  \cdot\right)  }\left(  K_{1}\right)
.$

\begin{lemma}
Let $\nu$ be an arbitrary initial state, with $N\left(  \nu\right)  <\infty. $
Then there exist $\bar{c}\left(  \nu\right)  <1$ and $\mathcal{T}<\infty, $
such that for every $t>\mathcal{T}$%
\[
\lambda_{\nu}\left(  t\right)  <\bar{c}\left(  \nu\right)  .
\]

\end{lemma}

\begin{proof}
The statement of the lemma is equivalent to the fact that $M=\limsup
_{t\rightarrow\infty}\lambda_{\nu}\left(  t\right)  <1.$ So suppose the
opposite, that $M\geq1.$ As we then know from Lemma \ref{T}, for every
$\mathcal{T}$ and every $\varepsilon>0$ we can find a segment $\left[
K,K+\mathcal{T}\right]  ,$ such that $\lambda_{\nu}\left(  t\right)
>1-\varepsilon$ for all $t\in\left[  K,K+\mathcal{T}\right]  .$ This, however,
contradicts to the statement $\left(  \ref{01}\right)  $ of Lemma \ref{la}.
\end{proof}

Since $S\left(  \nu\right)  <\infty,$ and the rate $\lambda_{\nu}\left(
t\right)  $ is uniformly bounded, the function $S\left(  \mu_{\nu,\lambda
_{\nu}}\left(  t\right)  \right)  $ is finite for every $t.$ It can grow, but
after the moment $\mathcal{T},$ obtained in the last Lemma, it stays bounded
from above by the value $S\left(  \mu_{\nu,\lambda_{\nu}}\left(
\mathcal{T}\right)  \right)  +S\left(  \nu_{\bar{c}}\right)  $, because of
Lemma \ref{domin}. So without loss of generality we can assume that the
initial state $\nu$ itself is such that $N\left(  \nu\right)  <\infty,$
$S\left(  \nu\right)  <\infty,$ while $\lambda_{\nu}\left(  t\right)  <\bar
{c}<1$ for all $t>0$ and so $S\left(  \mu_{\nu,\lambda_{\nu}}\left(  t\right)
\right)  \leq S\left(  \nu\right)  +S\left(  \nu_{\bar{c}}\right)  .$
Therefore the family of states $\left\{  \mu_{\nu,\lambda_{\nu}}\left(
t\right)  ,t\geq0\right\}  $ is weakly compact. Hence for every $\varepsilon
^{\prime} $ there exists $T$ such that for all $t$
\[
N\left(  \mu_{\mu\left[  _{\nu,\lambda_{\nu}\left(  \cdot\right)  }\left(
t\right)  \right]  ,m+\varepsilon}\left(  T\right)  \right)  \leq N\left(
\nu_{m+\varepsilon}\right)  +\varepsilon^{\prime}.
\]
This proves $\left(  \ref{165}\right)  .$

\section{Conclusions \label{con}}

In this paper we have proven the Poisson Hypothesis for the information
networks, for the case of the ''mean-field'' model of the network. We have
found the domain of its validity, and we will show in the forthcoming paper
\cite{RS} that beyond this domain Poisson Hypothesis is violated even in the
mean-field case, and the dependence of the initial condition does not vanish
with time. We strongly believe that the methods we have developed here -- in
particular, the self-averaging relation -- are relevant not only for
mean-field models, but also for more realistic ones.

The following problems appear as the natural continuation of our work.

\begin{itemize}
\item The study of Poisson Hypothesis for the case of the service times
forming an ergodic random process, rather than the sequence of i.i.d. random variables.

\item The study of PH for the case of customers of several identities, with
service times depending on their identity.

\item The study of PH for more general graphs.
\end{itemize}

We are going to work on these problems in the near future.

$\medskip$

\textbf{Acknowledgment. }\textit{We would like to thank our colleagues -- in
particular, Yu. Golubev, T. Liggett, O. Ogievetsky, G. Olshansky, S. Pirogov,
A. Vladimirov -- for valuable discussions and remarks, concerning this paper.
We are grateful to our referee for his remarks, which resulted in a
substantial improvement of the presentation of the paper. We are grateful also
to the Institute for Pure and Applied Mathematics (IPAM) at UCLA, for the
uplifting atmosphere and support during the Spring 2002 program on Large Scale
Communication Networks, where part of this work was done.}

\end{document}